\documentclass{amsart}

\usepackage{amsrefs,amsmath,amsfonts,amssymb,graphicx,latexsym,color}
\usepackage{hyperref}

\newtheorem{theorem}{Theorem}
\newtheorem{proposition}[theorem]{Proposition}

\newtheorem{lemma}[theorem]{Lemma}

\newcommand{\eqdef}{\overset{\mbox{\tiny{def}}}{=}}
\DeclareMathOperator{\esssup}{ess\,sup}

\theoremstyle{definition}

\numberwithin{equation}{section}
\numberwithin{theorem}{section}

\begin{document}
\title[Asymptotic Stability for Soft Potentials]{Asymptotic Stability of the Relativistic Boltzmann Equation for the Soft Potentials}
\author[Robert M. Strain]{Robert M. Strain}
\thanks{The authors research was partially supported by the NSF grant DMS-0901463.}

\address{University of Pennsylvania, Department of Mathematics, David Rittenhouse
Lab, 209 South 33rd Street, Philadelphia, PA 19104-6395, USA} 
\email{strain at math.upenn.edu}
\urladdr{http://www.math.upenn.edu/~strain/}

\begin{abstract}
In this paper it is  shown that unique solutions to the relativistic Boltzmann equation exist for all time and decay with any polynomial rate towards their steady state relativistic Maxwellian provided that the initial data starts out sufficiently close in $L^\infty_\ell$.  If the initial data  are continuous then so is the corresponding solution.  We work in the case of a spatially periodic box.  Conditions on the collision kernel are generic in the sense of \cite{MR933458}; this resolves the open question of global existence for the soft potentials.
\end{abstract}

\maketitle

\setcounter{tocdepth}{1}
\tableofcontents

\thispagestyle{empty}

\section{Introduction}

The relativistic Boltzmann equation is a fundamental model for fast moving particles; it can be written with appropriate initial conditions as
$$
p^\mu \partial_\mu F = \mathcal{C}(F,F).
$$
The  collision operator \cite{MR1026740,MR635279} is given by
\begin{equation}
\mathcal{C}(f,h)
=
 \frac {c}{2}\int_{\mathbb{R}^3} \frac{dq}{q_0}\int_{\mathbb{R}^3}\frac{dq^\prime}{q^\prime_0}\int_{\mathbb{R}^3}\frac{dp^\prime}{p^\prime_0}W(p, q | p^\prime, q^\prime) [f(p^{\prime })h(q^{\prime})-f(p)h(q)].
\notag
\end{equation}
The transition rate, $W(p, q | p^\prime, q^\prime)$, can be expressed as
\begin{equation*}
W(p, q | p^\prime, q^\prime) 
=
s\sigma(g, \theta) \delta^{(4)}(p^\mu+q^\mu-p^{\mu\prime}-q^{\mu\prime}),
\end{equation*}
where $\sigma(g, \theta)$ is the differential cross-section or scattering kernel; it measures the interactions between particles.   
The speed of light is the physical constant denoted
$
c.
$
Standard references in relativistic Kinetic theory include 
\cite{MR1898707,MR635279,MR1379589,MR0088362,nla.cat-vn2540748}.
The rest of the notation is defined in the sequel. 

\subsection{A Brief History of relativistic Kinetic Theory}

Early results include those on derivations  \cite{MR0004796}, local existence \cite{MR0213137}, and  linearized solutions \cite{MR933458,MR1031410}.

DiPerna-Lions renormalized weak solutions \cite{MR1014927} were shown to exist in 1992 by Dudy{\'n}ski and Ekiel-Je{\.z}ewska \cite{MR1151987}  globally in time for large data, using the causality of the relativistic Boltzmann equation 
\cite{MR818441,MR841735}.  See also \cite{MR1402446,MR1480243}
and
\cite{MR1714446,MR1676150}.
In particular \cite{MR1402446} 
proves the strong $L^1$ convergence to a relativistic Maxwellian, after taking a subsequence, for weak solutions with large initial data that is
not necessarily close to an equilibrium solution.
There are also results in the context of local 
\cite{MR2098116} 
and global \cite{strainNEWT} Newtonian limits, and near Vacuum results \cite{MR2217287,strainNEWT,MR2459827} and blow-up \cite{MR2102321} for the gain term only.
We also mention a study of the collision map
 and the pre-post collisional change of variables in \cite{MR1105532}.  
For a more discussions of historical results, we refer to \cite{strainNEWT}.

We review in more detail the results most closely related to those in this paper.  
In 1993, Glassey and Strauss \cite{MR1211782} proved for the first time global existence and uniqueness of smooth solutions which are initially close to a relativistic Maxwellian and in a periodic box.  They also established exponential convergence to Maxwellian. Their assumptions on the differential cross-section, $\sigma$, fell into the regime of ``hard potentials'' as discussed below.  
 In 1995, they extended that result to 
the whole space case \cite{MR1321370} where the convergence rate is polynomial.  More recent results with reduced restrictions 
on the 
cross-section were  proven in \cite{MR2249574}, using the energy method from \cite{MR1908664,MR2000470,MR2013332,MR2095473}; these results also apply to the hard potentials.  

For relativistic interactions--when particles are fast moving--an important physical regime is the ``soft potentials''; see \cite{DEnotMSI} for a physical discussion.  Despite their importance, prior to the results in this paper there were no existence results for the soft potentials in the context of strong nearby relativistic Maxwellian global solutions.
In 1988 a general physical assumption was given in 
\cite{MR933458}; see \eqref{hypSOFT} and \eqref{hypHARD}.  
In this paper we will prove global existence of unique $L^\infty$ near equilibrium  solutions 
to the relativistic Boltzmann equation  and rapid time decay  under the full physical assumption on the cross-section from \cite{MR933458}.  Our main focus is of course the soft potentials; and we do not require any angular cut-off.

\subsection{Notation}
Prior to discussing our main results, we will now define the notation of the problem carefully.
The momentum of a particle is denoted by $p^\mu$, $\mu = 0, 1, 2, 3$.  
Let the signature of the metric be $(- + + + )$.
Without loss of generality, we set the rest mass for each particle $m=1$.
The  momentum for each particle is restricted to the mass shell
$
p^\mu p_\mu = -c^2,
$
with
$
p^0>0.
$

Further with 
$
p \in \mathbb{R}^3,
$
we may write
$
p^\mu=(-p^0, p)
$  
and similarly
$
q^\mu=(-q^0, q).
$  
Then  the energy of a relativistic particle with momentum $p$ is $p^{0}=\sqrt{c^2+|p|^2}$.
We will use the Lorenz inner product which is given by
$$
p^\mu q_\mu \eqdef -p_0 q_0+\sum_{i=1}^3 p_i q_i.
$$
Now the streaming term of the relativistic Boltzmann equation is
$$
p^\mu \partial_\mu 
=
p_0\partial_t + p\cdot \nabla_x.
$$
We thus write the relativistic Boltzmann equation as
\begin{equation}
\partial _{t}F+\hat{p}\cdot \nabla_x F=\mathcal{Q}(F,F).
\label{RBF}
\end{equation}
Here $\mathcal{Q}(F,F) = \mathcal{C}(F,F)/p_0$, with $\mathcal{C}$
defined at the top of this paper.  

Above we consider $F=F(t,x,p)$ to be a  function of time $t\in [0,\infty)$, space $x\in \mathbb{T}^3$ and momentum $p\in \mathbb{R}^3$.  
   The normalized velocity of a  particle is denoted 
\begin{equation}
\hat{p}=c\frac{p}{p_0}=\frac{p}{\sqrt{1+|p|^2/c^2}}.
\label{normV}
\end{equation}
Steady states of this model are the well known J\"{u}ttner solutions, also known as the relativistic Maxwellian.
They are given by
$$
J (p)
\eqdef
\frac{\exp\left(-c p_0/(k_B T)\right)}{4\pi   c k_B T K_2(c^2/(k_B T))},
$$
where $K_2(\cdot)$ is the Bessel function 
$
K_2(z)\eqdef \frac {z^2}2 \int_1^\infty e^{-zt} (t^2-1)^{3/2} dt,
$ 
$T$ is the temperature and $k_B$ is Boltzmann's constant. 

In the rest of this paper, without loss of generality but for the sake of simplicity, we will now normalize all physical constants to one, including the
 speed of light  to $c = 1$.  So that in particular we denote the relativistic Maxwellian by
 \begin{equation}
J (p)
=
\frac{e^{- p_0}}{4\pi}.
\label{juttner}
\end{equation}
Henceforth we let $C$, and sometimes $c$ denote generic positive inessential constants whose value may change from line to line.
 
We will now define the quantity $s$,  which is 
 the square of the energy in the ``center of momentum''  system, $p+q=0$, as
\begin{eqnarray}
s
\eqdef
-(p^\mu+q^\mu) (p_\mu+q_\mu)
=
2\left( -p^\mu q_\mu+1\right)\ge 0.
\label{sDEFINITION}
\end{eqnarray}
The relative momentum is denoted
\begin{gather}
g 
\eqdef
\sqrt{(p^\mu-q^\mu) (p_\mu-q_\mu)}
=
\sqrt{2(-p^\mu q_\mu-1)}.
\label{gDEFINITION}
\end{gather}
Notice that $s=g^2+4$.  We warn the reader that this notation, which is used in \cite{MR635279}, may differ from other authors notation by a constant factor.

Conservation of momentum and energy for elastic collisions is expressed as
\begin{equation}
p^\mu+q^\mu=p^{\mu\prime}+q^{\mu\prime}.
\label{collisionalCONSERVATION}
\end{equation}
The angle $\theta$ is defined by
\begin{equation}
\cos\theta
\eqdef
(p^\mu - q^\mu) (p_\mu^\prime -q_\mu^\prime)/g^2.
\label{angle}
\end{equation}
This angle is well defined under \eqref{collisionalCONSERVATION}, see  \cite[Lemma 3.15.3]{MR1379589}.

We now consider the center of momentum expression for the collision operator below.  An alternate expression for the collision operator was derived in \cite{MR1211782}; see \cite{strainNEWT} for an explanation of the connection between the expression from \cite{MR1211782} and the one we give just now.
One may use Lorentz transformations
as described in \cite{MR635279} and \cite{strainPHD} to reduce the delta functions and obtain 
\begin{equation}
\mathcal{Q}(f,h)
=
\int_{\mathbb{R}^3}  dq ~
\int_{\mathbb{S}^{2}} d\omega ~
~ v_{\o} ~ \sigma (g,\theta ) ~ [f(p^{\prime })h(q^{\prime})-f(p)h(q)].
\label{collisionCM}
\end{equation}
where $v_{\o}=v_{\o}(p,q)$ is the M\o ller velocity given by
\begin{equation}
v_{\o}=
v_{\o}(p,q)
\eqdef
\sqrt{\left| \frac{p}{p_0}-\frac{q}{q_0}\right|^2-\left| \frac{p}{p_0}\times\frac{q}{q_0}\right|^2}
=
\frac{ g\sqrt{s}}{p_0 q_0}.
\label{moller}
\end{equation}
The post-collisional momentum in the expression (\ref{collisionCM}) can be written:
\begin{equation}
\begin{split}
p^\prime&=\frac{p+q}{2}+\frac{g}{2}\left(\omega+(\gamma-1)(p+q)\frac{(p+q)\cdot \omega}{|p+q|^2}\right), 
\\
q^\prime&=\frac{p+q}{2}-\frac{g}{2}\left(\omega+(\gamma-1)(p+q)\frac{(p+q)\cdot \omega}{|p+q|^2}\right),
\label{postCOLLvelCMsec2}
\end{split}
\end{equation}
where $\gamma =(p_0+q_0)/\sqrt{s}$.
  The energies are then
\begin{equation}
\begin{split}
p^{0\prime}&=\frac{p^0+q^0}{2}+\frac{g}{2\sqrt{s}}\omega\cdot (p+q), 
\\
q^{0\prime}&=\frac{p^0+q^0}{2}-\frac{g}{2\sqrt{s}}\omega\cdot (p+q).
\end{split}
\label{0prime}
\end{equation}
These clearly satisfy \eqref{collisionalCONSERVATION}.  The angle further satisfies
$
\cos\theta = k\cdot \omega
$
with $k = k(p,q)$ and $| k| = 1$.    
This is the expression for the collision operator that we will use.

For a smooth function $h(p)$ 
the collision operator satisfies
$$
\int_{\mathbb{R}^3}
dp
\begin{pmatrix}
      1   \\      p  \\ p_0
\end{pmatrix}\mathcal{Q}(h, h)(p)
=
0.
$$
By integrating the relativistic Boltzmann equation \eqref{RBF} and using this identity we obtain the conservation of mass, momentum and  energy for solutions as
$$
\frac d{dt}
\int_{\mathbb{T}^3}  dx ~
\int_{\mathbb{R}^3}   dp ~
\begin{pmatrix}
      1   \\      p  \\ p_0
\end{pmatrix}
 F(t)
=0.
$$
Furthermore the entropy of the relativistic Boltzmann equation is defined as
\begin{equation}
\mathcal{H}(t)
\eqdef
-\int_{\mathbb{T}^3}  dx ~
\int_{\mathbb{R}^3}   dp ~
 F(t, x, p)\ln F(t, x, p).  
\nonumber
\end{equation}
The celebrated Boltzmann
H-Theorem  is then formally
\begin{equation}
\frac d{dt}\mathcal{H}(t)
 \ge 0,  \nonumber
\end{equation}
which says that the entropy of solutions is increasing as time passes.   Notice that the steady state relativistic Maxwellians \eqref{juttner} maximize the entropy which formally implies convergence to \eqref{juttner} in large time.
It is this formal reasoning that our main results make mathematically rigorous in the context of perturbations of the relativistic Maxwellian for a general class of cross-sections.

\section{Statement of the Main Results}

We are now ready to discuss in detail our main results.  
We define the standard perturbation $f(t,x,p)$ to the relativistic Maxwellian \eqref{juttner} as 
$$
F \eqdef J +\sqrt{J } f. 
$$
With \eqref{collisionalCONSERVATION} we observe that the quadratic collision operator \eqref{collisionCM} satisfies
$$
\mathcal{Q}(J, J)
=0.
$$
 Then the relativistic Boltzmann equation \eqref{RBF} for the perturbation $f=f(t,x,p)$ takes the form
\begin{gather}
 \partial_t f + \hat{p}\cdot \nabla_x f + L (f)
=
\Gamma (f,f),
\quad
f(0, x, p)=f_0(x,p).
\label{rBoltz}
\end{gather}
 The linear operator $L( f)$, as defined in (\ref{L}), and the non-linear operator $\Gamma(f,f)$, defined  (\ref{gamma0}), are derived from an expansion of the relativistic Boltzmann collision operator \eqref{collisionCM}.
In particular, the linearized collision operator is given by 
\begin{gather}
 L(h)
 \eqdef 
 -J^{-1/2}\mathcal{Q}(J ,\sqrt{J} h)- J^{-1/2}\mathcal{Q}(\sqrt{J} h,J)
 \label{L}
 \\
  =
 \nu(p) h-K (h).
 \notag
\end{gather}
Above the multiplication operator takes the form
\begin{equation}
\nu(p) \eqdef
\int_{\mathbb{R}^3} ~  dq 
\int_{\mathbb{S}^{2}} ~ d\omega
~ v_{\o} ~ \sigma(g,\theta)~ J(q).
\label{nuDEF}
\end{equation}
The remaining integral operator is
\begin{gather}
K(h) 
\eqdef
\int_{\mathbb{R}^3} ~  dq 
\int_{\mathbb{S}^{2}} ~ d\omega
~ v_{\o} ~ \sigma(g,\theta)
\sqrt{J(q)}\left\{ \sqrt{J(q^{\prime})} ~ h(p^{\prime })+\sqrt{J(p^{\prime})} ~ h(q^{\prime})\right\}
\notag
\\
-\int_{\mathbb{R}^3} ~  dq 
\int_{\mathbb{S}^{2}} ~ d\omega
~ v_{\o} ~ \sigma(g,\theta)
~ \sqrt{J(q) J(p)} ~ h(q)
\label{compactK}
\\
=
K_2(h) - K_1(h).
\notag
\end{gather}
The non-linear part of the collision operator is defined as
\begin{multline}
\Gamma (h_1,h_2)
\eqdef
 J^{-1/2}\mathcal{Q}(\sqrt{J} h_1,\sqrt{J} h_2)
\label{gamma0}
\\
=
\int_{\mathbb{R}^3} ~  dq 
\int_{\mathbb{S}^{2}} ~ d\omega
~ v_{\o} ~ \sigma(g,\theta)~
 \sqrt{J(q)}
 [h_1(p^{\prime })h_2(q^{\prime})-h_1(p)h_2(q)].
\end{multline}
Without loss of generality we can assume that the mass, momentum, and energy conservation laws for the perturbation, $f(t,x,p)$, hold for all $t\ge0$ as
\begin{equation}
\int_{\mathbb{T}^3}  dx ~
\int_{\mathbb{R}^3}   dp ~
\begin{pmatrix}
      1   \\      p  \\ p_0
\end{pmatrix}
\sqrt{J(p)} ~ f(t,x,p) 
=
0.
\label{conservation}
\end{equation}
We now state our conditions on the collisional cross-section. 

\subsection*{\bf Hypothesis on the collision kernel:}   
{\it For soft potentials we assume the collision kernel in \eqref{collisionCM} satisfies 
the following growth/decay estimates 
\begin{equation}
\begin{split}
\sigma (g,\omega) 
 &
\lesssim  g^{-b} ~ \sigma_0(\omega),
\\ 
\sigma (g,\omega) 
 &
\gtrsim
\left( \frac{g}{\sqrt{s}}\right) g^{-b} ~ \sigma_0(\omega).
\end{split}
\label{hypSOFT}
\end{equation}
We also consider angular factors such that
$\sigma_0(\omega) \lesssim \sin^\gamma\theta$
with $\gamma> -2$.  Additionally $\sigma_0(\omega) \ge 0$ and $\sigma_0(\omega)$ should be non-zero on a set of positive measure.
We suppose further that $0 < b <\min(4,4+\gamma)$.

For hard potentials we make the assumption
\begin{equation}
\begin{split}
\sigma (g,\omega) 
 &
 \lesssim
  \left(  g^{a}+ g^{-b}\right) ~ \sigma_0(\omega) ,
\\ 
\sigma (g,\omega) 
 &
 \gtrsim
  \left( \frac{g}{\sqrt{s}}\right) g^{a}  ~ \sigma_0(\omega).
\end{split}
\label{hypHARD}
\end{equation}
In addition to the previous parameter ranges we consider $0\le a\le 2+\gamma$ and also $0 \le b <\min(4,4+\gamma)$ (in this case we allow the possibility of $b=0$).
}

\bigskip

This hypothesis includes the full range of conditions which were introduced in \cite{MR933458} as a general physical assumption on the kernel (of course we add the corresponding necessary lower bound in each case); see also \cite{DEnotMSI} for further discussions.

Prior to stating our main theorem, we will need to introduce the following
 mostly standard notation.
The notation $A \lesssim B$ will imply that a positive constant $C$ exists such that $A \leq C B$ holds uniformly over the range of parameters which are present in the inequality and  moreover that the precise magnitude of the constant is unimportant.  The notation $B \gtrsim A$ is equivalent to $A \lesssim B$, and $A \approx B$ means that both $A \lesssim B$ and $B \lesssim A$.   
 We work with the $L^\infty$ norm
\begin{equation*}
\| h\|_\infty
\eqdef
{\rm ess~sup}_{x\in \mathbb{T}^3, ~ p\in\mathbb{R}^3}  |h(x,p)|.
\end{equation*}
If we only wish to take the  supremum in the momentum variables we write
\begin{equation*}
| h|_\infty
\eqdef
{\rm ess~sup}_{p\in\mathbb{R}^3}  |h(p)|.
\end{equation*}
We will additionally use the following  standard $L^2$ spaces
\begin{equation*}
\| h\|_2
\eqdef
\sqrt{\int_{\mathbb{T}^3 } ~ dx ~ 
\int_{\mathbb{R}^3} ~ dp ~  
|h(x,p)|^2},
\quad
| h|_2
\eqdef
\sqrt{\int_{\mathbb{R}^3} ~ dp ~  
|h(p)|^2}.
\end{equation*}
Similarly in the sequel any norm represented by one set of lines instead of two only takes into account the momentum variables.
Next we define the norm which measures the (very weak) ``dissipation'' of the linear operator
\begin{equation*}
\| h\|_\nu
\eqdef
\sqrt{
\int_{\mathbb{T}^3 } ~ dx ~ 
\int_{\mathbb{R}^3} ~ dp ~  
 \nu(p) |h(x,p)|^2}.
\end{equation*}
The $L^2(\mathbb{R}^n)$ inner product is denoted $\langle \cdot, \cdot \rangle$.   We use $(\cdot, \cdot )$ to denote the  $L^2(\mathbb{T}^n \times \mathbb{R}^n)$ inner product.
Now, for $\ell\in\mathbb{R}$, we define the following weight function
\begin{equation}
w_{\ell}= w_\ell(p)
\eqdef 
\left\{
\begin{array}{cc}
p_0^{\ell b/2}, &  \text{for the soft potentials:} \quad \eqref{hypSOFT}\\
p_0^{\ell},  & \text{for the hard potentials:} \quad \eqref{hypHARD}.\\
\end{array}
\right.
\label{weight}
\end{equation}
For the soft potentials $w_1(p) \approx 1/\nu(p)$ (Lemma \ref{nuEST}).
We consider weighted spaces
\begin{gather*}
\| h\|_{\infty,\ell}
\eqdef
\| w_{\ell} h\|_{\infty},
\quad
\|  h \|_{2,\ell}
\eqdef
\| w_{\ell} h \|_{2},
\quad
\| h\|_{\nu,\ell}
\eqdef
\| w_{\ell} h\|_{\nu}.
\end{gather*}
Here as usual $L^\infty_{\ell}(\mathbb{T}^3 \times \mathbb{R}^3)$ is the Banach space with norm $\| \cdot \|_{\infty,\ell}$ etc.
We will also use the momentum  only counterparts of these spaces
\begin{gather*}
| h |_{\infty,\ell}
\eqdef
| w_{\ell} h|_{\infty},
\quad
|  h |_{2,\ell}
\eqdef
| w_{\ell} h |_{2},
\quad
| h|_{\nu,\ell}
\eqdef
| w_{\ell} h |_{\nu}.
\end{gather*}
We further need the following time decay norm
\begin{equation}
||| f |||_{k,\ell} \eqdef \sup_{s\ge 0} ~ (1+s)^k ~ \| f(s) \|_{\infty, \ell}.
\label{dNORM}
\end{equation}
We are now ready to state our main results.
We will first state our theorem for the soft potentials which is the main focus of this paper:

\begin{theorem}
\label{mainGLOBAL}  
(Soft Potential)
Fix $\ell >3/b$.
Given $f_0 = f_0(x,p) \in L^\infty_{\ell}(\mathbb{T}^3 \times \mathbb{R}^3)$ which satisfies \eqref{conservation} initially.  There is an $\eta>0$ such that if $\|f_0\|_{\infty, \ell} \le \eta$, then   there exists a unique global mild solution, $f = f(t,x,p)$, to equation \eqref{rBoltz} with soft potential kernel \eqref{hypSOFT}.
For any  $k\ge 0$, there is a $k'=k'(k)\ge 0$ such that 
$$
\| f \|_{\infty, \ell}(t) \le 
C_{\ell, k} (1+t)^{-k} \| f_0 \|_{\infty, \ell+k'}.
$$
These solutions are continuous if it is so initially.  We further have positivity, i.e. $F= \mu + \sqrt{\mu} f \ge 0$ if $F_0= \mu + \sqrt{\mu} f_0 \ge 0$.
\end{theorem}

We point out that $k'(0)=0$ in the above theorem, and $k'(k)\ge k$ can in general be computed explicitly from our proof.
Our approach also applies to the hard potentials and in that case we state the following theorem which can be proven using the same methods.

\begin{theorem}
\label{mainHARD}  
(Hard Potential)  
Fix $\ell >3/2$.
Given $f_0 = f_0(x,p) \in L^\infty_{\ell}(\mathbb{T}^3 \times \mathbb{R}^3)$ which satisfies \eqref{conservation} initially.  There is an $\eta>0$ such that if $\|f_0\|_{\infty, \ell} \le \eta$, then   there exists a unique global mild solution, $f = f(t,x,p)$, to equation \eqref{rBoltz} 
with hard potential kernel \eqref{hypHARD}
 which further satisfies for some $\lambda >0$ that
$$
\|f\|_{\infty, \ell}(t) \le C_{\ell} e^{-\lambda t} \|f_0\|_{\infty, \ell}.
$$
These solutions are continuous if it is so initially.  We further have positivity, i.e. $F= \mu + \sqrt{\mu} f \ge 0$ if $F_0= \mu + \sqrt{\mu} f_0 \ge 0$.
\end{theorem}

Previous results for the hard potentials are as follows.  In 1993 
Glassey and Strauss \cite{MR1211782} proved asymptotic stability such as  Theorem \ref{mainHARD}  in  $L^\infty_\ell$ with $\ell>3/2$.  They
consider collisional cross-sections which satisfy \eqref{hypHARD} for the parameters
$b\in [0,1/2)$, $a\in [0, 2-2b)$ and either $\gamma \ge 0$ or 
$$
\left| \gamma \right| 
< 
\min\left\{2-a, \frac{1}{2} -b, \frac{2-2b - a}{3}
\right\},
$$
which in particular implies $\gamma > -\frac{1}{2}$ if $b=0$ say.
They further assume a related growth bound on the derivative of the  cross-section
$
\left| \frac{\partial \sigma}{\partial g} \right|.
$
In \cite{MR2249574} this growth bound was removed while the rest of the assumptions on the cross-section from \cite{MR1211782} remained the same.   These results also sometimes work in smoother function spaces, and we note that we could also include space-time regularity to our solutions spaces.  

However for the relativistic Boltzmann equation
in \eqref{RBF} 
 the issue of adding momentum derivatives is more challenging.  
In recent years many new tools have been developed to solve these problems.  
A method was developed in non-relativistic kinetic theory to study the soft potential Boltzmann equation with angular cut-off by Guo in \cite{MR2013332}.  This approach makes crucial use of the momentum derivatives, and Sobolev embeddings to control the singular kernel of the collision operator.  Yet in the context of relativistic interactions, high derivatives of the post-collisional variables \eqref{postCOLLvelCMsec2} create additional high singularities which are hard to control.  Worse in the more common relativistic variables from \cite{MR1211782},  derivatives of the post-collisional momentum exhibit enough momentum growth to preclude hope of applying the method from \cite{MR2013332}; these growth estimates on the derivatives were known  much earlier in \cite{MR1105532}.

Notice also that the methods for proving time decay, such as \cite{MR2209761,MR2366140,MR2116276},  require working in the context of  smooth solutions.  We would also like to mention recent developments on Landau Damping
\cite{Mouhot:1173020} proving exponential decay with analytic regularity.  Furthermore we point out very recent results proving rapid time decay of smooth perturbative solutions to the Newtonian Boltzmann equation without the Grad angular cut-off assumption as in \cite{gsNonCut1,gsNonCut2,gsNonCutA}.
In this paper however we avoid smooth function spaces in particular because of the aforementioned problem created by the relativistic post-collisional momentum.
Other recent work \cite{G2} developed a framework to study near Maxwellian boundary value problems for the hard potential Newtonian Boltzmann equation in $L^\infty_\ell$.  In particular a key component of this analysis was to 
consider solutions to the Boltzmann equation \eqref{rBoltz} after linearization as
\begin{gather}
 \left( \partial_t  + \hat{p}\cdot \nabla_x  + L \right) f
=
0,
\quad
f(0,x,p)=f_{0}(x,p),
\label{rBlin}
\end{gather}
with $L$ defined in \eqref{L}.  The semi-group for this equation (relativistic or not) will satisfy a certain `A-Smoothing property' which was pioneered by Vidav \cite{MR0259662} about 40 years ago.  This A-Smoothing property which appears at the level of the second iterate of the semi-group, as seen below in \eqref{linearEXP}, has been  employed effectively  for instance in \cite{MR1211782,MR1379589,G2,MR1361017}. 
Related to this a  new compactness connected to a similar iteration has been observed in the `mixing lemma' of \cite{MR2284213,MR2230121,MR2082240}.  The key new step in \cite{G2} was to estimate the second iterate in $L^2$ rather than $L^\infty$ and then use a linear decay theory in $L^2$ which does not require regularity and is exponential for the hard potentials case that paper considered.    Note further that a method was developed in \cite{MR2209761,MR2366140} to prove rapid polynomial decay for the soft potential Newtonian Boltzmann and Landau equations; this is related to the articles \cite{MR2116276,MR576265,MR575897,MR677262}, all of which make use of smooth function spaces.  

In the present work  we adapt the method from \cite{MR2366140} to prove rapid $L^2$ polynomial decay of solutions to the linear equation \eqref{rBlin} without regularity, and we further adapt the $L^2$ estimate from \cite{G2} to control the second iterate.  This approach, for the soft potentials in particular, yields global bounds and slow decay, $O(1/t)$, of solutions to \eqref{rBlin} in $L^\infty_\ell$.  The details and complexity of this program are however intricate in the relativistic setting.  And fortunately, this slow decay is sufficient to just barely control the nonlinearity and prove global existence to the full non-linear problem as in Theorem \ref{mainGLOBAL}.

It is not clear how to apply the above methods to establish the rapid ``almost exponential'' polynomial decay from Theorem \ref{mainGLOBAL} in this low regularity $L^\infty_\ell$ framework.  
To prove the rapid decay in Theorem \ref{mainGLOBAL}  our key contribution is to perform a new high order expansion of the remainder term, $R_1(f)$, from \eqref{r1def}.  This is contained in Proposition \ref{r1rapidDecay} and its proof.  This term, $R_1(f)$, is the crucial problematic term which only appears to exhibit first order decay.

More generally, the main difficulty with proving rapid decay for the soft potentials is created by the high momentum values, where the time decay is diluted by the momentum decay.  This results in the generation of additional weights, typically one weight for each order of time decay.  At the same time the term $R_1(f)$ only allows us to absorb one weight, $w_1(p)$, and therefore only appears to allow one order of time decay.  In our proof of Proposition \ref{r1rapidDecay} we are able to overcome this apparent obstruction by performing a new high order expansion for $k\ge 2$ as
\begin{gather*}
R_k(f)
=
G_k(f)  + D_k(f) + N_k(f)  + L_k(f) + R_{k+1}(f).
\end{gather*}
(The expansion from $R_1(f)$ to $R_2(f)$ requires a slightly different approach.)  
At every level of this expansion we can peel of each of the terms $G_k(f)$, $D_k(f)$, $N_k(f)$,  and $L_k(f)$  which will for distinct reasons exhibit rapid polynomial decay to any order.  In particular we use an $L^2_\ell$ estimate for $L_k(f)$ which crucially makes use of the bounded velocities that come with special relativity.
On the other hand the last term $R_{k+1}(f)$ will be able to absorb $k+1$ momentum weights, and therefore it will be able to produce time decay up to the order $k+1$.  By continuing this expansion to any finite order, we are able to prove rapid polynomial decay.  We hope that this expansion will be useful in other relativistic contexts.

The rest of this paper is organized as follows.  In Section \ref{secL2} we prove $L^2_\ell$ decay of solutions to the linearized relativistic Boltzmann equation \eqref{rBlin}.
Then in Section \ref{secL0} we prove global $L^\infty_\ell$ bounds and slow decay of solutions to \eqref{rBlin}.
Following that in Section \ref{s:n0bd} we prove nonlinear $L^\infty_\ell$ bounds using the slow linear decay, and we thereby conclude global existence.
In the remaining Sections \ref{s:linRdecay}
and \ref{s:rapidD} we prove linear and non-linear rapid decay respectively.  
Then in the Appendix we give an exposition of a derivation of the kernel of the compact part of the linear operator from \eqref{compactK}.

\section{Linear $L^2$ Bounds and Decay}\label{secL2}

It is our purpose in this section to prove global in time $L^2_\ell$ bounds for solutions to 
the linearized Boltzmann equation \eqref{rBlin} with initial data  in the same space $L^2_\ell$.  We will then prove high order, almost exponential decay for these solutions.
We begin by stating a few important lemmas, and then we use them to prove the desired integral bounds \eqref{energyWest} and the decay it implies in Theorem \ref{decay2}.  We will prove these lemmas at the end of the section.

\begin{lemma}
\label{nuEST}  
Consider  \eqref{nuDEF} with the soft potential collision kernel \eqref{hypSOFT}.  Then
$$
\nu(p) 
\approx
 p_0^{-b/2}.
$$
More generally, 
$
\int_{\mathbb{R}^3} ~  dq 
\int_{\mathbb{S}^{2}} ~ d\omega
~ v_{\o} ~ \sigma(g,\theta)~ J^\alpha (q)
\approx
 p_0^{-b/2}
$
for any $\alpha >0$.
\end{lemma}

We will next look at the ``compact'' part of the linear operator $K$.  The most difficult part is 
 $K_2$ from \eqref{compactK}.  We will employ a splitting to cut out the singularity.  The new element of this splitting is the Lorentz invariant argument: $g$.  
Given a small $\epsilon>0$, choose a smooth cut-off function $\chi =\chi (g)$ satisfying
\begin{equation}
\chi (g)=
\left\{
\begin{array}{cl}
 1 & {\rm if } ~~  g \ge 2\epsilon
 \\
  0 &
  {\rm if } ~~ g \le \epsilon.
\end{array}
\right.
\label{cut}
\end{equation}
Now with \eqref{cut} and \eqref{compactK} we define
\begin{gather}
K_2^{1-\chi}(h) 
\eqdef
\int_{\mathbb{R}^3} ~  dq 
\int_{\mathbb{S}^{2}} ~ d\omega
~ \left( 1 - \chi(g) \right)
~ v_{\o} ~ \sigma(g,\theta)
\sqrt{J(q)} \sqrt{J(q^{\prime})} ~ h(p^{\prime })
\notag
\\
+
\int_{\mathbb{R}^3} ~  dq 
\int_{\mathbb{S}^{2}} ~ d\omega
~ \left( 1 - \chi(g) \right)
~ v_{\o} ~ \sigma(g,\theta)
\sqrt{J(q)}\sqrt{J(p^{\prime})} ~ h(q^{\prime}).
\label{kCUT}
\end{gather}
Define $K_1^{1-\chi}(h)$ similarly. 
We use the splitting
$
K
\eqdef
K^{1-\chi}
+
K^{\chi}.
$
A splitting with the same goals has been previously used for the Newtonian Boltzmann equation in \cite{MR2366140}.  
The advantage for soft potentials, on the singular region, is that one has exponential decay in all momentum variables.   Then on the region away from the singularity we are able to extract a modicum of extra decay which is sufficient for the rest of the estimates in this paper, see Lemma \ref{boundK2}   just below.

In the sequel we will use the Hilbert-Schmidt form for the non-singular part.  The following representation is derived in the Appendix:
$$
K_i^\chi (h) = 
\int_{\mathbb{R}^3}dq~ k_i^\chi(p,q) ~ h(q),
\quad
i=1,2.
$$
We will also record the kernel $k_i^\chi(p,q)$ below in  \eqref{k1def} and \eqref{k2def}.    We have

\begin{lemma}
\label{boundK2}  
Consider the soft potentials \eqref{hypSOFT}.  The kernel enjoys the estimate
$$
0\le k_2^{\chi}(p,q)
\le C_\chi \left( p_0 q_0 \right)^{-\zeta}
\left( p_0+ q_0 \right)^{-b/2}
e^{-c |p-q|},
\quad
C_\chi, c >0,
$$
with
$
\zeta
\eqdef
\min\left\{2-|\gamma|, 4-b,2\right\}/4
>0.
$
This estimate also holds for $k_1^\chi$.
\end{lemma}

We remark that for certain ranges of the parameters $\gamma$ and $b$ the decay in Lemma \ref{boundK2} could be improved somewhat.
In particular the term $k_1^\chi$ in \eqref{k1def} clearly yields exponential decay.
However what is written above is sufficient for our purposes.

We will use the decomposition given above and the decay of the kernels in Lemma \ref{boundK2}  
to establish the following lemma.

\begin{lemma}
\label{noKestimate}  
Fix any small $\eta >0,$ we may decompose  $K$ from \eqref{compactK} as 
$$
K = K_c + K_s,
$$
where $K_c$ is a compact operator on $L^2_\nu$.  In particular
for any $\ell\ge 0$, and for some $R=R(\eta)>0$ sufficiently large we have
$$
|\langle w^2_\ell  K_c h_1, h_2\rangle |
\le
C_{\eta} |{\bf 1}_{\le R} h_1 |_2 | {\bf 1}_{\le R} h_2 |_2.
$$
Above ${\bf 1}_{\le R}$ is the indicator function of the ball of radius $R$.
Furthermore
\begin{equation}
|\langle w^2_{\ell}  K_s h_1, h_2\rangle |
\le 
 \eta 
\left|  h_1\right| _{\nu,\ell }
 \left|   h_2\right| _{\nu,\ell  }.
\notag
\end{equation}
\end{lemma}

This estimate will be important for proving the coercivity of the linearized collision operator, $L$, away from its null space.  
More generally, from the H-theorem  $L$ is non-negative and for every fixed $(t,x)$ the null
space of $L$ is given by the five dimensional space \cite{MR1379589}:
\begin{equation}
{\mathcal N}\eqdef {\rm span}\left\{
\sqrt{J},
p_1 \sqrt{J}, p_2 \sqrt{J}, p_3 \sqrt{J}, 
p_0\sqrt{J}\right\}.
 \label{null}
\end{equation}
We define the orthogonal projection from $L^2(\mathbb{R}^3)$ onto the null space ${\mathcal N}$ by ${\bf P}$. 
Further expand ${\bf P} h$ as a linear combination of the basis in (\ref{null}): 
\begin{equation}
{\bf P} h
\eqdef
 \left\{a^h(t,x)+\sum_{j=1}^3 b_j^h(t,x)p_j+c^h(t,x)p_0\right\}\sqrt{J}.
\label{hydro}
\end{equation}
We can then decompose  $f(t,x,p)$ as
\[
f={\bf P}f+\{{\bf I-P}\}f.
\]
With this decomposition we have

\begin{lemma}
\label{lowerN}  $L\ge 0$.
$Lh = 0$ if and only if $h = {\bf P} h$.  
And
$\exists \delta_0>0$ such that
\begin{equation*}
\langle {L} h, h \rangle
\ge
\delta_0 
| \{ {\bf I - P } \} h |_\nu^2.
\end{equation*}
\end{lemma}

This last statement on coercivity holds as an operator inequality at the functional level.  The following lemma is a well-known statement of the linearized H-Theorem for solutions to \eqref{rBlin} which was shown in the non-relativistic case in \cite{G2}.

\begin{lemma}
\label{lowerA}  
Given initial data $f_0 \in L^2_\ell\left(\mathbb{T}^3\times \mathbb{R}^3\right)$ for some $\ell\ge0$,
which satisfies \eqref{conservation} initially. 
Consider the corresponding solution, $f$, to \eqref{rBlin} in the sense of distributions.
Then there is a universal constant $\delta_\nu>0$ such that
$$
\int_0^1~ ds ~ \| \{ {\bf I - P } \} f \|_\nu^2(s)
\ge 
\delta_\nu
\int_0^1 ~ ds ~ \| { \bf  P  } f \|_\nu^2(s).
$$
\end{lemma}

We will give just one more operator level inequality. 

\begin{lemma}
\label{weight2}  
Given $\delta\in(0,1)$ and $\ell \ge 0$,  there are constants $C, R >0$ such that
$$
\langle w^{2}_{\ell} L h, h \rangle
\ge 
\delta | h|_{\nu,\ell}^2
- C | {\bf 1}_{\le R} h|^2_2.
$$
\end{lemma}

Notice that Lemma \ref{weight2} follows easily from Lemma \ref{noKestimate}.
With these results, we can prove the following energy inequality for any 
$\ell \ge 0$:
\begin{equation}
 \|  f \|^2_{2,\ell}(t) + \delta_\ell \int_0^t ~ ds ~ \| f\|^2_{\nu,\ell}(s) 
 \le 
 C_\ell \| f_0\|^2_{2,\ell},
 \quad
\exists \delta_\ell, C_\ell >0,
\label{energyWest}
\end{equation}
as long as 
$
\|  f_0 \|^2_{2,\ell}
$
is finite.    We will prove this first for $\ell =0$, and then for arbitrary $\ell >0$.  In the first case we multiply \eqref{rBlin} with $f$ and integrate to obtain
$$
 \|  f \|^2_{2}(t) + \int_0^t ~ ds ~ (Lf,f)
=
 \| f_0\|^2_{2}.
$$
First suppose that $t\in\{1,2,\ldots\}$.
By Lemma \ref{lowerN} 
we have
\begin{multline*}
 \int_0^t ~ ds ~ (Lf,f)
 =
 \sum_{j=0}^{t-1}
  \int_0^1 ~ ds ~ (Lf,f)(s+j)
  \ge
   \sum_{j=0}^{t-1}
\delta_0 \int_0^1 ~ ds ~ \| \{ {\bf I - P } \} f\|^2_\nu(s+j)
\\
=
\frac{\delta_0}{2}    \sum_{j=0}^{t-1}
\int_0^1 ~ ds ~ \| \{ {\bf I - P } \} f\|^2_\nu(s+j)
+
\frac{\delta_0}{2}    \sum_{j=0}^{t-1}
\int_0^1 ~ ds ~ \| \{ {\bf I - P } \} f\|^2_\nu(s+j).
\end{multline*}
Then by Lemma \ref{lowerA}  
the second term above satisfies the lower bound
\begin{gather*}
\frac{\delta_0}{2}    \sum_{j=0}^{t-1}
\int_0^1 ~ ds ~ \| \{ {\bf I - P } \} f\|^2_\nu(s+j)
  \ge
\frac{\delta_0\delta _{\nu}}{2}    \sum_{j=0}^{t-1}
\int_0^1 ~ ds ~ \|  {\bf P } f\|^2_\nu(s+j).
\end{gather*}
This follows   in particular because   $f_{j}(s,x,v) \eqdef f(s+j,x,p)$ satisfies the linearized Boltzmann
equation \eqref{rBlin} on the interval $0\leq s\leq 1$.  Collecting the previous two estimates yields
\begin{multline*}
 \int_0^t ~ ds ~ (Lf,f)
  \ge
  \frac{\delta_0}{2}    \sum_{j=0}^{t-1}
\int_0^1 ~ ds ~ \| \{ {\bf I - P } \} f\|^2_\nu(s+j)
\\
+
\frac{\delta_0\delta _{\nu}}{2}    \sum_{j=0}^{t-1}
\int_0^1 ~ ds ~ \|  {\bf P } f\|^2_\nu(s+j)
\\
\ge
 \tilde{\delta} \sum_{j=0}^{t-1}
\int_0^1 ~ ds ~ \|   f\|^2_\nu(s+j)
=
\tilde{\delta}\int_0^t ~ ds ~ \| f\|^2_\nu(s),
\end{multline*}
with
$
\tilde{\delta} 
=
\frac{1}{2}
\min\left\{
\frac{\delta _{0} \delta _{\nu}}{2},
\frac{\delta _{0} }{2}
\right\}>0.
$
Plugging this estimate into the last one establishes the energy inequality \eqref{energyWest} for $\ell = 0$
and
$
t \in \{1,2,\dots \}.
$
For an arbitrary $t>0$, we choose
$
m \in \{ 0,1,2,\dots \}
$
such that
 $m\leq t\leq m+1$. 
 We then split the time integral as 
 $[0,t]=[0,m]\cup \lbrack m,t].$
For the time interval $[m,t]$ we have
\begin{equation}
\|f(t)\|_{2}^{2}+\int_{m}^{t}~ ds ~ (Lf,f)
=
\|f(m)\|_{2}^{2}.  
\label{mtESt}
\end{equation}
Since $L\ge 0$ by Lemma \ref{lowerN}, we  see that
$
\|f(t)\|_{2}^{2}
\le
\|f(m)\|_{2}^{2}.  
$
We then have
\begin{gather}
\| f \|_{2}^{2}(t)
+
\int_{m}^{t}~ ds ~ (Lf,f)
+
\frac{\tilde{\delta} }{2}
 \int_{0}^{m} ~ ds ~
\|f\|_{\nu }^{2}(s)
\leq 
C_k \| f_0 \|_{2}^{2}.
\notag
\end{gather}
Furthermore with Lemma \ref{weight2} for $C_\delta>0$ independent of $t$ we obtain   
\begin{multline}
 \int_m^t ~ ds ~ (Lf,f)
\ge 
\delta \int_m^t ~ ds ~ \| f\|^2_{\nu}(s)
-
C_{\delta} \int_m^t ~ ds ~ \| {\bf 1}_{\le R} f\|^2_{2}(s)
\\
\ge 
\delta \int_m^t ~ ds ~ \| f\|^2_{\nu}(s)
-
C_{\delta} ~ \sup _{m\le s \le t} \| f\|^2_{2}(s).
\label{lowerMM}
\end{multline}
However we have already shown that 
$
\sup _{m\le s \le t} \| f\|^2_{2}(s)
\le 
\|f\|_{2}^{2}(m)
\le
C \| f_0 \|_{2}^{2}.
$
Collecting the last few estimates for any $t> 0$ we have \eqref{energyWest} for $\ell =0$.  

For $\ell >0$, we multiply the equation
\eqref{rBlin} with $w_{\ell}^2 f$ and integrate to obtain
$$
 \|  f \|^2_{2,\ell}(t) + \int_0^t ~ ds ~ (w_{\ell}^2 Lf,f)
=
 \| f_0\|^2_{2,\ell}.
$$
In this case
using
Lemma \ref{weight2}    
we have
$$
 \int_0^t ~ ds ~ (w_{\ell}^2 Lf,f)
\ge 
\delta \int_0^t ~ ds ~ \| f\|^2_{\nu,\ell}(s)
-
C_{\delta,\ell} \int_0^t ~ ds ~ \| f\|^2_{\nu}(s).
$$
Adding this estimate to the line above it, we obtain
$$
 \|  f \|^2_{2,\ell}(t) 
 +
\delta \int_0^t ~ ds ~ \| f\|^2_{\nu,\ell}(s)
\le
 \| f_0\|^2_{2,\ell}
 +
 C_{\delta,\ell} \int_0^t ~ ds ~ \| f\|^2_{\nu}(s).
$$
Now the just proven integral inequality \eqref{energyWest}  in the case $\ell = 0$ (as an upper bound for the integral on the right hand side) establishes the claimed energy inequality \eqref{energyWest} for any $\ell >0$.
In fact we prove a more general time decay version of this inequality in \eqref{nTdecay}.
These will imply  the following rapid decay theorem.

\begin{theorem}
\label{decay2}  
Consider a solution $f(t,x,p)$ to the linear
Boltzmann equation (\ref{rBlin})
with data
$
\| f_0\|_{2,\ell+k} < \infty
$
 for some $\ell, k \ge 0$.  
Then
$$
 \|  f \|_{2,\ell}(t)\le C_{\ell,k}(1+t)^{-k}
 \|  f_0 \|_{2,\ell+k}.
$$
This allows ``almost exponential'' polynomial decay of any order.
\end{theorem}

Theorem \ref{decay2} is the main result of this section.  We now proceed to prove each of Lemma \ref{nuEST}, Lemma \ref{boundK2}, Lemma \ref{noKestimate}, Lemma \ref{lowerN}, Lemma \ref{lowerA},
and then
Theorem \ref{decay2} in order.  These proofs will complete this section on linear decay.  \\

\noindent {\it Proof of Lemma \ref{nuEST}.}  We will use the soft potential hypothesis
for the collision kernel
\eqref{hypSOFT}
to estimate
\eqref{nuDEF}.  For $\alpha >0$ we more generally consider 
$$
\nu_\alpha(p) 
\eqdef 
\int_{\mathbb{R}^3} ~ dq 
~ v_{\o} ~ J^\alpha (q) ~ \int_0^\pi ~d\theta~ \sigma(g,\theta)~\sin\theta.
$$
Initially we record the following pointwise 
estimates
\begin{equation}
\frac{|p-q|}{\sqrt{p_0 q_0}}
\le
g
\le
\left\{
\begin{array}{c}
|p-q| \\
 2 \sqrt{p_0q_0},\\
\end{array}
\right.
 \label{gESTgs}
\end{equation}
see \cite[Lemma 3.1]{MR1211782}.  
However notice that in \cite{MR1211782}  their ``$g$'' is actually defined to be 2 times our ``$g$''.
With the 
M{\o}ller velocity 
\eqref{moller},
we thus also have
$$
s=4+g^2
\lesssim
p_0 q_0,
\quad 
v_{\o} \lesssim 1.
$$
For 
$b\in(1,4)$, these estimates including \eqref{gESTgs} yield
$$
v_{\o} ~ g^{-b} 
=
\frac{\sqrt{s}}{p_0 q_0} g^{1-b}
\lesssim 
\frac{\sqrt{s}}{p_0 q_0}
\frac{(p_0 q_0)^{(b-1)/2}}{|p-q|^{b-1}}
\lesssim
\frac{(p_0 q_0)^{(b-2)/2}}{|p-q|^{b-1}}.
$$
We thus obtain
\begin{multline*}
\nu_\alpha (p)
\lesssim
\int_{\mathbb{R}^3} ~  dq 
~ 
\frac{(p_0 q_0)^{(b-2)/2}}{|p-q|^{b-1}}
~ J^\alpha (q)
\int_0^\pi ~d\theta~\sin^{1+\gamma}\theta
\\
\lesssim
p_0^{(b-2)/2}
\int_{\mathbb{R}^3} ~  dq 
~ 
\frac{J^{\alpha /2}(q)}{|p-q|^{b-1}}
\int_0^\pi ~d\theta~\sin^{1+\gamma}\theta
\\
\lesssim
 p_0^{(b-2)/2}
p_0^{1-b}
\approx
 p_0^{-b/2}.
\end{multline*}
We note that the angular integral is finite since $\gamma > -2$:
$$
\int_0^\pi ~d\theta~\sin^{1+\gamma}\theta = C_\gamma <\infty.
$$
In this case above and the cases below a key point is that
$$
\int_{\mathbb{R}^3} ~  dq 
~ 
\frac{J^{\alpha /2}(q)}{|p-q|^{\beta}}\approx p_0^{-\beta},
\quad
\forall \beta < 3.
$$
For 
$b\in(0,1)$, with \eqref{gESTgs} we alternatively have
$$
v_{\o} ~ g^{-b} 
=
\frac{\sqrt{s}}{p_0 q_0} g^{1-b}
\lesssim
\frac{(p_0 q_0)^{(1-b)/2}}{\sqrt{p_0 q_0}}
\lesssim
(p_0 q_0)^{-b/2}.
$$
Now in a slightly easier way than for the previous case we have 
$
\nu_\alpha (p)
\lesssim
 p_0^{-b/2}.
$

For the lower bound with $b\in(2,4)$, we use \eqref{hypSOFT}, \eqref{gESTgs} and the estimate
$$
v_{\o} ~ \left( \frac{g}{\sqrt{s}}\right) g^{-b} 
=
\frac{g^{2-b}}{p_0 q_0} 
\gtrsim
 \frac{(p_0 q_0)^{(2-b)/2}}{p_0 q_0}
\approx
(p_0 q_0)^{-b/2}.
$$
Alternatively for $b\in (0,2)$,  we use \eqref{gESTgs}  to get the estimate
$$
v_{\o} ~ \left( \frac{g}{\sqrt{s}}\right) g^{-b} 
=
\frac{g^{2-b}}{p_0 q_0} 
\gtrsim
\frac{1}{p_0 q_0}
\left(
\frac{|p-q|}{\sqrt{p_0 q_0}}
\right)^{2-b}
\approx
(p_0 q_0)^{(b-4)/2} 
|p-q|^{2-b}.
$$
We use these estimates to obtain
$
\nu_\alpha (p)
\gtrsim
 p_0^{-b/2}
$
as  for the upper bound.
\qed \\

Now that the proof of Lemma \ref{nuEST} is complete, we develop the necessary formulation  for the proof of Lemma \ref{boundK2}.  It is trivial to write the Hilbert-Schmidt form for the cut-off \eqref{cut} part of $K_1$ from \eqref{compactK} as
\begin{equation}
k_1^\chi(p,q) = 
\sqrt{J(q)J(p)}
~\chi(g)
\int_{ \mathbb{S}^{2}}
~ d\omega
~ v_{\o} ~ \sigma(g,\theta).
\label{k1def}
\end{equation}
Furthermore, we can write the Hilbert-Schmidt form for the cut-off \eqref{cut} part of $K_2$ from \eqref{compactK} in the following somewhat complicated integral form
\begin{multline}
k_2^\chi(p,q)
\eqdef
c_2
\frac{s^{3/2} }{gp_0q_0} 
\chi(g) ~
\int_0^{\infty} ~dy~
e^{-l\sqrt{1+y^2}} ~
\sigma\left(\frac{g}{\sin \left(\psi/2 \right)},\psi \right)
\\
\quad 
\times
\frac{y\left(1+\sqrt{1+y^2}\right)}{\sqrt{1+y^2}}
I_0 (j y).
\label{k2def}
\end{multline}
Here
 $c_2>0$, and
 the modified Bessel function of index zero is defined by
\begin{equation}
I_0(y)=\frac{1}{2\pi}\int_0^{2\pi} ~  e^{y \cos\varphi} ~ d\varphi.
\label{bessel0}
\end{equation}
We are also using the simplifying notation
\begin{equation}
\sin \left(\psi/2 \right)
\eqdef
\frac{\sqrt{2} g}{[g^2 - 4 + (g^2 + 4)\sqrt{1+y^2}]^{1/2}}.
\label{sinPSI}
\end{equation}
Additionally
$$
l = (p_0 + q_0)/2,
\quad
j = \frac{|p\times q|}{g}.
$$ 
This derivation for a pure $k_2$ operator appears to go back to \cite{MR635279,MR0471665}, where it was done in the case of the alternate linearization $F = J(1+f)$. 
The author gave a similar derivation with many details explained in full in \cite{strainPHD}, including the explicit form of the necessary Lorentz transformation; in particular equation (5.51) in this thesis.  For the benefit of the reader, we have provided the derivation of $k_2^\chi$ in the case of the linearization $F = J+\sqrt{J}f$ in the appendix to this paper.

Note that it is elementary to verify that \eqref{k1def} under \eqref{hypSOFT} satisfies the
estimate in Lemma \ref{boundK2}.   In the proof below we focus on the more involved estimate for \eqref{k2def}. 
\\

\noindent
{\it Proof of Lemma \ref{boundK2}.}
We consider $K_2$ from \eqref{compactK} with Hilbert-Schmidt form represented by the kernel
\eqref{k2def}.  
From \eqref{hypSOFT}, we have the bound
\begin{multline}
\sigma\left(\frac{g}{\sin \left(\psi/2 \right)},\psi \right)
\lesssim
\left(\frac{\sin \left(\psi/2 \right)}{g}\right)^b
\sin^\gamma \psi
\\
\lesssim
 g^\gamma \left(\frac{\sin \left(\psi/2 \right)}{g}\right)^{b+\gamma}
\cos^\gamma \left(\psi /2 \right).
\label{SIGest}
\end{multline}
We have just used the trigonometric identity
$$
\sin\psi = 2\sin \left(\psi /2 \right) \cos \left(\psi /2 \right).
$$
We estimate these angles in three cases.  In each of the cases below we will repeatedly use the following known \cite[p.317]{MR1211782} estimates 
\begin{equation}
\frac{y}{2\sqrt{1+y^2}}
\le 
\cos \left(\psi /2 \right) 
\le 
1.
\label{cosEST}
\end{equation}
The estimates above and below are proved for instance in \cite[Lemma 3.1]{MR1211782}:
\begin{equation}
\frac{1}{\sqrt{s} ( 1+ y^2 )^{1/4} }
\lesssim
\frac{\sin \left(\psi /2 \right) }{g}
\lesssim
\frac{1}{g ( 1+ y^2 )^{1/4}}.
\label{sinEST}
\end{equation}
Notice that in general
from \eqref{k2def} and \eqref{SIGest} we have the bound
\begin{multline}
k_2^\chi(p,q)
\lesssim
\frac{s^{3/2} }{gp_0q_0} 
 g^\gamma
\chi(g)
\int_0^{\infty} ~dy~
e^{-l\sqrt{1+y^2}}
y
I_0 (j y)
\\
\quad \times
\left(\frac{\sin \left(\psi/2 \right)}{g}\right)^{b+\gamma}
\cos^\gamma \left(\psi /2 \right).
\label{FIRSTk2est}
\end{multline}
We will estimate this upper bound in three cases.

{\bf Case 1}.
Take
 $\gamma= -\left| \gamma \right| < 0$ and $b+\gamma= b-\left| \gamma \right| <0$.
Then 
we have
\begin{multline*}
\left(\frac{\sin \left(\psi/2 \right)}{g}\right)^{b+\gamma}
\cos^\gamma \left(\psi /2 \right)
\lesssim
s^{-(b -\left| \gamma \right|)/2}
\left(\frac{1}{ ( 1+ y^2 )^{1/4} } \right)^{b-\left| \gamma \right|}
\left(\frac{y}{\sqrt{1+y^2}}\right)^{-\left| \gamma \right|}
\\
\lesssim
 s^{-(b-\left| \gamma \right|)/2} y^{-\left| \gamma \right|}
{ ( 1+ y^2 )^{-(b-3\left| \gamma \right|)/4} }
\\
\lesssim
 s^{-(b -\left| \gamma \right|)/2} y^{-\left| \gamma \right|}
{ ( 1+ y^2 )^{-(b+3\gamma)/4} }.
\end{multline*}
Above we have used  \eqref{cosEST} and \eqref{sinEST}.
In this case from \eqref{FIRSTk2est} we have
\begin{gather*}
k_2^\chi(p,q)
\lesssim
\frac{s^{3/2} }{gp_0q_0} 
 \frac{s^{-(b -\left| \gamma \right|)/2}}{g^{\left| \gamma \right|}}
\chi(g)
\int_0^{\infty} ~dy~
e^{-l\sqrt{1+y^2}}
y^{1+\gamma}
I_0 (j y)
{ ( 1+ y^2 )^{-(b+3\gamma)/4} }
\\
\le
 C_\epsilon
\frac{s^{(3+|\gamma| - b)/2} }{p_0q_0} 
\chi(g)
\int_0^{\infty} ~dy~
e^{-l\sqrt{1+y^2}}
y^{1-\left| \gamma \right|}
I_0 (j y)
{ ( 1+ y^2 )^{-(b-3\left| \gamma \right|)/4} }.
\end{gather*}
Note that we have just used the $\epsilon>0$ from \eqref{cut}.
From \eqref{hypSOFT},  $b \in (0,4-\left| \gamma \right|)$ and $\gamma \in (-2,0)$.  We have in this case 
$b \in (0,\left| \gamma \right|)$.
Hence $b-3\left| \gamma \right|\in (-6,0)$.

We evaluate the relevant integral above as
\begin{gather}
\label{relINT} 
\int_0^{\infty} ~dy~
e^{-l\sqrt{1+y^2}}
y^{1-\left| \gamma \right|}
I_0 (j y)
{ ( 1+ y^2 )^{(3\left| \gamma \right|-b)/4} }
=
\int_0^{1} ~dy~
+
\int_{1}^{\infty} ~dy
\\
\lesssim
\int_{0}^{1} ~dy~
e^{-l\sqrt{1+y^2}}
y^{1-\left| \gamma \right|}
I_0 (j y)
+
\int_{1}^{\infty} ~dy~
e^{-l\sqrt{1+y^2}}
y
I_0 (j y)
{ ( 1+ y^2 )^{(\left| \gamma \right|-b)/4} }.
\notag
\end{gather}
For the unbounded integral we have used the estimate
$
y^{-\left| \gamma \right|}
\lesssim
 ( 1+ y^2 )^{-2\left| \gamma \right|/4}.
$
In this case
$\left| \gamma \right|-b \in (0,2)$ since 
$0< b < \left| \gamma \right|$.

To estimate the remaining integrals above we use the precise theory of special functions, see e.g. \cite{MR0350075,MR950173}. 
We define
$$
\tilde{K}_\alpha(l,j)
\eqdef
\int_0^{\infty} ~dy~
e^{-l\sqrt{1+y^2}}
y
I_0 (j y)
{ ( 1+ y^2 )^{\alpha/4} }.
$$
Then for $\alpha \in [-2,2]$ from \cite[Corollary 1 and Corollary 2]{MR1211782} it is known that
\begin{gather}
\tilde{K}_\alpha(l,j)
\le
C l^{1+\alpha/2}
e^{-c |p-q|}.
\label{kEST}
\end{gather}
We also define
$$
\tilde{I}_\eta(l,j)
\eqdef
\int_0^{1} ~dy~
e^{-l\sqrt{1+y^2}}
y^{1-\eta}
I_0 (j y).
$$
Then for $\eta \in [0,2)$ from \cite[Lemma 3.6]{MR1211782} we have the asymptotic estimate
\begin{gather}
\tilde{I}_\eta(l,j)
\le
C
e^{-c \sqrt{l^2 - j^2}}
\le
Ce^{-c|p-q|/2}.
\label{iEST}
\end{gather}
We will use these estimates in each of the cases below.

Thus in this Case 1 by \eqref{iEST} and \eqref{kEST} the integral in \eqref{relINT} is 
\begin{gather*}
\lesssim
e^{-c|p-q|/2}
+
 l^{1+(\left| \gamma \right|-b)/2}
e^{-|p-q|/2}
\lesssim
 l^{1+(\left| \gamma \right|-b)/2}
e^{-c|p-q|/2}.
\end{gather*}
We may collect the last few estimates together to obtain
\begin{gather*}
k_2^\chi(p,q)
\le
 C_\epsilon
\frac{s^{(3+|\gamma| - b)/2}  }{p_0q_0} 
\chi(g)
\left( p_0+q_0\right)^{1+(\left| \gamma \right|-b)/2}
e^{-c|p-q|/2}.
\end{gather*}
Note that for any $\ell\in\mathbb{R}$ we have
\begin{gather*}
s^{\ell}
e^{-c|p-q|/2}
\le
C_\ell 
e^{-c|p-q|/4}.
\end{gather*}
This follows trivially from \eqref{gESTgs} and   $s = 4+g^2 \le 4+|p-q|^2$. 
Furthermore, for any $\ell \in [0,2]$, we {\it claim} the following estimate
\begin{gather}
\frac{\left( p_0+q_0\right)^{\ell} }{p_0q_0} 
e^{-c|p-q|/4}
\lesssim
\left( p_0q_0\right)^{\ell/2 - 1} 
e^{-c|p-q|/8}.
\label{sumPRODest}
\end{gather}
Using \eqref{sumPRODest} with $\ell = 1 + |\gamma|/2$, in this Case 1, we have the general estimate
\begin{gather*}
k_2^\chi(p,q)
\le
 C_\epsilon
\left( \sqrt{p_0q_0}\right)^{\left| \gamma \right|/2 - 1} 
\left( p_0+q_0\right)^{-b/2}
e^{-c|p-q|}.
\end{gather*}
This is the desired estimate in the current range of exponents for Case 1.

Before moving on to the next case, we establish the {\it claim}.  Suppose that
$$
\frac{1}{2} 	|q| \le |p| \le 2 |q|.
$$
In this case \eqref{sumPRODest} is obvious.  If 
$
\frac{1}{2} 	|q| \ge |p|,
$
then we have
\begin{equation}
|p-q|
\ge 
|q| - |p| 
\ge 
\frac{1}{2} 	|q|.
\label{split1314}
\end{equation}
Whence
\begin{gather}
\frac{\left( p_0+q_0\right)^{\ell} }{p_0q_0} 
e^{-c|p-q|/4}
\le
C_\ell 
\left( p_0q_0\right)^{-1} 
e^{-c|p-q|/8}
e^{-c|q|/64}.
\notag
\end{gather}
In the last splitting
$
|p| \ge 2 |q|,
$
then we alternatively have
\begin{equation}
|p-q|
\ge 
|p| - |q| 
\ge 
\frac{1}{2} 	|p|.
\label{split2314}
\end{equation}
Similarly in this situation 
\begin{gather}
\frac{\left( p_0+q_0\right)^{\ell} }{p_0q_0} 
e^{-c|p-q|/4}
\le
C_\ell 
\left( p_0q_0\right)^{-1} 
e^{-c|p-q|/8}
e^{-c|p|/64}.
\notag
\end{gather}
These last two stronger estimates establish \eqref{sumPRODest}. 
We move on to the next case.

{\bf Case 2}. 
We still consider
   $\gamma= -\left| \gamma \right| < 0$ but now $b+\gamma = b-\left| \gamma \right| \ge 0$.
We have
\begin{multline*}
\left(\frac{\sin \left(\psi/2 \right)}{g}\right)^{b+\gamma}
\cos^\gamma \left(\psi /2 \right)
\lesssim
g^{-b+\left| \gamma \right|}\left(\frac{1}{ ( 1+ y^2 )^{1/4}} \right)^{b-\left| \gamma \right|}
\left(\frac{y}{2\sqrt{1+y^2}}\right)^{-\left| \gamma \right|}
\\
\lesssim
 g^{-(b-\left| \gamma \right|)} y^{-\left| \gamma \right|} { ( 1+ y^2 )^{-(b-3\left| \gamma \right|)/4} }
\\
\le
C_\epsilon ~ y^{-\left| \gamma \right|} { ( 1+ y^2 )^{-(b-3\left| \gamma \right|)/4} }.
\end{multline*}
We have used $g\ge \epsilon$ on the support of $\chi(g)$ in \eqref{cut}.  We also used  \eqref{cosEST} and \eqref{sinEST}.
Then again from \eqref{FIRSTk2est} we have
\begin{gather*}
k_2^\chi(p,q)
\le
C_\epsilon
\frac{s^{3/2} }{p_0q_0} 
\int_0^{\infty} ~dy~
e^{-l\sqrt{1+y^2}}
y^{1-\left| \gamma \right|}
I_0 (j y){ ( 1+ y^2 )^{-(b-3\left| \gamma \right|)/4} }.
\end{gather*}
From \eqref{hypSOFT}, we have in this case that $\gamma \in (-2,0)$ and 
$b \in [\left| \gamma \right|,4-\left| \gamma \right|)$.
Hence for the exponent above $b-3\left| \gamma \right|\in (-4,4)$.

Then the relevant integral above is bounded by
\begin{multline*}
\int_0^{\infty} ~dy~
e^{-l\sqrt{1+y^2}}
y^{1-\left| \gamma \right|}
I_0 (j y){ ( 1+ y^2 )^{-(b-3\left| \gamma \right|)/4} }
=
\int_0^{1} ~dy~
+
\int_{1}^{\infty} ~dy
\\
\lesssim
\int_0^{1} ~dy~
e^{-l\sqrt{1+y^2}}
y^{1-\left| \gamma \right|}
I_0 (j y)
+
\int_{1}^{\infty} ~dy~
e^{-l\sqrt{1+y^2}}
y
I_0 (j y){ ( 1+ y^2 )^{(\left| \gamma \right|-b)/4} }.
\end{multline*}
Here we used the same estimates as in Case 1.
In this case $\left| \gamma \right|-b\in (-4,0)$.
Let $\zeta_2 = \max(-2, \left| \gamma \right|-b)$.  
Then in this case by \eqref{iEST} and \eqref{kEST} the above is
\begin{gather*}
\le
Ce^{-c|p-q|/2}
+
C l^{1+\zeta_2/2}
e^{-|p-q|/2}
\le
C l^{1+\zeta_2/2}
e^{-c|p-q|/2}.
\end{gather*}
We may collect the last few estimates together to obtain
\begin{gather*}
k_2^\chi(p,q)
\le
 C_\epsilon
\frac{s^{3/2} }{p_0q_0} 
 \left(p_0 + q_0\right)^{1+\zeta_2/2}
e^{-c|p-q|/2}
\le
 C_\epsilon
\frac{  \left(p_0 + q_0\right)^{1+\zeta_2/2}}{p_0q_0} 
e^{-c|p-q|/4}.
\end{gather*}
We will further estimate the quotient.

If $\zeta_2 = |\gamma| - b$ then 
$
1+|\gamma|/2 \in [1,2)
$
and 
\eqref{sumPRODest} implies
\begin{gather*}
\frac{  \left(p_0 + q_0\right)^{1+\zeta_2/2}}{p_0q_0} 
e^{-c|p-q|/4}
=
\frac{  \left(p_0 + q_0\right)^{1+|\gamma|/2}}{p_0q_0} 
 \left(p_0 + q_0\right)^{-b/2}
 e^{-c|p-q|/4}
 \\
 \le
  \left(p_0q_0 \right)^{(|\gamma|-2)/4 }
 \left(p_0 + q_0\right)^{-b/2}
 e^{-c|p-q|/8}.
\end{gather*}
Alternatively, if $\zeta_2 = -2$ then 
$
1+\zeta_2/2 =0
$
and
\eqref{sumPRODest} implies
\begin{gather*}
\frac{  \left(p_0 + q_0\right)^{1+\zeta_2/2}}{p_0q_0} 
e^{-c|p-q|/4}
=
\frac{  \left(p_0 + q_0\right)^{b/2}}{p_0q_0} 
 \left(p_0 + q_0\right)^{-b/2}
 e^{-c|p-q|/4}
 \\
\lesssim
  \left(p_0q_0 \right)^{(b-4)/4 }
 \left(p_0 + q_0\right)^{-b/2}
 e^{-c|p-q|/8}.
\end{gather*}
In either situation
\begin{gather*}
k_2^\chi(p,q)
\le
C_\epsilon
  \left(p_0q_0 \right)^{-\zeta }
 \left(p_0 + q_0\right)^{-b/2}
 e^{-c|p-q|/8},
\end{gather*}
with
$
\zeta
\eqdef
\min\left\{2-|\gamma|, 4-b\right\}/4
>0.
$

{\bf Case 3}. 
In this last case
$\gamma = \left| \gamma \right| \ge 0$ and $b+\gamma \ge 0$.  From 
 \eqref{cosEST} and \eqref{sinEST}:
\begin{multline*}
\left(\frac{\sin \left(\psi/2 \right)}{g}\right)^{b+\gamma}
\cos^\gamma \left(\psi /2 \right)
\lesssim
g^{-b-\left| \gamma \right|}\left(\frac{1}{ ( 1+ y^2 )^{1/4}} \right)^{b+\left| \gamma \right|}
\\
\lesssim
 g^{-(b+\left| \gamma \right|)}  { ( 1+ y^2 )^{-(b+\left| \gamma \right|)/4} }
\\
\le
C_\epsilon  { ( 1+ y^2 )^{-(b+\left| \gamma \right|)/4} }.
\end{multline*}
We have again used $g\ge \epsilon$ on the support of $\chi(g)$ in \eqref{cut}.
In this case from \eqref{FIRSTk2est} we have
\begin{gather*}
k_2^\chi(p,q)
\le
C_\epsilon
\frac{s^{3/2}g^{\left| \gamma \right|}  }{p_0q_0} 
\int_0^{\infty} ~dy~
e^{-l\sqrt{1+y^2}}
y
I_0 (j y)
{ ( 1+ y^2 )^{-(b+\left| \gamma \right|)/4} }.
\end{gather*}
From \eqref{hypSOFT},  
$b \in [0,4)$
and
 $b+\left| \gamma \right| \in [0,4+\left| \gamma \right|)$.
Let
$\zeta_3 = \min(2, b+\left| \gamma \right|)\ge 0$.  

Again using \eqref{kEST} we have
\begin{gather*}
\int_0^{\infty} ~dy~
e^{-l\sqrt{1+y^2}}
y
I_0 (j y)
{ ( 1+ y^2 )^{-(b+\left| \gamma \right|)/4} }
\le
C l^{1-\zeta_3/2}
e^{-c|p-q|/2}.
\end{gather*}
Hence
\begin{gather*}
k_2^\chi(p,q)
\le
 C_\epsilon
\frac{s^{3/2}g^{\left| \gamma \right|}  }{p_0q_0} 
\left(p_0 + q_0 \right)^{1-\zeta_3/2}
e^{-c|p-q|/2}
\le
 C_\epsilon
\frac{\left(p_0 + q_0 \right)^{1-\zeta_3/2} }{p_0q_0} 
e^{-c|p-q|/4}.
\end{gather*}
If $\zeta_3 = 2$ this estimate can be handled exactly as in Case 2.  If
$\zeta_3 = b+|\gamma|$:
\begin{multline*}
\frac{\left(p_0 + q_0 \right)^{1-\zeta_3/2} }{p_0q_0} 
e^{-c|p-q|/4}
=
\frac{\left(p_0 + q_0 \right)^{1-|\gamma|/2} }{p_0q_0} 
\left(p_0 + q_0 \right)^{-b/2}
e^{-c|p-q|/4}
\\
\le
\frac{\left(p_0 + q_0 \right)}{p_0q_0} 
\left(p_0 + q_0 \right)^{-b/2}
e^{-c|p-q|/4}
\\
\lesssim
  \left(p_0q_0 \right)^{-1/2 }
\left(p_0 + q_0 \right)^{-b/2} 
e^{-c|p-q|/16}.
\end{multline*}
We have again used the estimate \eqref{sumPRODest}.
In all of the cases we see that $k_2^\chi(p,q)$ satisfies the claimed bound from Lemma \ref{boundK2}
with
$
\zeta
=
\min\left\{2-|\gamma|, 4-b,2\right\}/4
>0.
$
We could obtain a larger $\zeta$ in Case 3 if it was needed.  
\qed \\

With the estimate for the Hilbert-Schmidt form just proven in Lemma \ref{boundK2}, we will now prove the decomposition from Lemma \ref{noKestimate}.  \\

\noindent {\it Proof of Lemma \ref{noKestimate}.} 
We recall the splitting 
$
K
=
K^{1-\chi}
+
K^{\chi}
$
from 
\eqref{kCUT}.
For $K^{\chi}$ we have the kernel $k^\chi = k^\chi_2 - k^\chi_1$ from \eqref{k1def} and \eqref{k2def}.
For a given $R \ge 1$, choose another smooth cut-off function $\phi_R = \phi_R(p,q)$ satisfying 
\begin{gather}
\label{smmothPr}
\phi_R \equiv 1, 
\quad
\text{if} 
\quad
|p|+|q| \le R/2, \quad \left| \phi_R \right| \le 1,
\\
\text{supp} (\phi_R)
\subset
\left\{(p,q) \left| ~
|p|+|q| \le R\right.
\right\}.
\notag
\end{gather}
We will use this cut-off with several different $R$'s in the cases below.
Now we split the kernels $k^\chi(p,q)$ of the operator $K^\chi$  into
\begin{gather*}
k^\chi(p,q) 
= 
k^\chi(p,q)  \phi_R(p,q)
+
k^\chi \left( 1-\phi_R\right)
\\
= 
k_{c}^\chi(p,q)
+
k_{s}^\chi(p,q).
\end{gather*}
We further define 
$$
K_s \eqdef K^{1-\chi}
+
K^{\chi}_s,
$$
where 
$
K^{\chi}_s(h) \eqdef
\int_{\mathbb{R}^3}dq ~ 
k_s^\chi(p,q) ~ h(q).
$
Then the compact part is given by
$$
K_c \eqdef
K^{\chi}_c,
$$
where
$
K^{\chi}_c(h) \eqdef
\int_{\mathbb{R}^3}dq ~ 
k_c^\chi(p,q) ~ h(q).
$
Note that the compactness of $K_c(h)$ is evident from the
integrability of the kernel.  In the following we will show that the operators $K_c$ and $K_s$ satisfy the estimates claimed in Lemma \ref{noKestimate}.

First off, for $K_c$, from the Cauchy-Schwartz inequality  we have
\begin{multline*}
\left|
\langle w^2_\ell  K_c(h_1), h_2\rangle 
\right|
\le
\int_{\mathbb{R}^3}dq ~  \int_{\mathbb{R}^3}dp ~ 
w^2_\ell(p) ~
\left| k_c^\chi(p,q) \right| 
\left|
 h_1(q)
h_2(p)
\right|
\\
\le
\left( \int dq   dp ~ 
w^2_\ell(p) 
\left| k_c^\chi(p,q) \right| ~
\left|
 h_1(q)
\right|^2
\right)^{1/2}
\\
\times
\left( \int dq   dp ~ 
w^2_\ell(p) 
\left| k_c^\chi(p,q) \right| 
\left|
 h_2(p)
\right|^2
\right)^{1/2}.
\end{multline*}
From the definition of $k_c^\chi(p,q)$
and
Lemma \ref{boundK2}, we see that
$$
w^2_\ell(p) \left| k_c^\chi(p,q) \right| 
\le 
C_R ~ e^{-c|p-q|} ~ {\bf 1}_{\le R}(p) ~ {\bf 1}_{\le R}(q),
$$
where ${\bf 1}_{\le R}$ is the indicator function of the ball of radius $R$ centered at the origin as defined in
Lemma \ref{noKestimate}.  By combining the last few estimates we clearly have
the claimed estimate for $K_c$ from Lemma \ref{noKestimate}.

In the remainder of this proof we estimate 
$
K_s= K^{1-\chi}
+
K^{\chi}_s.
$
For
$
K^{\chi}_s 
$
we have
\begin{gather*}
\left|
\langle w^2_\ell  K^{\chi}_s (h_1), h_2\rangle 
\right|
\le
\int_{\mathbb{R}^3}dq ~  \int_{\mathbb{R}^3}dp ~ 
w^2_\ell(p) ~
\left|  k_s^\chi(p,q) \right|  ~
\left|
 h_1(q)
h_2(p)
\right|.
\end{gather*}
With the definition of $k_s^\chi(p,q)$
and
Lemma \ref{boundK2}, we obtain
$$
w^2_\ell(p) \left|  k_s^\chi(p,q) \right|
=
w^2_\ell(p) \left|  k^\chi(p,q) \right|(1-\phi_R)
\lesssim
\frac{ w^2_\ell(p)}{R^{\zeta}} ~ (p_0+q_0)^{-b/2} ~ e^{-c|p-q|},
$$
Furthermore, we {\it claim} that
\begin{gather}
\label{weightCONT}
w_\ell^2(p)e^{-c|p-q|}
\lesssim
w_\ell(p)w_\ell(q)e^{-c|p-q|/2}.
\end{gather}
 By combining the last few estimates including \eqref{weightCONT}, with Lemma \ref{nuEST}, we  have
\begin{multline*}
\left|
\langle w^2_\ell  K^{\chi}_s (h_1), h_2\rangle 
\right|
\lesssim
\frac{1}{R^{\zeta}} 
\int_{\mathbb{R}^3}dq ~  \int_{\mathbb{R}^3}dp ~ 
w_\ell(p)w_\ell(q) ~
\frac{e^{-c|p-q|}}{(p_0 q_0)^{b/4}}
\left|
 h_1(q)
h_2(p)
\right|
\\
\lesssim
\frac{| h_1 |_{\nu,\ell} |  h_2 |_{\nu,\ell}}{R^{\zeta}}.
\end{multline*}
Since $\zeta>0$ we conclude our estimate here by choosing $R>0$ sufficiently large. Notice that the size of this $R$ above clearly depends upon $\epsilon>0$ from \eqref{cut}.

To prove the {\it claim} in \eqref{weightCONT}, we use the same general strategy which was used to prove \eqref{sumPRODest}.  Indeed, if 
$
\frac{1}{2} 	|q| \ge |p|
$
then because of \eqref{split1314} we have
\begin{gather}
w^2_\ell(p)
e^{-c|p-q|/2}
e^{-c|p-q|/2}
\le
w^2_\ell(q)
e^{-c|q|/4}
e^{-c|p-q|/2}
\le C_\ell
e^{-c|p-q|/2},
\notag
\end{gather}
which is better than \eqref{weightCONT}.
Alternatively if 
$
|p| \ge 2 |q|,
$
then with \eqref{split2314} we have
\begin{gather}
w^2_\ell(p)
e^{-c|p-q|/2}
e^{-c|p-q|/2}
\le
w^2_\ell(p)
e^{-c|p|/4}
e^{-c|p-q|/2}
\le
C_\ell
e^{-c|p-q|/2}.
\notag
\end{gather}
The only remaining case is 
$
\frac{1}{2} 	|q| \le |p| \le 2 |q|
$
for which the estimate \eqref{weightCONT} is obvious.

The last term to estimate is $K^{1-\chi}= K^{1-\chi}_2 - K^{1-\chi}_1$.  Notice that 
$$
K^{1-\chi}_1(h) = \int_{\mathbb{R}^3} dq ~ k_1^{1-\chi}(p,q) ~ h(q),
$$
where from \eqref{compactK} and  \eqref{cut}
\begin{equation*}
k_{1}^{1-\chi}(p,q) = 
\left( 1-\chi(g)\right)
\sqrt{J(q)J(p)}
\int_{ \mathbb{S}^{2}}
~ d\omega
~ v_{\o} ~ \sigma(g,\theta).
\end{equation*}
To estimate $K^{1-\chi}_1$ we apply Cauchy-Schwartz to obtain
\begin{gather*}
\left| \langle w_\ell^2  K_{1}^{1-\chi}(h_1), h_2\rangle  \right|
\le
\int_{\mathbb{R}^3}dq ~  \int_{\mathbb{R}^3}dp ~ 
w^2_\ell(p) ~
k_1^{1-\chi}(p,q) ~
\left|
 h_1(q)
h_2(p)
\right|
\\
\le
\left(
\int dp dq~ w_\ell^2(p)
k_1^{1-\chi}(p,q) 
| h_1(q) |^2
\right)^{1/2}
\left(
\int dp dq~ w_\ell^2(p)
k_1^{1-\chi}(p,q) 
| h_2(p) |^2
\right)^{1/2}.
\end{gather*}
We will estimate the kernel of each term above.   We further split
$$
k_1^{1-\chi}(p,q) 
=
k_1^{1-\chi}(p,q) \phi_R(p,q)
+
k_1^{1-\chi}\left( 1 -  \phi_R \right)
=
k_{1S}^{1-\chi}
+
k_{1L}^{1-\chi}.$$
The value of $R\ge 1$ used here is independent of the case considered previously.  The $R$ here will be independent of $\epsilon$.
From \eqref{hypSOFT} and \eqref{moller}, in general we have 
\begin{equation}
\label{generalK}
\int_{ \mathbb{S}^{2}}
~ d\omega
~ v_{\o} ~ \sigma(g,\theta)
\lesssim
 \frac{\sqrt{s}}{p_0 q_0}
g^{1-b} \int_0^{\pi} ~ d\theta ~ \sin^{1+\gamma}\theta
\lesssim \frac{\sqrt{s}}{p_0 q_0}
g^{1-b} .
\end{equation}
For $g\le 2\epsilon$ as  in \eqref{cut}, with \eqref{gESTgs} we conclude that
\begin{equation}
|p-q| \le 2\epsilon \sqrt{p_0 q_0}.
\label{pqBDDD}
\end{equation}
Furthermore, on the support of $\phi_R$ we notice that additionally
$
|p-q| \le 4\epsilon R.
$
Then if $b \in (1,4)$, with the formula for $k_{1}^{1-\chi}(p,q)$, \eqref{generalK} and then \eqref{gESTgs},  we obtain
\begin{multline*}
\int~ dp ~ w_\ell^2(p)
~k_{1S}^{1-\chi}(p,q) 
\le
C_R 
\int_{|p-q|\le 4\epsilon R}~ dp ~ 
\frac{\sqrt{s}}{p_0 q_0}
g^{1-b}
\sqrt{J(q)J(p)}
\\
\le
C_R 
\int_{|p-q|\le 4\epsilon R}~ dp ~ 
g^{1-b}
J^{1/4}(q)J^{1/4}(p)
\\
\le
C_R 
\int_{|p-q|\le 4\epsilon R}~ dp ~ 
\left(\frac{p_0 q_0}{|p-q|} \right)^{b-1}
J^{1/4}(q)J^{1/4}(p)
\le
C R^{4-b} \epsilon^{4-b}
J^{1/8}(q).
\end{multline*}
The last inequality above follows easily from 
$
\left(p_0 q_0\right)^{b-1}
J^{1/8}(q)J^{1/8}(p) \le C
$
and also 
$$
\int_{|p-q|\le 4\epsilon R}~ dp ~ 
\frac{J^{1/8}(p)}{|p-q|^{b-1}} 
\le C R^{4-b} \epsilon^{4-b}.
$$
Thus when $b \in (1,4)$ we have
\begin{gather}
\left| \langle w_\ell^2  K_{1}^{1-\chi}(h_1), h_2\rangle  \right|
\le
C_R  \epsilon^{4-b} | J^{1/16} h_1 |_{2} | J^{1/16} h_2 |_{2}.
\label{keyRbound}
\end{gather}
This is much stronger than the desired estimate for $\epsilon= \epsilon(R) >0$ chosen sufficiently small.
Alternatively if $b \in [0,1]$ then with \eqref{generalK} we have
\begin{gather*}
\int_{\mathbb{R}^3} dp ~ w_\ell^2(p)
~k_{1S}^{1-\chi}(p,q) 
\le
C 
\int_{|p-q|\le 2\epsilon R} dp ~ 
J^{1/4}(q)J^{1/4}(p)
\le
C R^3\epsilon^{3}
J^{1/4}(q).
\end{gather*}
Thus when $b \in [0,1]$ we have
$
\left| \langle w_\ell^2  K_{1}^{1-\chi}(h_1), h_2\rangle  \right|
\le
C_R  \epsilon^{3} | J^{1/8} h_1 |_{2} | J^{1/8} h_2 |_{2}.
$
This concludes our estimate for the part containing $k_{1S}^{1-\chi}(p,q)$ for any fixed $R$ after choosing $\epsilon= \epsilon(R) >0$ small enough (depending on the size of $R$).

For the term involving $k_{1L}^{1-\chi}(p,q)$ the estimate is much easier.  In this case
\begin{multline}
\int_{\mathbb{R}^3} ~ dp ~ w_\ell^2(p)
~k_{1L}^{1-\chi}(p,q) 
\lesssim
\int_{\mathbb{R}^3}~ dp ~ g^{1-b} \left( 1 -  \phi_R(p,q) \right)J^{1/4}(q)J^{1/4}(p)
\\
\lesssim
e^{-R/16}.
\label{RboundINDEP}
\end{multline}
The same estimates hold for the other term in the inner product above.  These estimates are independent of $\epsilon$.  We thus obtain the desired estimate for this term in the same way as for the last term, here we first choose $R>0$ sufficiently large.

The last term to estimate is $K^{1-\chi}_2$ from \eqref{kCUT}.   
With \eqref{pqBDDD} we see that
$$
p_0 \le |p-q| + q_0 \le 2\epsilon \sqrt{p_0 q_0} + q_0
\le
\epsilon p_0 +\left(1 + \epsilon \right) q_0.
$$
The first inequality in this chain  can be found in \cite[Inequality A.1]{MR933458}.
We conclude $p_0 \lesssim q_0$ and similarly $q_0 \lesssim p_0$.  For $0<\epsilon<1/4$ say the constant in these inequalities can be chosen to not depend upon $\epsilon$.
Furthermore from \eqref{0prime}, if $g \le 2\epsilon$ and $\epsilon$ is small (say less than $1/8$), then it is easy to show that
\begin{equation}
p_0' \ge \frac{p_0 + q_0}{4},
\quad 
q_0' \ge \frac{p_0 + q_0}{4}.
\label{0primeBELOW}
\end{equation}
These post-collisional energies are also clearly bounded from above by $p_0$ and $q_0$.   So that all of these variables are comparable on \eqref{kCUT}.
We thus have
\begin{multline*}
\left| \langle 
K_2^{1-\chi}(h_1), h_2 \rangle \right| 
\lesssim
\int_{\mathbb{R}^3\times\mathbb{R}^3\times \mathbb{S}^{2}} ~ d\omega dq dp
~ \left( 1 - \chi(g) \right)
~ v_{\o} ~ \sigma(g,\theta)
e^{-c q_0 - c p_0}
\notag
\\
\qquad
\times \left( \left|  h_1(p^{\prime }) \right| + \left|  h_1(q^{\prime }) \right| \right) \left|  h_2(p) \right|.
\notag
\end{multline*}
With Cauchy-Schwartz we obtain
\begin{multline*}
\lesssim
\left(
\int_{g\le 2\epsilon} ~ d\omega dq dp ~
~ v_{\o} ~ \sigma(g,\theta)
e^{-c q_0 - c p_0}
\left|  h_1(p^{\prime }) \right|^2
\right)^{1/2}
\\
\times
\left(
\int_{g\le 2\epsilon} ~ d\omega dq dp ~
~ v_{\o} ~ \sigma(g,\theta)
e^{-c q_0 - c p_0}
\left|  h_2(p) \right|^2
\right)^{1/2}
\\
+
\left(
\int_{g\le 2\epsilon} ~ d\omega dq dp ~
~ v_{\o} ~ \sigma(g,\theta)
e^{-c q_0 - c p_0}
\left|  h_1(q^{\prime }) \right|^2
\right)^{1/2}
\\
\times
\left(
\int_{g\le 2\epsilon} ~ d\omega dq dp ~
~ v_{\o} ~ \sigma(g,\theta)
e^{-c q_0 - c p_0}
\left|  h_2(p) \right|^2
\right)^{1/2}.
\end{multline*}
From \eqref{generalK} and the arguments just below it, for any small $\eta>0$ we can estimate
\begin{gather*}
\int_{g\le 2\epsilon} ~ d\omega dq ~
~ v_{\o} ~ \sigma(g,\theta)
e^{-c q_0 - c p_0}
\le 
\eta e^{ - c p_0/2}.
\end{gather*}
Above of course we have $\eta = \eta(\epsilon)\to 0$ as $\epsilon \to 0$, and by symmetry the same estimate holds if the roles of $p$ and $q$ are reversed.
Since the kernels of the integrals above are invariant with respect to the relativistic pre-post collisional change of variables \cite{MR1105532}, which is justified for \eqref{postCOLLvelCMsec2}, we may apply it as
$$
dp dq = \frac{p_0' q_0'}{p_0 q_0} dp' dq'. 
$$
Putting all of this together with \eqref{collisionalCONSERVATION} we have
\begin{gather*}
\left| \langle 
K_2^{1-\chi}(h_1), h_2 \rangle \right| 
\le 
\eta(\epsilon)
\left(
\int dp ~
e^{ - c p_0}
\left|  h_1(p) \right|^2
\right)^{1/2}
\left(
\int dp ~
e^{ - c p_0}
\left|  h_2(p) \right|^2
\right)^{1/2}.
\end{gather*}
For $\epsilon>0$ small enough, this is stronger than the estimate which we wanted to prove.  We have now completed the proof of this lemma.
\qed \\

We will at this time use the prior lemma to prove the next lemma.\\

\noindent {\it Proof of Lemma \ref{lowerN}.}    Most of this lemma is standard, see e.g. \cite{MR1379589}.  We only prove the coercive lower bound for the linear operator.  Assuming the converse grants a sequence of functions $h^n(p)$ satisfying 
$
{\bf P} h^n = 0,
$ 
$| h^n |_\nu^2 = \langle \nu h^n, h^n \rangle = 1$ and
$$
\langle {L} h^n, h^n \rangle
=
| h^n |_\nu^2 
- 
\langle {K} h^n, h^n \rangle
\le \frac{1}{n}.
$$ 
Thus $ \{ h^n \} $ is weakly compact in $| \cdot |_\nu$ with limit point $h^0$.  
By weak lower-semi continuity $| h^0 |_\nu \le 1$.  Furthermore,
$$
\langle {L} h^n, h^n \rangle
=
1
- 
\langle {K} h^n, h^n \rangle.
$$ 
We {\it claim} that 
$$
\lim_{n\to\infty}
\langle {K} h^n, h^n \rangle
=
\langle {K} h^0, h^0 \rangle.
$$
The {\it claim} will follow from the prior Lemma \ref{noKestimate}.  
This {\it claim} implies
$$
0 = 1 - \langle {K} h^0, h^0 \rangle.
$$
Or equivalently
$$
\langle {L} h^0, h^0 \rangle
=
| h^0 |_\nu^2 
- 
1.
$$
Since $L\ge 0$, we have $| h^0 |_\nu^2 
=
1$ which implies $h^0 = {\bf P} h^0$.  On the other hand since
 $h^n = \{{\bf I- P}\} h^n$ the weak convergence implies $h^0 = \{{\bf I- P}\} h^0$.
This  is a contradiction to $| h^0 |_\nu^2 
=
1$. 

We now establish the {\it claim}.  For any small $\eta>0$, we split $K = K_c + K_s$ as in Lemma \ref{noKestimate}.  Then
$
\left| \langle {K}_s h^n, h^n \rangle\right| \le \eta.
$
Also $K_c$ is a compact operator in $L^2_\nu$ so that
$$
\lim_{n\to\infty} | K_c h^n - K_c h^0 |_\nu = 0.
$$
We conclude by first choosing $\eta$ small and then sending $n\to\infty$.
\qed \\

We are now ready to prove Lemma \ref{lowerA}.  We point out that
similar estimates, but with strong Sobolev norms,
have been established in recent years \cite{MR1946444,MR1908664,MR2000470,MR2100057} via  the macroscopic
equations for the coefficients $a,b$ and $c$.
We will use the approach from \cite{G2}, which exploits the hyperbolic nature of the transport operator, to prove our Lemma \ref{lowerA} in the low regularity $L^2$ setting.  
\\

\noindent {\it Proof of Lemma \ref{lowerA}.}  
We use the method of contradiction, if Lemma \ref{lowerA}   is not valid then for any $k\ge 1$ we can find a sequence of normalized solutions to \eqref{rBlin} which we denote by $f_{k}$ that satisfy
$$
\int_{0}^{1} ~ ds ~ \|\{\mathbf{I}-\mathbf{P}\}f_{k}\|_{\nu }^{2}(s)
\leq 
\frac{1}{k}
\int_{0}^{1} ~ ds ~\|\mathbf{P}f_{k}\|_{\nu }^{2}(s).
$$
Equivalently the normalized function
$$
Z_{k}(t,x,p)\eqdef \frac{f_{k}(t,x,p)}{\sqrt{\int_{0}^{1} ~ ds ~\|\mathbf{P}f_{k}\|_{\nu }^{2}(s)}},
$$
satisfies 
$$
\int_{0}^{1} ~ ds ~ \|\mathbf{P}Z_{k}\|_{\nu }^{2}(s) =  1,
$$ 
and 
\begin{equation}
\int_{0}^{1} ~ ds ~ \|\{\mathbf{I}-\mathbf{P}\}Z_{k}\|_{\nu }^{2}(s)
\leq 
\frac{1}{k}.
\label{1/n}
\end{equation}
Moreover, from \eqref{conservation} the following integrated conservation laws hold
\begin{gather}
\int_0^1~ ds ~ 
\int_{\mathbb{T}^3}  dx ~
\int_{\mathbb{R}^3}   dp ~
\begin{pmatrix}
      1   \\      p  \\ p_0
\end{pmatrix}
~
\sqrt{J(p)}~
Z_k(s,x,p) 
=
0. 
 \label{conserveZ}
\end{gather}
Furthermore, since $f_{k}$ satisfies \eqref{rBlin}, so does $Z_k$.
Clearly
\begin{gather}
\sup_{k\ge 1}\int_{0}^{1}~ ds ~ \|Z_{k}\|_{\nu }^{2}(s)
\lesssim 1.
\label{boundZk}
\end{gather}
Hence there exists $Z(t,x,p)$ such that 
$$
 Z_k(t,x,p)
 \rightharpoonup 
Z(t,x,p),
\quad
\;{\rm as }\;k\rightarrow \infty , 
$$
weakly with respect to the inner product $\int_0^1~ ds ~ (\cdot, \cdot)_\nu$
of the norm $\int_{0}^{1} ~ ds ~ \|\cdot \|_{\nu }^{2}$.
Furthermore, from (\ref{1/n}) we know that
\begin{equation}
\int_{0}^{1} ~ ds ~ \|\{\mathbf{I}-\mathbf{P}\}Z_{k}\|_{\nu }^{2}(s)\rightarrow 0.
\label{i-pto0}
\end{equation}%
We conclude that $\{\mathbf{I}-\mathbf{P}\}Z_{k}\rightarrow \{\mathbf{I}-\mathbf{P}\}Z$
 and $\{\mathbf{I}-\mathbf{P}\}Z=0$ from (\ref{i-pto0}).
It is then straightforward to verify that
\begin{equation*}
\mathbf{P}Z_{k}\rightarrow \mathbf{P}Z\text{ weakly in }
\int_{0}^{1} ~ ds ~ \|\cdot \|_{\nu }^{2}.
\end{equation*}%
  Hence  
\begin{equation}
Z(t,x,p)
=
\mathbf{P}Z
=
\{a(t,x)+p\cdot b(t,x)+p_0 c(t,x)\}\sqrt{J}.  \label{zabc}
\end{equation}%
At the same time notice that
$
LZ_{k}=L\{\mathbf{I}-\mathbf{P}\}Z_{k}
$ 
and we have \eqref{i-pto0}.
Send
$k\rightarrow \infty $ in \eqref{rBlin} for $Z_{k}$ to obtain, in the sense of distributions, that
\begin{equation}
\partial _{t}Z+\hat{p}\cdot \nabla _{x}Z=0.  \label{zvlasov}
\end{equation}%
At this point our main strategy is to show, on the one hand, $Z$ has to be zero from (\ref{i-pto0}), the periodic boundary conditions, and the hyperbolic transport equation \eqref{zvlasov}, and \eqref{conserveZ}. On the other hand, $%
Z_{k}$ will be shown to converge strongly to $Z$ in 
$
\int_{0}^{1}~ ds ~\|\cdot \|_{\nu}^{2}
$  
with the help of the averaging lemma \cite{MR1003433} in the relativistic formulation \cite{MR2119935}
and 
$
\int_{0}^{1}~ ds ~\|Z\|_{\nu}^{2}(s)> 0.
$ 
This would be a contradiction.

{\bf Strong Convergence.}
We begin by proving the strong convergence, and then later we will prove that the limit is zero.
Split $Z_k(t,x,p)$ as 
\[
Z_k(t,x,p)={\bf P} Z_k+\{{\bf I-P}\}Z_k=\sum_{j=1}^5\langle Z_k(t,x,\cdot
),e_j\rangle e_j(p)+\{{\bf I-P}\}Z_k, 
\]
where $e_j(p)$ are an orthonormal basis for \eqref{null} in $\| \cdot \|_\nu$.

To prove the strong convergence in 
$
\int_0^1 ~ ds ~ 
\|\cdot\|_\nu^2,
$
recalling \eqref{i-pto0},  
we will show
\begin{equation}
\sum_{1\le j\le 5}\int_0^1 ~ ds ~ 
\|\langle  Z_k,e_j\rangle e_j-\langle  Z,e_j\rangle e_j\|_\nu^2(s)
\rightarrow
0.
\notag
\end{equation}
Since $e_j(p)$ are smooth with exponential decay when $p\rightarrow \infty ,$
it suffices to prove
\begin{equation}
\int_0^1 ~ ds ~
\int_{{\mathbb T}^3} ~ dx ~ 
|\langle  Z_k,e_j\rangle -\langle
 Z,e_j\rangle |^2
 \rightarrow 0.  \label{strong0}
\end{equation}
We will now establish \eqref{strong0} using the averaging lemma.

Choose any small $\eta >0$ and a smooth cut off
function $\chi _1(t,x,p)$ in $(0,1)\times {\mathbb T}^3\times {\mathbb R}^3$ such
that $\chi _1(t,x,p)\equiv 1$ in $[\eta ,1-\eta ]\times {\mathbb T}^3\times
\left\{|p|\le \frac 1\eta \right\}$
and  
$\chi _1(t,x,p)\equiv 0$ outside $[\eta/2 ,1-\eta/2 ]\times {\mathbb T}^3\times
\left\{|p|\le \frac 2\eta \right\}$.
Split
\begin{equation}
\langle  Z_k(t,x,\cdot ),e_j\rangle 
=
\langle \left(1-\chi_1\right)  Z_k(t,x,\cdot ),e_j\rangle +\langle \chi _1 Z_k(t,x,\cdot ),e_j\rangle .  \label{splitn}
\end{equation}
For the first term above, notice that 
\begin{multline*}
\int_0^1 ~ ds ~
\int_{{\mathbb T}^3} ~ dx ~ 
|\langle \left(1-\chi_1\right)
\left| Z_k-Z \right|,e_j\rangle |^2
\\
\lesssim
\int_0^1\int_{{\mathbb T}^3\times {\mathbb R}^3} 
\left(1-\chi_1\right)^2
| Z_k|^2|e_j|^2
+
\int_0^1\int_{{\mathbb T}^3\times {\mathbb R}^3} 
\left(1-\chi_1\right)^2
| Z|^2|e_j|^2
\\
\lesssim
\int_{0\le s\le \eta }+\int_{1-\eta \le s\le 1}+\int_{|p|\ge 1/\eta
}.
\end{multline*}
Since $e_j = e_j(p)$ has exponential decay in $|p|$ we have the crude bound
$$
|e_j(p)|\le C ~\eta, 
\quad 
{\rm for }
\quad 
|p|\ge 1/\eta.
$$
Thus
all three integrals above can be bounded by 
\begin{equation}
C\eta  \sup_{0\le s \le 1}\left(
 \|Z_k\|_{2}^2(s)+\|Z\|_{2}^2(s) 
 \right)
 \le C
 \left(
 \|Z_k(0) \|_{2}^2+\|Z(0) \|_{2}^2
 \right) \eta
 \le C \eta,
\label{three}
\end{equation}
which will hold for $Z_k$ uniformly in $k$.
These bounds follow from \eqref{rBlin} and \eqref{zvlasov}.

The second term 
in (\ref{splitn}), $\langle \chi _1
Z_k(t,x,\cdot ),e_j\rangle $, is actually uniformly bounded
in $H^{1/4}([0,1]\times {\mathbb T}^3).$ 
To prove this, notice that  
\eqref{rBlin} implies that $\chi _1 Z_k$ satisfies 
\begin{equation}
[\partial _t+\hat{p}\cdot \nabla _x]\left(\chi _1 Z_k\right)=-\chi
_1L[ Z_k]+ Z_k[\partial _t+\hat{p}\cdot\nabla _x]\chi _1.
\label{average}
\end{equation}
The goal is to show that each term on the right hand side of (\ref{average})
is uniformly bounded in $L^2([0,1]\times {\mathbb T}^3\times {\mathbb R}^3)$.
This would imply the $H^{1/4}$ bound by the averaging lemma.
It clearly follows from \eqref{boundZk} that
\[
 Z_k[\partial _t+ \hat{p}\cdot \nabla _x]\chi _1\in L^2([0,1]\times 
{\mathbb T}^3\times {\mathbb R}^3). 
\]
Furthermore, it follows from
Lemma \ref{nuEST} and
Lemma \ref{noKestimate} with $\ell =0$ that
\begin{gather*}
\int_{\mathbb{R}^3} dp~ |\chi _1L[ Z_k]|^2
\le 
\int_{\mathbb{R}^3} dp~ |\chi _1\nu(p)  Z_k|^2
+
\int_{\mathbb{R}^3} dp~ |\chi _1K(Z_k)|^2
\lesssim
 | Z_k(t,x)|_\nu^2.
\end{gather*}
Thus the right hand side of 
(\ref{average}) 
is uniformly bounded in $L^2([0,1]\times {\mathbb T}^3\times 
{\mathbb R}^3).$ By the averaging lemma \cite{MR1003433,MR2119935} it follows that 
\[
\langle \chi _1 Z_k(t,x,\cdot ),e_j\rangle 
=
\int_{{\mathbb R}^3} dp ~ \chi _1(t,x,p) Z_k(t,x,p)e_j(p)
\in 
H^{1/4}([0,1]\times {\mathbb T}^3). 
\]
This holds uniformly in $k$, which implies up to a subsequence that
\[
\langle \chi _1 Z_k(t,x,\cdot ),e_j\rangle \rightarrow
\langle \chi _1 Z(t,x,\cdot ),e_j\rangle 
\;{\rm in }\;
L^2([0,1]\times {\mathbb T}^3). 
\]
Combining this last convergence with \eqref{three} concludes the proof of \eqref{strong0}.

As a consequence of this strong convergence we have
$$
\int_0^1 ~ ds ~ \| Z_k- Z\|_\nu^2(s) \rightarrow 0,
$$
which implies that
\begin{equation}
\int_0^1 ~ ds ~ \| {\bf P} Z\|_\nu ^2(s) = 1.  \notag
\end{equation}
Now if we can show that at the same time ${\bf P} Z = 0$, then we have a contradiction.

{\bf The Limit Function $Z(t,x,p) = 0$.}  By analysing the equations satisfied by $Z$, we will show that $Z$ must be trivial.
We will now derive the macroscopic equations for ${\bf P}Z$'s coefficients $a ,b$ and 
$c. $ 
Since $\{ {\bf I- P} \}Z = 0$, we see that ${\bf P}Z$ solves \eqref{zvlasov}.  We plug 
the expression for ${\bf P}Z$ in  \eqref{zabc} into the equation \eqref{zvlasov}, and expand in the basis \eqref{null} to obtain
\[
\left\{\partial^0 a+
\frac{p_j}{p_0}\left\{ \partial^j a\right\}
+\frac{p_jp_i}{p_0}\left\{ \partial^i b_j \right\}+p_j\left\{\partial^0 b_j+\partial^j c\right\}
+p_0 \left\{ \partial^0 c \right\} \right\}J^{1/2}(p)
=
0.
\]
where $\partial ^0=\partial _t$ and $\partial ^j=\partial _{x_j}$. 
By a comparison of coefficients, we obtain the important relativistic
macroscopic equations for 
$a(t,x),$ $b_i(t,x)$ and  $c(t,x)$: 
\begin{eqnarray}
\partial ^0c 
&=& 0,  \label{bi} \\
\partial ^i c+\partial^0 b_i
&=& 0, \label{c} \\
(1-\delta_{ij})\partial ^ib_j+\partial ^jb_i
&=& 0, \label{bij} \\
\partial^ia
&=& 0,
\label{ai} \\
\partial ^0a
&=& 0,  \label{adot}
\end{eqnarray}
which hold in the sense of distributions.  

We will show that these equations 
\eqref{bi} - \eqref{adot}, combined with the periodic boundary conditions imply that any solution to \eqref{bi} - \eqref{adot} is a constant.  Then the conservation laws \eqref{conserveZ} will imply that the constant can only be zero.

We deduce from \eqref{adot} and \eqref{ai}
that 
\begin{gather*}
a(t,x) = a(0,x), 
\quad
a.e.~  x,t
\\
a(s,x_1) = a(s, x_2),
\quad
a.e. ~  s, x_1, x_2.
\end{gather*}
Thus $a$ is a constant for almost every $(t,x)$.  From \eqref{bi}, we have 
$c(t,x) = c(x)$ for a.e. $t$.
Then from \eqref{c} for some spatially dependent function $\tilde{b}_i(x)$ we have
$$
b_i(t,x) = \partial^i c(x) t + \tilde{b}_i(x).
$$
From \eqref{bij} and the above
$$
0 = \partial^i b_i(t,x) = \partial^i \partial^i c(x) t + \partial^i \tilde{b}_i(x).
$$
which implies 
$$
\partial^i \partial^i c(x) = 0,
\quad 
\partial^i \tilde{b}_i(x) = 0.
$$
Similarly if $i\ne j$ we have
$$
0 = \partial^j b_i(t,x)+ \partial^i b_j(t,x) =\left(  \partial^j \partial^i c(x) +  \partial^i \partial^j c(x) \right) t + \partial^j \tilde{b}_i(x) + \partial^i \tilde{b}_j(x),
$$
so that
$$
\partial^j \partial^i c(x) = - \partial^i \partial^j c(x) ,
\quad 
\partial^j \tilde{b}_i(x) = - \partial^i \tilde{b}_j(x),
$$
which implies 
$
\partial^ i c(x) = c_i ,
$
and $c(x)$ is a polynomial.  
By the periodic boundary conditions $c(x) = \tilde{c}\in\mathbb{R}$.

We further observe that $b_i(t,x) = b_i(x)$ is a constant in time a.e. from \eqref{c} and the above.  From \eqref{bij} again
$\partial^i b_i = 0$ so that trivially $\partial^i\partial^i b_i = 0$.  Moreover \eqref{bij} further implies that 
$\partial^j\partial^j b_i = 0$.  Thus for each $i$,  $b_i(x)$ is a periodic polynomial,  which must be a constant: $b_i(x) = b_i\in\mathbb{R}$.

We compute from \eqref{zabc} that
\begin{gather*}
\int_{{\mathbb R}^3} ~
 dp ~ p_i
J^{1/2}(p)~
Z(t,x,p) 
=
b_i\int_{{\mathbb R}^3} ~
 dp ~ p_i^2
J(p),
\quad 
i=1,2,3, 
\\
\int_{{\mathbb R}^3} ~
 dp ~ 
J^{1/2}(p)~
Z(t,x,p) 
=
a\int_{{\mathbb R}^3} ~
 dp ~ 
J(p)
+
c
\int_{{\mathbb R}^3} ~
 dp ~ 
p_0 J(p),
\\
\int_{{\mathbb R}^3} ~
 dp ~ p_0
J^{1/2}(p)~
Z(t,x,p) 
=
a\int_{{\mathbb R}^3} ~
 dp ~ p_0
J(p)
+
c
\int_{{\mathbb R}^3} ~
 dp ~ 
p_0^2 J(p).
\end{gather*}
As in \cite{MR2100057} 
we define 
\begin{gather*}
\rho_1=\int_{\mathbb{R}^3} J(p) dp=1,
\quad
\rho_0=\int_{\mathbb{R}^3} p_0J(p) dp,
\quad
\rho_2 =\int_{\mathbb{R}^3} p_0^2J(p) dp.
\end{gather*}
Now the matrix given by 
$
\begin{pmatrix}
      \rho_1 & \rho_0   \\      \rho_0  & \rho_2
\end{pmatrix}
$
is invertible because $\rho_0^2 < \rho_1 \rho_2$.
It then follows from the conservation law
\eqref{conserveZ} which is satisfied by the limit function $Z(t,x,p)$
that the constants $a$, $b_i$, $c$ must indeed be zero.
\qed \\

We have now completed all of the $L^2$ energy estimates for the linearized relativistic Boltzmann equation \eqref{rBlin}.  We will now use \eqref{energyWest}, Lemma \ref{lowerN},
Lemma \ref{lowerA},
and Lemma \ref{weight2}
to prove
 Theorem \ref{decay2}.
This will be the final proof in this section.

\bigskip

\noindent
{\it Proof of Theorem \ref{decay2}.}
For $k\ge 0$, we define the time weight function by
\begin{equation}
P_k(t) \eqdef
(1+t)^k.
\label{timeW}
\end{equation}
For a solution $f(t,x,p)$ to the linear
Boltzmann equation (\ref{rBlin}), with $P_k(t)$ from \eqref{timeW}, $P_k(t)f(t)$ satisfies the equation
\begin{equation}
\left(
\partial _{t}+\hat{p}\cdot \nabla _{x}+L
\right)
\left(P_k(t)f(t)\right)
-
k P_{k-1}(t)f(t)
=0.  \label{lambdaf}
\end{equation}
For the moment, suppose that $t=m$  and 
$m\in \{1,2,3,\ldots\}$.
For the time interval $[0,m]$,
we multiply $P_k(t)f(t)$
with (\ref{lambdaf}) and take the $L^{2}$ energy estimate over $0\leq s\leq m$ to obtain
\begin{multline}
\label{syleA}
P_{2k}(m)
\|f\|_{2}^{2}(m)
+
\int_{0}^{m} ~ ds ~
P_{2k}(s)(Lf,f)
\\
-
k
\int_{0}^{m} ~ ds ~ P_{2k-1}(s)
\|f\|_{2}^{2}(s)
=
\|f_0\|_{2}^{2}.
\end{multline}
We divide the time interval into $\cup _{j=0}^{m-1}[j,j+1)$ and also 
$f_{j}(s,x,p)\eqdef f(j+s,x,p)$ for $j\in\{0,1,2,\ldots,m-1\}$. We have
\begin{gather*}
P_{2k}(m)
\|f(m)\|_{2}^{2}
+
\sum_{j=0}^{m-1}
\int_{0}^{1} ~ ds ~
\left\{ 
P_{2k}(j+s)
(Lf_{j},f_{j})
-
k P_{2k-1}(j+s)
\|f_{j}\|_{2}^{2}(s)
\right\}
\notag
\\
=
\|f_0\|_{2}^{2}.  
\end{gather*}
Clearly $f_{j}(s,x,p)\,$ satisfies the same linearized Boltzmann
equation \eqref{rBlin} on the interval $0\leq s\leq 1$.
Notice that on this time interval
\begin{gather}
P_{2k}(j+s)
\ge P_{2k}(j),
~
P_{2k-1}(j+s)
\le
\tilde{C}_k P_{2k-1}(j),
~
\forall
k\ge 1/2,
~
s\in[0,1].
\label{timeWb} 
\end{gather}
These estimates are uniform in $j$, so that
\begin{multline*}
P_{2k}(m)
\|f(m)\|_{2}^{2}
+
\sum_{j=0}^{m-1}
\int_{0}^{1} ~ ds ~
\left\{ 
P_{2k}(j)
(Lf_{j},f_{j})(s)
-
\tilde{C}_k k P_{2k-1}(j)
\|f_{j}\|_{2}^{2}(s)
\right\}
\notag
\\
\le
\|f_0\|_{2}^{2}.  
\end{multline*}
Moreover,
by Lemma \ref{lowerN}, we have
\begin{equation*}
(Lf_{j},f_{j})
\geq 
\delta _{0}\|\{\mathbf{I-P}\}f_{j}\|_{\nu }^{2}
=
\frac{\delta _{0}}{2} \|\{\mathbf{I-P}\}f_{j}\|_{\nu }^{2}
+\frac{\delta _{0}}{2} \|\{\mathbf{I-P}\}f_{j}\|_{\nu }^{2}.
\end{equation*}
Furthermore, 
with Lemma \ref{lowerA} applied  to each $f_{j}(s,x,p)$ we obtain
\begin{equation}
\frac{\delta _{0}}{2}
\sum_{j=0}^{m-1} P_{2k}(j) \int_{0}^{1} ~ ds ~
\| \{\mathbf{I - P}\}f_{j} \|_{\nu }^{2}
\geq 
\frac{\delta _{0}\delta_\nu}{2}
\sum_{j=0}^{m-1} P_{2k}(j) \int_{0}^{1} ~ ds ~
\|\mathbf{P}f_{j}\|_{\nu }^{2}(s).  
\label{deltalower}
\end{equation}
We combine \eqref{deltalower} with the estimate above it to conclude 
\begin{multline*}
\sum_{j=0}^{m-1}
P_{2k}(j)
\int_{0}^{1} ~ ds ~
(Lf_{j},f_{j})(s)
\ge
\frac{\delta _{0}}{2} \sum_{j=0}^{m-1}
P_{2k}(j)
\int_{0}^{1} ~ ds ~ \|\{\mathbf{I-P}\}f_{j}\|_{\nu }^{2}(s)
\\
+
\frac{\delta _{0}\delta_\nu}{2}\sum_{j=0}^{m-1}
P_{2k}(j)
\int_{0}^{1} ~ ds ~ \|\mathbf{P}f_{j}\|_{\nu }^{2}(s)
\\
\ge
\tilde{\delta}  \sum_{j=0}^{m-1}
P_{2k}(j)
\int_{0}^{1} ~ ds ~ \| f_{j}\|_{\nu }^{2}(s).
\end{multline*}
where
$
\tilde{\delta} 
=
\frac{1}{2}\min\left\{
\frac{\delta _{0} \delta _{\nu}}{2},
\frac{\delta _{0} }{2}
\right\}.
$
Define $C_k = \tilde{C}_k k$.
With this lower bound 
\begin{gather*}
P_{2k}(m)
\| f \|_{2}^{2}(m)
+
\sum_{j=0}^{m-1}
 \int_{0}^{1}  ds 
 \left\{
\tilde{\delta}
 P_{2k}(j)
\|f_{j}\|_{\nu }^{2}
-
 C_k
  P_{2k-1}(j)
\|f_{j} \|_{2}^{2} 
\right\}(s)
\leq 
\| f_0 \|_{2}^{2}.
\end{gather*}
Next, for $\lambda>0$ sufficiently small we introduce the following splitting
\begin{gather}
E_{\lambda,j} = 
\left\{p\left|
p_0^{b/2} <  \lambda \left(1+j\right)
\right\}\right.,
\quad
E^c_{\lambda,j}  = 
\left\{p\left|
p_0^{b/2} \ge  \lambda \left(1+j\right)
\right\}\right..
\label{mainSplit3}
\end{gather}
We incorporate this splitting into our energy inequality as follows
\begin{multline*}
P_{2k}(m)
\| f \|_{2}^{2}(m)
+
\frac{\tilde{\delta} }{2}
\sum_{j=0}^{m-1}
 P_{2k}(j)
 \int_{0}^{1} ~ ds ~
\|f_{j}\|_{\nu }^{2}(s)
\\
+
\sum_{j=0}^{m-1}
 \int_{0}^{1} ~ ds ~
 \left\{
\frac{\tilde{\delta} }{2}
\| \sqrt{P_{2k}(j)}  f_{j}{\bf 1}_{E_{\lambda,j} }\|_{\nu }^{2}(s)
-
 C_k
\|  \sqrt{P_{2k-1}(j)}  f_{j} {\bf 1}_{E_{\lambda,j} } \|_{2}^{2} (s)
\right\}
\\
\leq 
\| f_0 \|_{2}^{2}
+
\sum_{j=0}^{m-1}
 \int_{0}^{1} ~ ds ~
 C_k
\|   \sqrt{P_{2k-1}(j)} f_{j} {\bf 1}_{E_{\lambda,j} ^c} \|_{2}^{2} (s),
\end{multline*}
where 
$
{\bf 1}_{E_{\lambda,j} }
$
is the usual indicator function of the set 
$
{E_{\lambda,j} }.
$
On $E_{\lambda,j} $, with the help of Lemma \ref{nuEST}, we have
$$
-1 \ge -\lambda  \left(1+j\right)p_0^{-b/2}
\ge -C_\nu \lambda  \left(1+j\right)\nu(p),
\quad 
C_\nu > 0.
$$
Hence
\begin{multline*}
\frac{\tilde{\delta} }{2}
\| \sqrt{P_{2k}(j)}  f_{j}{\bf 1}_{E_{\lambda,j} }\|_{\nu }^{2}(s)
-
 C_k
\|  \sqrt{P_{2k-1}(j)}  f_{j} {\bf 1}_{E_{\lambda,j} } \|_{2}^{2} (s)
\\
\ge
\left(
\frac{\tilde{\delta} }{2}
-
C_k C_\nu \lambda
\right)
\| \sqrt{P_{2k}(j)}  f_{j}{\bf 1}_{E_{\lambda,j} }\|_{\nu }^{2}(s).
\end{multline*}
We choose $\lambda>0$ small enough such that
$$
C_\lambda
\eqdef
\frac{\tilde{\delta} }{2}
-
C_k C_\nu \lambda
>0.
$$
Then we have the following useful energy inequality
\begin{multline*}
P_{2k}(m)
\| f \|_{2}^{2}(m)
+
\frac{\tilde{\delta} }{2}
\sum_{j=0}^{m-1}
 P_{2k}(j)
 \int_{0}^{1} ~ ds ~
\|f_{j}\|_{\nu }^{2}(s)
\\
+
C_\lambda
\sum_{j=0}^{m-1}
 \int_{0}^{1} ~ ds ~
\| \sqrt{P_{2k}(j)}  f_{j}{\bf 1}_{E_{\lambda,j} }\|_{\nu }^{2}(s)
\\
\leq 
\| f_0 \|_{2}^{2}
+
\sum_{j=0}^{m-1}
 \int_{0}^{1} ~ ds ~
 C_k
\|   \sqrt{P_{2k-1}(j)} f_{j} {\bf 1}_{E_{\lambda,j} ^c} \|_{2}^{2} (s).
\end{multline*}
The last term on the left side of the inequality is positive and we  discard it from the energy inequality.  
For the right side of the energy inequality, on the complementary set $E_{\lambda,j} ^c$, using Lemma \ref{nuEST} again, we have
$$
P_{2k-1}(j) 
\le
\left(\frac{p_0^{b/2}}{\lambda} 
\right)^{2k-1}
\le 
\frac{C_\nu }{\lambda^{2k-1} }
~ \nu(p) ~ w_{2k}(p).
$$
Thus we bound the time weights with velocity weights and the dissipation norm
\begin{gather*}
\sum_{j=0}^{m-1}
 \int_{0}^{1} ~ ds ~
 C_k
\|  \sqrt{ P_{2k-1}(j)} f_{j} {\bf 1}_{E_{\lambda,j} ^c} \|_{2}^{2} (s)
\le
C
\sum_{j=0}^{m-1}
 \int_{0}^{1} ~ ds ~
\|   f_{j} {\bf 1}_{E_{\lambda,j} ^c} \|_{\nu, k}^{2} (s).
\end{gather*}
We switch back to $f_{j}(t,x,p)=f(t+j,x,p)$ and use 
\eqref{timeWb} 
to deduce 
\begin{gather*}
P_{2k}(m)
\| f \|_{2}^{2}(m)
+
\delta_k 
 \int_{0}^{m} ~ ds ~
  P_{2k}(s)
\|f\|_{\nu }^{2}(s)
\leq 
\| f_0 \|_{2}^{2}
+
C
 \int_{0}^{m} ~ ds ~
\|   f  \|_{\nu, k}^{2} (s).
\end{gather*}
We can obtain an upper bound for the right side above using the regular energy inequality in
\eqref{energyWest} to achieve
\begin{gather}
P_{2k}(m)
\| f \|_{2}^{2}(m)
+
\delta_k 
 \int_{0}^{m} ~ ds ~
  P_{2k}(s)
\|f\|_{\nu }^{2}(s)
\leq 
C_k \| f_0 \|_{2,k}^{2}.
\notag
\end{gather}
We have thus established our desired energy inequality from Theorem \ref{decay2} for any
$
m \in \{ 0,1,2,\dots \}
$
and 
$\ell = 0$.
For an arbitrary $t>0$, we choose
$
m \in \{ 0,1,2,\dots \}
$
such that
 $m\leq t\leq m+1$. 
 We then split the time integral as 
 $[0,t]=[0,m]\cup \lbrack m,t].$
 
For the time interval $[m,t],$ we have the $L^{2}$ energy estimate as in \eqref{mtESt}.
Since $L\ge 0$ by Lemma \ref{lowerN}, we may use 
\eqref{timeWb} and \eqref{mtESt} to see that
\begin{equation*}
P_{2k}(t)\|f(t)\|_{2}^{2}
\le
C_k P_{2k}(m)
\|f(m)\|_{2}^{2},\quad \forall t \in[m,m+1].  
\end{equation*}
Since \eqref{syleA} holds for any time $t$ (not necessarily an integer), we can use the estimate above together with \eqref{syleA}, as in \eqref{lowerMM}, using Lemma \ref{weight2}, 
%
for any $t> 0$ to obtain
\begin{gather}
P_{2k}(t)
\| f \|_{2}^{2}(t)
+
\delta_k 
 \int_{0}^{t} ~ ds ~
  P_{2k}(s)
\|f\|_{\nu }^{2}(s)
\leq 
C_k \| f_0 \|_{2,k}^{2}.
\label{nTdecay}
\end{gather}
This proves our time decay  theorem for $\ell =0$.  For general $\ell>0$ this estimate can be proven in exactly the same way, except in this case we use Lemma \ref{weight2}  in the place of Lemma \ref{lowerN} and 
Lemma \ref{lowerA} as we did in the proof of \eqref{energyWest}.
\qed \\

This concludes our discussion of $L^2$ estimates for the linear Boltzmann equation.  In the next section we use these 
$L^2$ estimates to prove $L^\infty$ estimates.

\section{Linear $L^\infty$ Bounds and slow Decay}\label{secL0}

In this section we will prove global in time uniform bounds for solutions to the linearized equation \eqref{rBlin} in 
$L^\infty([0,\infty)\times\mathbb{T}^3\times \mathbb{R}^3 )$, and slow polynomial decay in time.  
We express solutions, $f(t,x,p)$, to \eqref{rBlin} with the semigroup $U(t)$ as 
\begin{equation}
f(t,x,p) = \{U(t)f_0\}(x,p),
\label{Udef}
\end{equation}
with initial data given by
$$
\{U(0)f_0\}(x,p) = f_0(x,p).
$$
Our goal in this section will be to prove the following.

\begin{theorem}
\label{decay0l}
Given $\ell >3/b$ and $k\in[0,1]$.  Suppose that
$
f_0 \in L^\infty_{\ell+k} (\mathbb{T}^3\times \mathbb{R}^3 )
$
 satisfies \eqref{conservation} initially,
then under \eqref{hypSOFT} the semi-group satifies
$$
\|\{U(t)f_0\} \|_{\infty,\ell}
\le
C (1+t)^{-k}\| f_0\|_{\infty,\ell+k}.
$$
Above the positive constant $C= C_{\ell, k}$ only depends on $\ell$ and $k$.  
\end{theorem}

The first step towards proving Theorem \ref{decay0l} is an appropriate decomposition.  
Initially we consider solutions to the linearization of \eqref{rBoltz} with the compact operator $K$ removed from \eqref{rBlin}. This equation is given by
\begin{gather}
\left( \partial_t  + \hat{p}\cdot \nabla_x  + \nu(p) \right) f
=
0,
\quad
f(0,x,p) =  f_0(x,p).
\label{rBlinWO}
\end{gather}
Let the semigroup 
$
G(t)f_0
$ 
denote the solution to this system \eqref{rBlinWO}.
Explicitly
$$
G(t)f_0(x,p) \eqdef e^{-\nu(p) t } f_0(x-\hat{p}t, p).
$$
For soft potentials \eqref{hypSOFT}, with Lemma \ref{nuEST},  this formula does not imply exponential decay in $L^\infty$ for high momentum values.  However, as we will see in Lemma \ref{hESTIMATES} below, this formula does imply that one can trade between arbitrarily high polynomial decay rates and additional polynomial momentum weights on the initial data.

More generally we consider solutions to the full linearized system \eqref{rBlin}, which are expressed with the semi-group \eqref{Udef}.
By the Duhamel formula 
\begin{gather*}
\left\{  U(t) f_0 \right\}(x,p)
=
G(t)f_0(x,p)
+
\int_0^t  ~ ds_1 ~ G(t-s_1 )K\left\{  U(s_1 ) f_0 \right\}(x,p).
\end{gather*}
We employ the splitting
$
K
=
K^{1-\chi}
+
K^{\chi}
$
which is defined with the cut-off function \eqref{cut} and \eqref{kCUT}.
We then further expand out
\begin{multline*}
\left\{  U(t) f_0 \right\}(x,p)
=
G(t)f_0(x,p)
+
\int_0^t  ~ ds_1  ~ G(t-s_1 )K^{1-\chi} \left\{  U(s_1 ) f_0 \right\}(x,p)
\\
+
\int_0^t  ~ ds_1  ~ G(t-s_1 )K^{\chi} \left\{  U(s_1 ) f_0 \right\}(x,p).
\end{multline*}
We  further iterate the Duhamel formula
of the last term, as did Vidav \cite{MR0259662}:
\begin{gather}
U(s_1) = G(s_1) + \int_0^{s_1} ~ ds_2 ~ G(s_1-s_2) K U(s_2). 
\label{duhamelDOUBLE} 
\end{gather}
This will grant the so-called A-Smoothing property.  Notice below that we only iterate on the $K^\chi$ term, which is different from Vidav.
Plugging this Duhamel formula into the previous expression yields a more elaborate formula
\begin{multline*}
\left\{  U(t) f_0 \right\}(x,p)
=
G(t)f_0(x,p)
+
\int_0^t  ~ ds_1 ~ G(t-s_1)K^{1-\chi}\left\{  U(s_1) f_0 \right\}(x,p)
\\
+
\int_0^t  ~ ds_1 ~G(t-s_1)K^{\chi}  G(s_1) f_0(x,p)
\\
+
\int_0^t ~ ds_1 ~  
\int_0^{s_1} ~ ds_2 ~
G(t-s_1)K^{\chi}
G(s_1-s_2)
K
\left\{  U(s_2) f_0 \right\}(x,p).
\end{multline*}
However this is not quite yet  in the form we want.  To get the final form, we once again split 
the compact operator 
$
K
=
K^{1-\chi}
+
K^{\chi}
$
in the last term to obtain
\begin{gather}
\notag
\left\{  U(t) f_0 \right\}(x,p) =
G(t)f_0(x,p)
+
\int_0^t  ~ ds_1 ~ G(t-s_1)K^{1-\chi} \left\{  U(s_1) f_0 \right\}(x,p)
\\
+
\int_0^t  ~ ds_1 ~G(t-s_1)K^{\chi}  G(s_1) f_0 (x,p)
\label{linearEXP}
\\
\notag
+
\int_0^t ~ ds_1 ~  
\int_0^{s_1} ~ ds_2 ~
G(t-s_1)K^{\chi}
G(s_1-s_2)
K^{1-\chi}
\left\{  U(s_2) f_0 \right\}(x,p)
\\
+
\int_0^t ~ ds_1 ~  
\int_0^{s_1} ~ ds_2 ~
G(t-s_1)K^{\chi}
G(s_1-s_2)
K^{\chi}
\left\{  U(s_2) f_0 \right\}
(x,p)
\notag
\\
\eqdef
H_1(t,x,p)
+
H_2(t,x,p)
+
H_3(t,x,p)
+
H_4(t,x,p)
+
H_5(t,x,p),
\notag
\end{gather}
where
\begin{eqnarray*}
H_1(t,x,p)
& \eqdef &
e^{-\nu(p) t } f_0(x-\hat{p}t, p),
\\
H_2(t,x,p)
& \eqdef &
\int_0^t  ~ds_1 ~ e^{-\nu(p) (t-s_1) }
K^{1-\chi} \left\{  U(s_1) f_0 \right\}(y_1,p),
\\
H_3(t,x,p)
& \eqdef &
\int_0^t  ~ds_1 ~ e^{-\nu(p) (t-s_1) }
\int_{\mathbb{R}^3} dq_1~ k^\chi(p,q_1)~ e^{-\nu(q_1) s_1 }
f_0(y_1-\hat{q}_1s_1,q_1).
\end{eqnarray*}
Just above and below we will be using the following short hand notation
\begin{eqnarray}
\notag
y_1 & \eqdef & x- \hat{p}(t-s_1),
\\
y_2 & \eqdef & y_1-\hat{q}_1(s_1-s_2)
=
x- \hat{p}(t-s_1)-\hat{q}_1(s_1-s_2).
\label{zDEF}
\end{eqnarray}
We are also using the notation
$
q_{10} = \sqrt{1+|q_{1}|^2}
$
and
$
q_1= (q_{11},q_{12},q_{13})\in\mathbb{R}^3
$
with
$
\hat{q}_1 = q_1/q_{10}.
$
Furthermore the next term is
\begin{multline*}
H_4
(t,x,p)
\eqdef
\int_{\mathbb{R}^3} dq_1~ k^\chi(p,q_1)
\int_0^t  ~ds_1 \int_0^{s_1}  ~ds_2 ~
e^{-\nu(p) (t-s_1) }
 e^{-\nu(q_1) (s_1- s_2) }
 \\
 \times
K^{1-\chi}
\left\{  U(s_2) f_0 \right\}
(y_2,q_1).
\end{multline*}
Lastly, we may also expand out the fifth component as
\begin{multline}
H_5(t,x,p)
=
\int_{\mathbb{R}^3} dq_1~ k^\chi(p,q_1)~ 
\int_{\mathbb{R}^3} dq_2~ k^\chi(q_1,q_2)~ 
\int_0^t  ~ds_1 ~ e^{-\nu(p) (t-s_1) }
\\
\times
\int_0^{s_1} ~ ds_2 ~
e^{-\nu(q_1) (s_1-s_2) }
\left\{  U(s_2) f_0 \right\}(y_2,q_2).
\label{h5}
\end{multline}
We will estimate each of these five terms individually.  
In Lemma \ref{hESTIMATES}  below we will show that the first and third term exhibit rapid polynomial decay.
Then after that, in Lemma \ref{hESTIMATE24}, we  show that the second and fourth terms can be bounded by the time decay norm \eqref{dNORM} multiplied by an arbitrarily small constant.
For the last term, in Lemma \ref{hESTIMATE3}, we will show that $H_5$ can be estimated by \eqref{dNORM}  times a small constant plus the $L^2_{-j}$ norm of the semi-group (for any $j>0$) multiplied by a large constant.
After stating each of these lemmas, we will put these estimates together to prove our key decay estimate on the semi-group for solutions to \eqref{rBlin} in $L^\infty_\ell$ as stated above in Theorem \ref{decay0l}.  Once this is complete we give the proofs of the three key Lemmas \ref{hESTIMATES}, \ref{hESTIMATE24} and \ref{hESTIMATE3} at the end of this section.

\begin{lemma}
\label{hESTIMATES}
Given $\ell \ge 0$,
for any  $k\ge 0$ we have
\begin{equation*}
\left|
w_\ell(p)
H_1(t,x,p)
\right|
+
\left|
w_\ell(p)
H_3(t,x,p)
\right|
\le C_{\ell,k} (1+t)^{-k}\| f_0\|_{\infty,\ell+ k}.
\end{equation*}
\end{lemma}

Next

\begin{lemma}
\label{hESTIMATE24}
Fix  $\ell \ge 0$.
For any small $\eta>0$,  which relies upon the small $\epsilon>0$ from \eqref{cut}, and any $k\ge 0$  we have
\begin{equation*}
\left|
w_\ell(p)
H_2(t,x,p)
\right|
+
\left|
w_\ell(p)
H_4(t,x,p)
\right|
\le \eta(1+t)^{-k} ||| f |||_{k,\ell}.
\end{equation*}
\end{lemma}

The estimates in the lemma above will be used to obtain upper bounds for the the $L^\infty_\ell$ norm of the semi-group.  The final lemma in this series is below

\begin{lemma}
\label{hESTIMATE3}
Fix $\ell \ge 0$, choose any (possibly large) $j>0$.
For any small $\eta>0$, which depends upon \eqref{cut}, and any $k\ge 0$ we have the estimate
\begin{multline*}
\left|
w_\ell(p)
H_5(t,x,p)
\right|
\le \eta(1+t)^{-k} ||| f |||_{k,\ell}
+
C_\eta 
\int_0^t ds ~ e^{-\eta (t-s)} \|f\|_{2,-j}(s)
\\
+
w_\ell(p)
\left|
R_1(f)(t)
\right|.
\notag
\end{multline*}
By the $L^2$ decay theory from Theorem \ref{decay2}, and also Proposition \ref{BasicDecay}, we have
$$
\int_0^t ds ~ e^{-\eta (t-s)} \|f\|_{2,-j}(s)
\le
C_\eta (1+t)^{-k}\|f_0 \|_{2,k}
\le
C_\eta (1+t)^{-k}\|f_0 \|_{\infty,\ell+k}.
$$
The above estimates hold for any $k\ge 0$ and $\ell> 3/b$ (as in \eqref{ELLestINEQ} just below).  On the other hand, for the last term
 if we restrict $k\in[0,1]$ then $\forall \eta>0$ we have
$$
w_\ell(p)
\left|
R_1(f)
\right|
\le
\eta (1+t)^{-k}||| f |||_{k,\ell}.
$$
Above $R_1$ is defined in \eqref{r1def} during the course of the proof.
\end{lemma}

These estimates would imply almost exponential decay except for the problematic term $R_1(f)(t)$, which only appears to decay to first order.  This will be discussed in more detail in Section \ref{s:linRdecay}, where it is shown that this term can decay to any polynomial order by performing a new high order expansion.

We  now show that the above lemmas grant a uniform bound and slow decay for solutions to \eqref{rBlin} in $L^\infty_\ell$.  We are using the semi-group notation
$
f(t) = \{U(t)f_0\}.
$
Lemmas \ref{hESTIMATES}, \ref{hESTIMATE24}, and \ref{hESTIMATE3} together imply that for any $\eta >0$ and $k\in [0,1]$ we have
$$
\|f \|_{\infty,\ell}(t)
\le
C_{\ell,k} (1+t)^{-k}\| f_0\|_{\infty,\ell+ k}
+
\eta (1+t)^{-k} ||| f |||_{k,\ell}
+
C_\eta (1+t)^{-k}\|f_0 \|_{2,k}.
$$
Equivalently
$$
 ||| f |||_{k,\ell}
\le
C_{\ell,k} \| f_0\|_{\infty,\ell+ k}
+
C_{1/2} \|f_0 \|_{2,k}
\le
C_{\ell,k} \| f_0\|_{\infty,\ell+ k}.
$$
The last estimate holds when we choose $\ell>3/b$, with \eqref{weight}, as follows
\begin{multline}
\|f_0 \|_{2,k} = \sqrt{\int_{\mathbb{R}^3} dp ~ w_{k}^2(p) |f_0(p)|^2}
\le  
\|f_0 \|_{\infty,\ell+k} \sqrt{\int_{\mathbb{R}^3} dp ~ w_{-\ell}^2(p) }
\\
\lesssim
\|f_0 \|_{\infty,\ell+k}.  
\label{ELLestINEQ}
\end{multline}
With this inequality, we have
 the desired decay rate for the $L_\ell^\infty$ norm of solutions to the linear equation \eqref{rBlin},
which proves Theorem \ref{decay0l} subject to  Lemmas \ref{hESTIMATES}, \ref{hESTIMATE24}, and \ref{hESTIMATE3}.  We now prove those lemmas.

Along this course we  will repeatedly use the following basic decay estimate

\begin{proposition}
\label{BasicDecay}
Suppose without loss of generality that $\lambda \ge \mu \ge 0$.  Then
$$
\int_0^t  \frac{ds}{(1+ t - s)^{\lambda} (1+ s)^{\mu}}
\le
\frac{C_{\lambda,\mu}(t)}{ (1+t)^{\rho}},
$$
where $\rho = \rho(\lambda,\mu) =\min\{ \lambda+\mu -1, \mu\}$ and
$$
C_{\lambda,\mu}(t) = 
C
\left\{
\begin{array}{cc}
1 & \mbox{if} ~ \lambda \ne 1,
\\
\log(2+t) & \mbox{if} ~ \lambda = 1.
\end{array}
\right.
$$
\end{proposition}

Furthermore, we will use the following 
basic estimate from the Calculus
\begin{equation}
e^{-ay} (1+y)^k \le \max\{1, e^{a-k} k^k a^{-k}\},
\quad
a,y, k \ge 0.
\label{ElemCalc}
\end{equation}
We will now write an elementary proof of this basic time decay estimate in Proposition \ref{BasicDecay}.  This result is not difficult and known, however 
we provide a short proof 
for the sake of completeness and because 
we have not seen a proof in the literature.  \\

\noindent {\it Proof of Proposition \ref{BasicDecay}}.  We will consider the cases $\mu =0$ and $\mu = \lambda$ separately.  Then the general result will then be established by interpolation.

{\bf Case 1: $\mu = 0$}.  If $\lambda \ne 1$ we have
\begin{gather*}
\int_0^t  \frac{ds}{(1+ t - s)^{\lambda}}
=
\int_0^t  \frac{ds}{(1+  s)^{\lambda}}
=
\frac{1}{\lambda - 1}\left\{1- (1+  t)^{1-\lambda} \right\}
\le
C (1+t)^{-\rho(\lambda,0)}.
\end{gather*}
Note $\rho = 0$ if $\lambda > 1$ and $\rho = \lambda -1<0$ otherwise.
Alternatively if $\lambda =1$
\begin{gather*}
\int_0^t  \frac{ds}{1+  s}
=
\log (1+  t).
\end{gather*}
This completes our study of the first case $\mu = 0$.

{\bf Case 2: $\mu = \lambda$}.  We split the integral as
\begin{gather*}
\int_0^t  \frac{ds}{(1+ t - s)^{\lambda} (1+  s)^{\lambda}}
=
\int_0^{t/2}  
+
\int_{t/2}^t.
\end{gather*}
For the first integral
\begin{gather*}
\int_0^{t/2}  \frac{ds}{(1+ t - s)^{\lambda} (1+  s)^{\lambda}}
\le
(1+ t/2)^{-\lambda} \int_0^{t/2}  \frac{ds}{ (1+  s)^{\lambda}}.
\end{gather*}
Now from Case 1, we can estimate the remaining integral as
\begin{gather*}
 \int_0^{t/2}  \frac{ds}{ (1+  s)^{\lambda}}
 \le
C_{\lambda,\lambda}(t) (1+t)^{\max\{0,1-\lambda \}},
\end{gather*}
which conforms with the claimed decay.
The second half of the integral can be estimated in exactly the same way as the first.

{\bf Case 3: $0<\mu < \lambda$}.  By H{\"o}lder's inequality, we have
\begin{multline*}
\int_0^t  \frac{ds}{(1+ t - s)^{\lambda} (1+ s)^{\mu}}
=
\int_0^t  \frac{ds}{(1+ t - s)^{\frac{\lambda}{p}+\frac{\lambda}{q}} (1+ s)^{\mu}}
\\
\le
\left( 
\int_0^t  \frac{ds}{(1+ t - s)^{\lambda}}
\right)^{1/p}
\left( 
\int_0^t  \frac{ds}{(1+ t - s)^{\lambda} (1+ s)^{q\mu}}
\right)^{1/q}.
\end{multline*}
Above 
$
q=\frac{\lambda}{\mu}
$
and
$
p=\frac{\lambda}{\lambda-\mu}
$
with
$
\frac{1}{p}+\frac{1}{q}
=1.
$
Therefore the above is
\begin{gather*}
\le
C
\left( 
\int_0^t  \frac{ds}{(1+ t - s)^{\lambda}}
\right)^{(\lambda - \mu)/\lambda}
\left( 
\int_0^t  \frac{ds}{(1+ t - s)^{\lambda} (1+ s)^{\lambda}}
\right)^{\mu/\lambda}.
\end{gather*}
By the previous two cases, this is 
\begin{multline*}
\le
\left( 
C_{\lambda,0}(t) (1+ t )^{-\rho({\lambda,0})}
\right)^{(\lambda - \mu)/\lambda}
\left( 
C_{\lambda,\lambda}(t)
(1+ t )^{-\rho({\lambda,\lambda})}
\right)^{\mu/\lambda}
\\
=
C_{\lambda,0}(t)^{(\lambda - \mu)/\lambda}
C_{\lambda,\lambda}(t)^{\mu/\lambda}
 (1+ t )^{-\rho({\lambda,0})(\lambda - \mu)/\lambda-\rho({\lambda,\lambda})\mu/\lambda}.
\end{multline*}
The proposition  follows by adding the exponents.
\qed \\

We are now ready to proceed to the \\

\noindent {\it Proof of Lemma \ref{hESTIMATES}.}
  We start with $H_1$.  
From Lemma \ref{nuEST} and
\eqref{ElemCalc}
we have
\begin{equation}
e^{-\nu(p) t }
\le
C_{k} 
p_0^{kb/2}
(1+t)^{-k}
\le
C_{k} 
w_k(p)
(1+t)^{-k},
\quad
\forall t, k >0.
\label{polyE}
\end{equation}
Here we use the notation from \eqref{weight}.
This procedure grants high polynomial time decay on the solution if we admit similar high polynomial momentum decay on the initial data.
In particular we have shown
$$
\left|
w_\ell(p) H_1(t,x,p)
\right|
=
\left|
w_\ell(p) ~ e^{-\nu(p) t } f_0(x-\hat{p}t, p)
\right|
\le
C (1+t)^{-k}
\| f_0 \|_{\infty,\ell+ k},
$$
which is the desired estimate for $H_1$.

We finish off this lemma by estimating $H_3$.  Notice that we trivially have
$$
 e^{-\nu(p) (t-s_1) }
  e^{-\nu(q) s_1 }
  \le
  e^{-\nu(\max\{|p|,|q|\}) t },
  $$
  where 
  $
  \nu(\max\{|p|,|q|\})
$
is $\nu$ evaluated at 
$
\max\{|p|,|q|\}.
$
We have
\begin{multline*}
\left| 
w_\ell (p)
H_3(t,x,p)
\right|
\le
w_\ell (p)
\int_0^t  ~ds_1 ~ e^{-\nu(p) t }
\int_{|p| \ge |q_1|} dq_1 \left| k^\chi(p,q_1) \right| 
\sup_{y\in\mathbb{T}^3}
\left| 
f_0(y,q_1)
\right|
\\
+
w_\ell (p)
\int_0^t  ~ds_1 ~ 
\int_{|p| < |q_1|} dq_1~ \left| k^\chi(p,q_1) \right| ~ e^{-\nu(q_1) t }
\sup_{y\in\mathbb{T}^3}
\left| 
f_0(y,q_1)
\right|.
\end{multline*}
We will estimate the second term, and we  remark that the first term can be handled in exactly the same way.
As in the previous estimate for $H_1$, we use  \eqref{polyE}  to obtain
$$
e^{-\nu(q_1) t }
\sup_{y\in\mathbb{T}^3} 
\left| 
f_0(y,q_1)
\right|
\le
C_k(1+t)^{-k-1} w_{k+1}(q_1)\sup_{y\in\mathbb{T}^3} \left| f_0(y,q_1) \right|.
$$
Next we use the estimate for $k^\chi(p,q)$ from Lemma \ref{boundK2}.
When using this estimate
 we may suppose
$
|q_1| \le 2|p|.
$
For otherwise, if say $|q_1| \ge 2 |p|$, then as in \eqref{split1314}
we have 
$
|p-q_1| \ge |q_1| /2
$
which leads directly to
\begin{gather}
w_{k+1}(q_1) w_\ell (p)
e^{-c|p-q_1|/2}
\le
Cw_{k+1+\ell}(q_1) 
e^{-c|q_1|/4}
\le C.
\label{comparablePQ}
\end{gather}
In this case we easily obtain an estimate better than \eqref{desired1} below.   In particular 
\begin{multline*}
\int_0^t  ~ds_1 ~ 
w_\ell (p)
\int_{|p| < |q_1|} dq_1~{\bf 1}_{|q_1| \ge 2|p|}~ \left| k^\chi(p,q_1) \right| ~ e^{-\nu(q_1) t }
\sup_{y\in\mathbb{T}^3}\left| f_0(y,q_1)
\right|
\\
\lesssim
(1+t)^{-k} \| f_0 \|_{\infty,-j},
\quad \forall k>0, ~ j >0.
\end{multline*}
Thus in the following  we assume 
$
 |p| < |q_1| \le 2 |p|.
$
On this region we may plug in the last few estimates including \eqref{polyE} to obtain
\begin{multline}
\int_0^t  ~ds_1 ~ 
w_\ell (p)
\int_{|p| < |q_1|\le 2|p|} dq_1~ \left| k^\chi(p,q_1) \right|  ~ e^{-\nu(q_1) t }
\sup_{y\in\mathbb{T}^3}
\left| f_0(y,q_1) \right|
\\
\le
C_{k} \frac{t}{(1+t)^{k+1}}
\int_{|p| < |q_1|\le 2|p|} dq_1~
w_\ell (q_1) ~ w_{k+1 }(q_1) ~ \left| k^\chi(p,q_1) \right| ~
\sup_{y\in\mathbb{T}^3} \left| f_0(y,q_1) \right|
\\
\le
C_{k}(1+t)^{-k} \| f_0\|_{\infty, \ell+ k}
~ w_{1 }(p)
\int_{|p| < |q_1|\le 2|p|} dq_1~ \left| k^\chi(p,q_1) \right|.
\label{desired1}
\end{multline}
Then from Lemma \ref{boundK2} we clearly have the following bound 
$$
w_{1 }(p)
\int_{|p| < |q_1|\le 2|p|} dq_1~ \left| k^\chi(p,q_1) \right|
\le C.
$$
This completes the time decay estimate for $H_3$ and our proof of  the lemma.
\qed \\

Our next aim is to prove Lemma \ref{hESTIMATE24}.  To do this we will use the following

\begin{lemma}
\label{boundKinfX}  
Fix any $\ell \ge 0$ and any $j>0$.
Then given any small $\eta >0$, which depends upon $\chi$ in \eqref{cut}, the following estimate holds
$$
\left| w_\ell(p) K^{1-\chi}(h)(p) \right| 
\le
\eta e^{- cp_0} \| h\|_{\infty,-j}.
$$
Above the constant $c>0$ is independent of $\eta$.
\end{lemma}

\noindent {\it Proof of Lemma \ref{boundKinfX}.}
We consider $K^{1-\chi}(h)$ as defined in \eqref{kCUT}.
From \eqref{0primeBELOW},
for $\epsilon$ in \eqref{cut} chosen sufficiently small ($\epsilon$ smaller than $1/4$ is sufficient),
 we see that
 $$
\sqrt{ J(q) J(p')} + \sqrt{ J(q) J(q')}
\lesssim
e^{-c q_0 - c p_0},
\quad 
\exists c>0.
 $$
 With that estimate,
and additionally the conservation of energy  \eqref{collisionalCONSERVATION}, we have
\begin{multline}
\left| w_\ell(p) K_2^{1-\chi}(h)(p) \right| 
\lesssim
\int_{\mathbb{R}^3\times \mathbb{S}^{2}} ~ d\omega dq 
~ \left( 1 - \chi(g) \right)
~ v_{\o} ~ \sigma(g,\theta)
e^{-c q_0 - c p_0}
\notag
\\
\qquad
\times \left( \left|  h(p^{\prime }) \right| + \left|  h(q^{\prime }) \right| \right)
\notag
\\
\lesssim
\int_{\mathbb{R}^3\times \mathbb{S}^{2}} ~ d\omega dq 
~ \left( 1 - \chi(g) \right)
~ v_{\o} ~ \sigma ~ 
e^{-\frac{c}{2} q_0 - \frac{c}{2} p_0}
 \left(
  e^{ - \frac{c}{2} p_0^{\prime }}\left|  h(p^{\prime}) \right| 
+ 
e^{-\frac{c}{2} q_0^{\prime }}
\left|  h(q^{\prime }) \right| \right)
\notag
\\
\le 
\eta
e^{ - \frac{c}{4}  p_0}
\|  h \|_{\infty,-j}.
\notag
\end{multline}
In the last inequality we have used the following estimate for any small $\eta >0$
\begin{equation}
\int_{\mathbb{R}^3\times \mathbb{S}^{2}} ~ d\omega dq 
~ \left( 1 - \chi(g) \right)
~ v_{\o} ~ \sigma (g, \theta) ~ 
e^{-\frac{c}{2} q_0 - \frac{c}{2} p_0}
\le \eta ~ e^{ - \frac{c}{4} p_0}.
\label{keyEPbound}
\end{equation}
The proof of this bound uses \eqref{generalK}, but also
the exact strategy used in the paragraph containing \eqref{generalK} and the paragraph just below it.  The  idea is to use the splitting \eqref{smmothPr} inside the integral in \eqref{keyEPbound}:
$$
1 = \phi_R(p,q) + (1-\phi_R(p,q) ).
$$
First when $|p|+|q|$ is large, this is the term containing $(1-\phi_R(p,q) )$, we have a bound for the integral in \eqref{keyEPbound} which is of the form $C e^{-cR}$ where the constant $C, c >0$ are independent of both $R$ and $\epsilon$ similar to \eqref{RboundINDEP}.  On the other hand, when $|p|+|q|\le R$, this is the term containing $\phi_R(p,q)$, we have a bound which is of the form $C_R \epsilon^p$.  
Here  $p = p(b)$ can be chosen to be $p=3$ if $b\in [0,1]$ and $p=4-b$ if $b\in(1,4)$.  
This second estimate is similar to \eqref{keyRbound}.
Since this strategy is already preformed in detail nearby \eqref{generalK}, we will not re-write the details.

Since $p>0$ we can first choose $R>>1$ sufficiently large, and then choose $\epsilon>0$ sufficiently small so that the constant $\eta>0$ in 
\eqref{keyEPbound}
can be chosen arbitrarily small.
 This yields the desired estimate.  We remark that the estimate for 
$
K_1^{1-\chi}
$
can be shown in the same exact way, it is in fact slightly easier.
\qed   \\

With Lemma \ref{boundKinfX} and Proposition \ref{BasicDecay} in hand, we proceed to the \\

\noindent {\it Proof of Lemma \ref{hESTIMATE24}.}
We begin with the estimate for $H_2$.  Using Lemma \ref{boundKinfX}, for any small $\eta'>0$ we have
\begin{multline*}
\left|
w_\ell(p) H_2(t,x,p) 
\right|
=
w_\ell(p)\left| 
\int_0^t  ~ds_1 ~ e^{-\nu(p) (t-s_1) }
K^{1-\chi} \left\{  U(s_1) f_0 \right\}(y_1,p)
\right|
\\
\le
\eta'  ~ e^{-c p_0}
 \int_0^t  ~ds_1 ~ e^{-\nu(p) (t-s_1) } ~ \|\left\{  U(s_1) f_0 \right\}\|_{\infty,-j},
 \quad \forall j>0
 \\
\le
\eta'  ~ e^{-c p_0} ||| f |||_{k,\ell}
 \int_0^t  ~ds_1 ~ e^{-\nu(p) (t-s_1) } (1+s_1)^{-k}.
\end{multline*}
The norm  is from \eqref{dNORM} for $k\ge 0$.
As in \eqref{polyE}
for any $\lambda >\max\{1,k\}$ 
we have
\begin{multline*}
\left|
w_\ell(p) H_2(t,x,p) 
\right|
\le
\eta'  ~ w_{\lambda}(p) e^{-c p_0} ||| f |||_{k,\ell}
 \int_0^t  ~ds ~ (1+t-s)^{-\lambda} (1+s)^{-k}
\\
\le
\eta ~ (1+t)^{-k}
||| f |||_{k,\ell},
\end{multline*}
which follows from Proposition \ref{BasicDecay}.
This is the desired estimate for $H_2$.

For $H_4$ we once again use Lemma \ref{boundKinfX}, for any small $\eta'>0$, to obtain
\begin{multline*}
\left|
w_\ell(p)
H_4(t,x,p)
\right|
\le
\int_{\mathbb{R}^3} dq_1~ k^\chi(p,q_1)
\int_0^t  ~ds_1 \int_0^{s_1}  ~ds_2 ~
e^{-\nu(p) (t-s_1) }
\\
\times
 e^{-\nu(q_1) (s_1- s_2) }
w_\ell(p)
\left|
K^{1-\chi}\left(
\left\{  U(s_2) f_0 \right\}
\right)
(y_2,q_1)
\right|
\\
\le
\eta' ||| f |||_{k,\ell}
\int_{\mathbb{R}^3} dq_1~ k^\chi(p,q_1)
\frac{w_\ell(p)}{w_\ell(q_1)}
e^{-cq_{10}}
\\
\times
 \int_0^t  ~ds_1 ~ 
\int_0^{s_1}  ~d s_2 ~ e^{-\nu(p) (t-s_1) } e^{-\nu(q_1) (s_1- s_2) }(1+s_2)^{-k}.
\end{multline*}
We recall that $q_{10} =\sqrt{1+|q_1|^2}$.
For the time decay, from Proposition \ref{BasicDecay} with
 \eqref{ElemCalc}
 as in
\eqref{polyE}
we notice that
\begin{multline}
 \int_0^t  ~ds_1 ~ 
\int_0^{s_1}  ~d s_2 ~ e^{-\nu(p) (t-s_1) } e^{-\nu(q_1) (s_1- s_2) }(1+s_2)^{-k}
\\
\lesssim
w_{ \lambda}(p) w_{\lambda}(q_1)
\int_0^t  ~ds_1 ~ (1+t-s_1)^{-\lambda}
\int_0^{s_1}  ~ds_2 ~ (1+s_1-s_2)^{-\lambda}(1+s_2)^{-k}
\\
\lesssim
 w_{\lambda}(p) w_{\lambda}(q_1)
(1+t)^{-k}.
\notag
\end{multline}
Above we have taken $\lambda >\max\{1,k\}$.
Combining these estimates yields
\begin{gather*}
\left|
w_\ell(p)
H_4 
\right|
\lesssim
 \eta' (1+t)^{-k} ||| f |||_{k,\ell} 
\int_{\mathbb{R}^3} dq_1~ k^\chi(p,q_1)
w_{\ell + \lambda}(p) w_{-\ell + \lambda}(q_1)
e^{-cq_{10}}.
\end{gather*}
To estimate the remaining integral and weights
we split into three cases.
If either
$
2|q_1| \le |p|,
$
or
$
|q_1| \ge 2|p|,
$
then we bound all the weights and the remaining momentum integral by a constant as in 
\eqref{comparablePQ}.
Alternatively if
$
\frac{1}{2} |q_1|
\le
|p| 
\le 
2|q_1|,
$
then the desired estimate is obvious since we have strong exponential decay in both $p$ and $q_1$.  
In either of these  cases we have the estimate for $H_4$. 
\qed \\

We will finish this section with a proof of the crucial Lemma \ref{hESTIMATE3}.\\

\noindent {\it Proof of Lemma \ref{hESTIMATE3}.}
We now turn to the proof of our estimate for $H_5$.  Recall the definition of $H_5$ from \eqref{h5}
with
$y_2$
 defined in \eqref{zDEF}.
 We will utilize rather extensively the estimate for $k^\chi$ from Lemma \ref{boundK2}.
We now further split 
\begin{equation}
H_5(t,x,p) =  H_5^{high}(t,x,p)+ H_5^{low}(t,x,p),
\label{h5sum}
\end{equation}
and estimate each term on the right individually.    For $M>>1$ 
we define
\begin{equation}
{\bf 1}_{high} 
\eqdef
{\bf 1}_{|p|> M}  {\bf 1}_{|q_1|\le M} 
+ 
 {\bf 1}_{|q_1|> M}.
 \label{1high}
\end{equation}
Notice
$
{\bf 1}_{high} +{\bf 1}_{|p|\le M}  {\bf 1}_{|q_1|\le M} 
=1.
$
Now the first term in the expansion is 
\begin{gather*}
H_5^{high}(t,x,p)
\eqdef
\int_{\mathbb{R}^3}  dq_1 ~ k^\chi(p,q_1)~ 
\int_{\mathbb{R}^3} dq_2~ k^\chi(q_1,q_2)~  {\bf 1}_{high} ~
\int_0^t  ~ds_1 ~ e^{-\nu(p) (t-s_1) }
\\
\times
\int_0^{s_1} ~ ds_2 ~
e^{-\nu(q_1) (s_1- s_2) }
\left\{  U(s_2) f_0 \right\}(y_2,q_2).
\end{gather*}
We use \eqref{ElemCalc}, as in \eqref{polyE}, and Lemma \ref{nuEST} to see that for any $\lambda \ge 0$ we have
\begin{multline*}
\int_0^t  ~ds_1 ~ e^{-\nu(p) (t-s_1) }
\int_0^{s_1} ~ ds_2 ~
e^{-\nu(q_1) (s_1- s_2) }
\\
\le
C_\lambda w_\lambda (p) w_\lambda (q_1)
\int_0^t  ~ds_1 ~ 
\int_0^{s_1} ~ ds_2 ~
(1+(t-s_1))^{-\lambda} (1+(s_1 - s_2))^{-\lambda}.
\end{multline*}
When either $|p|> M$ or $|q_1|> M$, by Lemma \ref{boundK2}, we have the bound
$$
\left| k^{\chi}(p,q_1)  \right|
\le C M^{-\zeta} 
\left( p_0+ q_{10} \right)^{-b/2}
e^{-c |p-q_1|}.
$$
If either $|p| \ge 2|q_1|$ or $|q_1| \ge 2|p|$
 then as in \eqref{comparablePQ} we have
$$
w_{\ell+\lambda}(p) w_\lambda(q_1) e^{-c|p-q_1|} 
\le C.
$$
Thus by combining the last few estimates we have
\begin{multline*}
w_\ell(p) \int_{\mathbb{R}^3}  dq_1 ~ \left| k^{\chi}(p,q_1)  \right|~ 
\int_{\mathbb{R}^3} dq_2~  \left| k^{\chi}(q_1,q_2)  \right|~ {\bf 1}_{high} ~
\left( {\bf 1}_{|p| \ge 2|q_1|}+{\bf 1}_{|p| \le \frac{1}{2}|q_1|} \right)
\\
\times
\int_0^t  ~ds_1 ~ e^{-\nu(p) (t-s_1) }
\int_0^{s_1} ~ ds_2 ~
e^{-\nu(q_1) (s_1- s_2) }
\left|
\left\{  U(s_2) f_0 \right\}(y_2,q_2)
\right|
\\
\le
\frac{C_{\lambda}}{M^{\zeta+b/2}}||| f|||_{k,\ell} \int_{\mathbb{R}^3}  dq_1 ~ e^{-c|p-q_1|}~ 
\int_{\mathbb{R}^3} dq_2~ e^{-c|q_1-q_2|}~ 
\\
\times
\int_0^t  ~ds_1 ~ 
\int_0^{s_1} ~ ds_2 ~
\frac{(1+s_2)^{-k}}{(1+(t-s_1))^{\lambda} (1+(s_1 - s_2))^{\lambda} }.
\end{multline*}
With Proposition \ref{BasicDecay}, for any $k\ge 0$ and $\lambda> \max\{k,1\}$ the previous term is bounded from above by
  \begin{gather*}
\le
\frac{C_{k,\lambda}}{M^{\zeta+b/2}}~ (1+t)^{-k} ~ ||| f|||_{k,\ell}.
\end{gather*}
This is the desired estimate for $M>>1$ chosen sufficiently large.

We now consider the remaining part of $H_5^{high}$.
As in the previous estimates and \eqref{comparablePQ}, if either 
 $|q_2 | \ge 2|q_1|$ or
 $|q_1| \ge 2|q_2|$ then for any $k\ge 0$ we have
 \begin{multline*}
w_\ell(p) \int_{\mathbb{R}^3}  dq_1 ~ \left| k^{\chi}(p,q_1)  \right|~ 
\int_{\mathbb{R}^3} dq_2~  \left| k^{\chi}(q_1,q_2)  \right|~
 \left( {\bf 1}_{|q_2| \ge 2|q_1|}+{\bf 1}_{|q_2| \le \frac{1}{2}|q_1|} \right)
\\
\times
{\bf 1}_{\frac{1}{2}|p| \le |q_1|\le 2|p|}
{\bf 1}_{high} 
\int_0^t  ~ds_1 ~ e^{-\nu(p) (t-s_1) }
\int_0^{s_1} ~ ds_2 ~
e^{-\nu(q_1) (s_1- s_2) }
\\
\times
\left| 
\left\{  U(s_2) f_0 \right\}(y_2,q_2)
\right|
\le
\frac{C_{k,\lambda}}{M^{2\zeta+b}}||| f|||_{k,\ell} \int_{\mathbb{R}^3}  dq_1 ~ e^{-c|p-q_1|}~ 
\int_{\mathbb{R}^3} dq_2~ e^{-c|q_1-q_2|}~ 
\\
\times
\int_0^t  ~ds_1 ~ 
\int_0^{s_1} ~ ds_2 ~
\frac{(1+s_2)^{-k}}{(1+(t-s_1))^{\lambda} (1+(s_1 - s_2))^{\lambda} }
\\
\le
\frac{C_{k,\lambda}}{M^{2\zeta+b/2}}~ (1+t)^{-k} ~ ||| f|||_{k,\ell}.
\end{multline*}
Above we have used exactly the same estimates as in the prior case.  Both of the last two terms have a suitably small constant in front if $M$ is  sufficiently large.

Thus the remaining part of $H_5^{high}$ to estimate is 
\begin{gather}
\notag
R_1(f)(t)\eqdef \int_{\mathbb{R}^3}  dq_1 ~ k^\chi(p,q_1)~ 
\int_{\mathbb{R}^3} dq_2~ k^\chi(q_1,q_2)~ 
{\bf 1}_{\frac{1}{2}|p| \le |q_1|\le 2|p|} {\bf 1}_{\frac{1}{2}|q_1| \le |q_2|\le 2|q_1|} 
\\
\times
{\bf 1}_{high} ~
\int_0^t  ~ds_1 ~ e^{-\nu(p) (t-s_1) }
\int_{0}^{s_1} ~ ds_2 ~
e^{-\nu(q_1) (s_1- s_2) }
\left\{  U(s_2) f_0 \right\}(y_2,q_2).
\label{r1def}
\end{gather}
It is only this term which slows down the time decay rate. In this current proof we will only make a basic argument to show that this term can naively exhibit first order decay.  Since all the momentum variables are comparable, we have
\begin{multline}
\notag
\left| R_1(f)(t) \right| 
\le 
\int_{\mathbb{R}^3}  dq_1 ~ \left| k^{\chi}(p,q_1)  \right|~ 
\int_{\mathbb{R}^3} dq_2~  \left| k^{\chi}(q_1,q_2)  \right|~
\\
\times
{\bf 1}_{\frac{1}{2}|p| \le |q_1|\le 2|p|} {\bf 1}_{\frac{1}{2}|q_1| \le |q_2|\le 2|q_1|}
{\bf 1}_{high}
\\
\times
\int_0^t  ~ds_1 ~ e^{-c\nu(q_1) (t-s_1) }
\int_0^{s_1} ~ ds_2 ~
e^{-c\nu(q_1) (s_1-s_2) }
\left| 
\left\{  U(s_2) f_0 \right\}(y_2,q_2)
\right|.
\notag
\end{multline}
Next using similar techniques as in the previous two estimates, including Proposition \ref{BasicDecay} twice, we obtain the following upper bound for any $k\in[0,1]$:
\begin{multline*}
w_\ell(p)\left| R_1(f)(t) \right| 
\le 
\frac{C_{k,\lambda}}{M^{2\zeta}}||| f|||_{k,\ell}  \int  dq_1 ~ q_{10}^{-b/2 - \zeta}e^{-c|p-q_1|}~ 
\\
\times
\int dq_2~  q_{10}^{-b/2 - \zeta}e^{-c|q_1-q_2|}~
w_{2+2\delta}(q_1)
\int_0^{t} ~ \frac{ ds_1 }{ (1+(t - s_1))^{1+\delta} }
\\
\times
\int_0^{s_1} ~ 
\frac{ ds_2 }{ (1+(s_1 - s_2))^{1+\delta} (1+s_2)^{k}}
\le 
\frac{C_{k,\lambda}}{M^{2\zeta}}
(1+t)^{-k}
||| f|||_{k,\ell}.
\end{multline*}
In the last line we used the fact that we have chosen $\delta >0$ to satisfy $\delta < \zeta$ where $\zeta>0$ is defined in the statement of Lemma \ref{boundK2}.  In Section \ref{s:linRdecay} we will examine this term at length to show that 
\eqref{r1def} actually decays ``almost exponentially.''

We are ready to define the second term in our splitting of $H_5$.  It must be 
\begin{multline*}
H_5^{low}(t,x,p)
\eqdef
{\bf 1}_{|p|\le M} 
\int_{|q_1|\le M}  dq_1 ~ k^\chi(p,q_1)~ 
\int_{\mathbb{R}^3} dq_2~ k^\chi(q_1,q_2)~ 
\\
\times
\int_0^t  ~ds_1 ~ e^{-\nu(p) (t-s_1) }
\int_0^{s_1} ~ ds_2 ~
e^{-\nu(q_1) (s_1- s_2) }
\left\{  U(s_2) f_0 \right\}(y_2,q_2).
\end{multline*}
For any small $\kappa>0$, we further split this term into two terms, one of which is
\begin{multline*}
H_5^{low,\kappa}(t,x,p)
\eqdef
\int_{|q_1|\le M}  dq_1 ~ k^\chi(p,q_1)~ 
\int_{\mathbb{R}^3} dq_2~ k^\chi(q_1,q_2)~ 
\int_0^{\kappa}  ~ds_1 ~ e^{-\nu(p) (t-s_1) }
\\
\times
{\bf 1}_{|p|\le M} 
\int_0^{s_1} ~ ds_2 ~
e^{-\nu(q_1) (s_1- s_2) }
\left\{  U(s_2) f_0 \right\}(y_2,q_2)
\\
+
\int_{|q_1|\le M}  dq_1 ~ k^\chi(p,q_1)~ 
\int_{\mathbb{R}^3} dq_2~ k^\chi(q_1,q_2)~ 
\int_{\kappa}^t  ~ds_1 ~ e^{-\nu(p) (t-s_1) }
\\
\times
{\bf 1}_{|p|\le M} 
\int_{s_1-\kappa}^{s_1} ~ ds_2 ~
e^{-\nu(q_1) (s_1- s_2) }
\left\{  U(s_2) f_0 \right\}(y_2,q_2).
\end{multline*}
The other term in this latest splitting is defined just below as $H_5^{low,2}$.
On this temporal integration domain $s_1- s_2 \le \kappa$.  Since we are proving uniform bounds, it is safe to assume when proving decay that $t\ge 1$ for instance.  

Since $p$ and $q_1$ are both bounded by $M$, from Lemma \ref{nuEST} 
we have
\begin{gather}
{\bf 1}_{|p|\le M}  {\bf 1}_{|q_1|\le M}  ~
 e^{-\nu(p) (t-s_1) -\nu(q_1) (s_1- s_2) }
\le
e^{-C(t- s_2)/M^{b/2} }.
\label{MdecayEST}
\end{gather}
Then for the first term in $H_5^{low,\kappa}$ above multiplied by $w_\ell(p)$  we have the bound
\begin{multline*}
w_\ell(p)\int_{|q_1|\le M} dq_1 ~ \left| k^{\chi}(p,q_1)  \right|~ 
\int_{\mathbb{R}^3} dq_2~  \left| k^{\chi}(q_1,q_2)  \right|~
\int_0^{\kappa}  ~ds_1 ~ e^{-\nu(p) (t-s_1) }
\\
\times
{\bf 1}_{|p|\le M} 
\int_0^{s_1} ~ ds_2 ~
e^{-\nu(q_1) (s_1- s_2) }
\left|\left\{  U(s_2) f_0 \right\}(y_2,q_2)\right|
\\
\le
C_M ||| f |||_{0,\ell} 
\int_0^{\kappa}  ~ds_1 ~\int_0^{s_1} ~ ds_2 ~
 e^{-C(t- s_2)/M^{b/2} }
 \\
\le
C_M\kappa^2  ||| f |||_{0,\ell}  ~ e^{-Ct/M^{b/2} } e^{C\kappa /M^{b/2} }
\le
C_M \kappa^2  (1+t)^{-k}  ||| f |||_{0,\ell}.
\end{multline*}
We have just used \eqref{ElemCalc}.  We obtain the desired estimate for the above terms by first choosing $M$ large, and second choosing $\kappa=\kappa(M)>0$ sufficiently small.

For the second term in $H_5^{low,\kappa}$ multiplied by $w_\ell(p)$  for any $k\ge 0$ we have
\begin{multline*}
w_\ell(p)
\int_{|q_1|\le M} dq_1 ~ \left| k^{\chi}(p,q_1)  \right|~ 
\int_{\mathbb{R}^3} dq_2~  \left| k^{\chi}(q_1,q_2)  \right|~
\int_{\kappa}^t  ~ds_1 ~ e^{-\nu(p) (t-s_1) }
\\
\times
{\bf 1}_{|p|\le M} 
\int_{s_1-\kappa}^{s_1} ~ ds_2 ~
e^{-\nu(q_1) (s_1- s_2) }
\left|\left\{  U(s_2) f_0 \right\}(y_2,q_2)\right|
\\
\le
C_M ||| f |||_{k,\ell}
\int_{\kappa}^t  ~ds_1 ~ \int_{s_1-\kappa}^{s_1} ~ ds_2 ~
 e^{-C(t- s_1)/M^{b/2} }e^{-C(s_1- s_2)/M^{b/2} }(1+s_2)^{-k}.
\end{multline*}
Since $s_2\in [s_1 - \kappa, s_1]$ and $\kappa\in(0,1/2)$, then 
$
(1+s_2) \ge 
\left(\frac{1}{2}+s_1\right).
$
We have
\begin{gather*}
\le
C_M \kappa ||| f |||_{k,\ell}
\int_{\kappa}^t  ~ds_1 ~ 
 e^{-C(t- s_1)/M^{b/2} }(1+s_1)^{-k}
 \\
 \le
C_M \kappa (1+t)^{-k} ||| f |||_{k,\ell}.
\end{gather*}
In the last step we have used \eqref{ElemCalc} and Proposition \ref{BasicDecay}.  
We conclude the desired estimate for $H_5^{low,\kappa}$ by first choosing $M$ large, and  then $\kappa>0$ sufficiently small.

The only remaining part of 
$
H_5^{low}(t,x,p)
$
to be estimated  is given by
\begin{multline*}
H_5^{low,2}(t,x,p)
\eqdef
\int_{|q_1|\le M}  dq_1 ~ k^\chi(p,q_1)~ 
\int_{\mathbb{R}^3} dq_2~ k^\chi(q_1,q_2)~ 
\int_{\kappa}^t  ~ds_1 ~ e^{-\nu(p) (t-s_1) }
\\
\times
{\bf 1}_{|p|\le M} 
\int_0^{s_1-\kappa} ~ ds_2 ~
e^{-\nu(q_1) (s_1- s_2) }
\left\{  U(s_2) f_0 \right\}(y_2,q_2).
\end{multline*}
With all of  the  $H_5^{high}$ and  $H_5^{low}$ terms defined  above, we remark that  \eqref{h5sum} holds.

We now estimate $H_5^{low,2}$.
Since $p$ and $q_1$ are both bounded by $M$, from Lemma \ref{nuEST} 
we still have \eqref{MdecayEST}.
For any $j\ge 0$,
 Lemma \ref{boundK2} implies the following bound
 $$
\int_{|q_1|\le M}  dq_1 ~
\int_{\mathbb{R}^3} dq_2~
\left|  w_{j}(q_2)  k^\chi(p,q_1)  k^\chi(q_1,q_2) \right|^2
\le
C_M.
 $$
 Indeed if  $|q_2| \ge 2 |q_1|$ then as in \eqref{comparablePQ} we can prove this bound.  Alternatively if 
 $|q_2| \le 2 |q_1|$ then $w_{j}(q_2) \le C w_{j}(q_1) \le CM^{jb/2}$  and the bound above also holds.
 
We use the above and Cauchy-Schwartz to estimate the momentum integrals:
\begin{multline*}
\left|
\int_{|q_1|\le M}  dq_1 ~ k^\chi(p,q_1)~ 
\int_{\mathbb{R}^3} dq_2~ k^\chi(q_1,q_2)~ 
\left\{  U(s_2) f_0 \right\}(y_2,q_2)
\right|
\\
\lesssim
\left(
\int_{|q_1|\le M}  dq_1 ~
\int_{\mathbb{R}^3} dq_2~
\left|  w_{j}(q_2)  k^\chi(p,q_1)  k^\chi(q_1,q_2) \right|^2
\right)^{1/2}
\\
\times
\left(
\int_{|q_1|\le M}  dq_1 ~
\int_{\mathbb{R}^3} dq_2~
\left| w_{-j}(q_2) \left\{  U(s_2) f_0 \right\}(y_2,q_2)\right|^2
\right)^{1/2}
\\
\le
C_M
\left(
\int_{|q_1|\le M}  dq_1 ~
\int_{\mathbb{R}^3} dq_2~
\left| w_{-j}(q_2) \left\{  U(s_2) f_0 \right\}(y_2,q_2) \right|^2
\right)^{1/2}.
\end{multline*}
That step was significant to yield rapid polynomial momentum decay.
 We   change  variables $q_1\to y_2$ on the $dq_1$ integration with $y_2$ given by
 \eqref{zDEF}.
Then
\begin{equation}
\left(\frac{d y_2 }{d q_1} \right)_{mn}
=
-(s_1-s_2)\left(\frac{\delta_{mn}q_{10}^2 - q_{1m} q_{1n}}{q_{10}^3} \right).
\label{jacobian11}
\end{equation}
This is a $3\times 3$ matrix with two eigenvalues equal to $-\frac{(s_1-s_2)}{q_{10}}$,
and a third eigenvalue given by
$-(s_1-s_2)\frac{q_{10}^2-|q_1|^2}{q_{10}^3} = \frac{-(s_1-s_2)}{q_{10}^3}$.
Thus the Jacobian is
$$
\left| \frac{d y_2 }{d q_1} \right| = \frac{\left| (s_1-s_2)\right|^3}{q_{10}^5} \ge C\frac{\kappa^3}{M^5}.
$$
This lower bound holds on the set $|q_1|\le M$, $s_2 \in [0, s_1-\kappa]$ and $s_1\in [\kappa, t]$
so that $s_1 - s_2 \ge \kappa$.
After application of this change of variables we have
\begin{multline*}
\left(
\int_{|q_1|\le M}  dq_1 ~
\int_{\mathbb{R}^3} dq_2~
\left| w_{-j}(q_2) \left\{  U(s_2) f_0 \right\}(y_2,q_2) \right|^2
\right)^{1/2}
\\
\lesssim
\left(\frac{M^5}{\kappa^3} \right)^{1/2}
\left(
\int_{|y_2-x|\le C (t-s_2)}   dy_2 ~
\int_{\mathbb{R}^3} dq_2~
\left| w_{-j}(q_2) \left\{  U(s_2) f_0 \right\}(y_2,q_2) \right|^2
\right)^{1/2}
\\
\lesssim
\left(\frac{M^5}{\kappa^3} \right)^{1/2}
\{1+(t-s_2)^{3/2}\}
\left(
\int_{{\mathbb{T}^3}}  dy_2 ~
\int_{\mathbb{R}^3} dq_2~
\left| w_{-j}(q_2) \left\{  U(s_2) f_0 \right\}(y_2,q_2) \right|^2
\right)^{1/2}
\\
= C(M,\kappa)\{1+(t-s_2)^{3/2}\} \| \left\{  U(s_2) f_0 \right\}\|_{2,-j} .
\end{multline*}
Putting together all of these estimates, in particular using \eqref{MdecayEST}, we have shown
\begin{multline*}
\left|
w_\ell(p)
H_5^{low,2}
\right|
\le
C_{\kappa, M} 
\int_{\kappa}^t  ~ds_1 ~
\int_0^{s_1-\kappa} ~ ds_2 ~
 e^{-C(t- s_1)/M^{b/2} }
e^{-C(s_1- s_2)/M^{b/2} }
\\
\times
\{1+(t-s_2)^{3/2}\} \| \left\{  U(s_2) f_0 \right\}\|_{2,-j}
\\
\le
C_{\kappa, M} 
\int_{0}^t  ~ds_1~   e^{-\frac{C}{2}(t- s_1)/M^{b/2} }
\int_0^{t} ~ ds_2 ~
e^{-\frac{C}{2}(t- s_1)/M^{b/2} }
e^{-\frac{C}{2}(s_1- s_2)/M^{b/2} }
\\
\times
\{1+(t-s_2)^{3/2}\} \| \left\{  U(s_2) f_0 \right\}\|_{2,-j}.
\end{multline*}
Notice that the first exponential controls the $s_1$ time integral, and the second and third exponential control the remaining time integral as follows
\begin{multline*}
\left|
w_\ell(p)
H_5^{low,2}
\right|
\le
C_{\kappa, M} 
\int_0^{t} ~ ds_2 ~
e^{-\frac{C}{2}(t- s_2)/M^{b/2} }
\{1+(t-s_2)^{3/2}\} \| \left\{  U(s_2) f_0 \right\}\|_{2,-j}
 \\
\le
C_{\kappa, M} 
\int_0^{t} ~ ds_2 ~
e^{-\frac{C}{4}(t- s_2)/M^{b/2} }
 \| \left\{  U(s_2) f_0 \right\}\|_{2,-j}.
\end{multline*}  
That was the last case.
Adding up the individual estimates for $H_5^{high}$ and  $H_5^{low}$  in \eqref{h5sum} completes our proof after first choosing $M$ large enough and then second choosing $\kappa$ sufficiently small.	
\qed \\

\section{Nonlinear $L^\infty$ Bounds and slow Decay} 
\label{s:n0bd}

Suppose $f=f(t,x,p)$ solves \eqref{rBoltz}
with initial condition $f(0,x,p) =  f_0(x,p)$.  
We may express mild solutions to this problem \eqref{rBoltz} in the form
\begin{equation}
f(t,x,p) = \{U(t)f_0\}(x,p)+N[ f, f](t,x,p),
\label{nonLinProbSol}
\end{equation}
where we have used the notation
$$
N[ f_1, f_2](t,x,p)
\eqdef
\int_0^t ~ ds ~ \{U(t-s) \Gamma [ f_1(s), f_2(s)]\}(x,p).
$$
Here as usual $U(t)$ is the semi-group \eqref{Udef} which represents solutions to the linear problem 
\eqref{rBlin}.  The main result of this section is to prove the following

\begin{theorem}
\label{existenceandminordecay}
Choose $\ell >3/b$, $k\in (1/2, 1]$.
Consider the following initial data $f_0 = f_0(x,p) \in L^\infty_{\ell+k}(\mathbb{T}^3 \times \mathbb{R}^3)$ which satisfies \eqref{conservation} initially.
  There is an $\eta>0$ such that if $\|f_0\|_{\infty, \ell+k} \le \eta$, then   there exists a unique global in time mild solution \eqref{nonLinProbSol}, $f = f(t,x,p)$, to equation \eqref{rBoltz} which satisfies
$$
\|f\|_{\infty, \ell}(t) \le C_{\ell, k}(1+t)^{-k} \|f_0\|_{\infty, \ell+k}.
$$
These solutions are continuous if it is so initially.  We further have positivity, i.e. $F= \mu + \sqrt{\mu} f \ge 0$, if $F_0= \mu + \sqrt{\mu} f_0 \ge 0$.
\end{theorem}

In Theorem \ref{rapidDECAYb}, which is proven in the Section \ref{s:rapidD},  we will show that these solutions exhibit rapid polynomial decay to any order.
 Notice that our main Theorem \ref{mainGLOBAL} will follow directly from  Theorem \ref{existenceandminordecay}
and Theorem \ref{rapidDECAYb}.
To prove this current Theorem \ref{existenceandminordecay} we will use the following non-linear estimate.

\begin{lemma}
\label{nonlin0}  
Considering the non-linear operator defined in \eqref{gamma0} with \eqref{hypSOFT},
we have the following pointwise estimates  
$$
\left| w_\ell(p) \Gamma (h_1, h_2)(p) \right|
\lesssim  \nu(p)
\|  h_1 \|_{\infty, \ell}  \|  h_2 \|_{\infty, \ell}. 
$$
These hold for any $\ell \ge 0$.
Furthermore,
$
\left\|  \Gamma (h_1, h_2) \right\|_{\infty, \ell+1}
\lesssim
\|  h_1 \|_{\infty, \ell}  \|  h_2 \|_{\infty, \ell}. 
$
\end{lemma}

The lemma above combined with Proposition \ref{BasicDecay} will be important tools in our proof of Theorem \ref{existenceandminordecay}.  We now give a simple proof.\\

\noindent {\it Proof of Lemma \ref{nonlin0}.}   We recall \eqref{nuDEF}, \eqref{gamma0}, and \eqref{weight}. 
For $\ell \ge 0$, it follows from \eqref{collisionalCONSERVATION} that
$$
w_\ell (p) 
\lesssim p_0^{\ell b/2}
\lesssim (p_0')^{\ell b/2}  (q_0')^{\ell b/2}
\lesssim w_\ell (p') w_\ell (q').
$$
A proof of this estimate above was given in \cite[Lemma 2.2]{MR1211782}.
Thus 
\begin{multline*}
w_\ell(p) \left| \Gamma (h_1,h_2) \right|
\lesssim
\int_{\mathbb{R}^3\times \mathbb{S}^{2}} ~ d\omega dq 
~ v_{\o} ~ \sigma(g,\theta)~
 \sqrt{J(q)} ~ 
 w_\ell (p') w_\ell (q') \left| h_1(p^{\prime })h_2(q^{\prime}) \right|
 \\
 +
\int_{\mathbb{R}^3\times \mathbb{S}^{2}} ~ d\omega dq 
~ v_{\o} ~ \sigma(g,\theta)~
 \sqrt{J(q)} ~ 
w_\ell(p) \left| h_1(p) h_2(q) \right|
 \\
\lesssim
  \|  h_1 \|_{\infty, \ell}  \|  h_2 \|_{\infty, \ell} 
  \int_{\mathbb{R}^3\times \mathbb{S}^{2}} ~ d\omega dq 
~ v_{\o} ~ \sigma(g,\theta)~
 J^{1/2}(q)
  \\
\lesssim
 \nu(p) \|  h_1 \|_{\infty, \ell}  \|  h_2 \|_{\infty, \ell}.
\end{multline*}
The last inequality above follows directly from Lemma \ref{nuEST} since both the integral and $\nu(p)$ have the same asymptotic behavior at infinity.
That yields the first estimate.  
For the second estimate we notice 
from the first estimate that
$
w_{\ell+1}(p) \left| \Gamma (h_1,h_2) \right|
\lesssim
w_{1}(p)  \nu(p) \|  h_1 \|_{\infty, \ell}  \|  h_2 \|_{\infty, \ell}.
$
But $w_{1}(p)  \nu(p) \lesssim 1$ from 
Lemma \ref{nuEST} and \eqref{weight}.  This
 completes the proof.
\qed \\

We now proceed to the\\

\noindent {\it Proof of Theorem \ref{existenceandminordecay}}.  We will prove 
Theorem \ref{existenceandminordecay}
in three steps.   The first step gives existence, uniqueness and slow decay via the contraction mapping argument.  The second step will establish continuity, and the last step shows positivity.

{\bf Step 1. Existence and Uniqueness}.  
When proving existence of mild solutions to \eqref{nonLinProbSol} it is natural to consider the mapping
$$
M[f]
\eqdef
\{U(t)f_0\}(x,p)+N[ f, f](t,x,p).
$$
With the norm \eqref{dNORM},
we will show that this is a contraction mapping on  the space 
$$
M_{k,\ell}^{R} \eqdef \{f \in L^\infty([0,\infty) \times \mathbb{T}^3 \times \mathbb{R}^3) :  ||| f |||_{k,\ell} \le R \},
\quad 
R>0.
$$
We first estimate the non-linear term $N[ f, f]$ defined in the equation display below \eqref{nonLinProbSol}.
We apply Theorem \ref{decay0l}, with $\ell >3/b$,  and $k\in(1/2,1]$  to obtain
\begin{multline*}
w_{\ell}(p)
\left|
N[ f_1, f_2](t,x,p)
\right|
\lesssim
\int_0^t ~ ds ~
w_{\ell}(p)
\left|
 \{U(t-s) \Gamma [ f_1(s), f_2(s)]\}(x,p)
\right|
\\
\lesssim
\int_0^t ~  \frac{ds}{(1+t-s)^k}
\left\|
  \Gamma [ f_1(s), f_2(s)]
\right\|_{\infty,\ell+k}.
\end{multline*}
  Next Lemma \ref{nonlin0} allows us to bound the above by
\begin{gather*}
\lesssim
\int_0^t ~  \frac{ds}{(1+t-s)^k}
\left\|  f_1(s) \right\|_{\infty,\ell+k -1}
\left\|  f_2(s) \right\|_{\infty,\ell+k -1}.
\end{gather*}
From Proposition \ref{BasicDecay} and the decay norm \eqref{dNORM} we see that the last line is
\begin{gather*}
\lesssim
||| f_1 |||_{k,\ell}
||| f_2 |||_{k,\ell}
\int_0^t ~  \frac{ds}{(1+t-s)^k (1+s)^{2k}}
\\
\lesssim
(1+t)^{-k}
||| f_1 |||_{k,\ell}
||| f_2 |||_{k,\ell}.
\end{gather*}
We have shown 
$$
|||  N[ f_1, f_2] |||_{k,\ell}
\lesssim
||| f_1 |||_{k,\ell}
||| f_2 |||_{k,\ell}.
$$
To handle the linear semigroup,  $U(t)$, we again use Theorem \ref{decay0l} to obtain
$$
||| M[ f] |||_{k,\ell}
\le
 C_{\ell,k}\left(\left\|  f_0 \right\|_{\infty,\ell+k}
 +
||| f |||_{k,\ell}^2
 \right).
$$
We conclude that $M[\cdot]$ maps $M_{k,\ell}^{R}$ into itself for $0<R$ chosen sufficiently small and e.g.
$
\left\|  f_0 \right\|_{\infty,\ell+k} \le \frac{R}{2 C_{\ell,k}}.
$
To obtain a contraction, we consider the difference
$$
M[ f_1] - M[ f_2] 
=
N[ f_1-f_2, f_1]
+
N[ f_2, f_1-f_2].
$$
Then as in the previous estimates we have
$$
|||  M[ f_1] - M[ f_2] |||_{k,\ell} 
\le
C_{\ell,k}^*
\left( 
|||  f_1 |||_{k,\ell}
 +
|||  f_2 |||_{k,\ell}
 \right)
|||  f_1 - f_2 |||_{k,\ell}.
$$
With these estimates, the existence and uniqueness of solutions to \eqref{rBoltz} follows from the contraction mapping principle on $M_{k,\ell}^{R}$ when $R>0$ is suitably small.

 {\bf Step 2. Continuity}.  We perform the estimates from Step 1 on the space
 $$
M_{k,\ell}^{R,0} \eqdef \{f \in C^0 ([0,\infty) \times \mathbb{T}^3 \times \mathbb{R}^3) :  ||| f |||_{k,\ell} \le R \},
\quad 
R>0.
$$
As in Step 1, 
we have a uniform in time contraction mapping on $M_{k,\ell}^{R,0}$ for suitable $R$.  
Furthermore $M[f]$ is continuous if $f\in M_{k,\ell}^{R,0}$  and $f_0$ is continuous.  
Since the convergence is uniform, the limit will be continuous  globally in time.
This argument is standard and we refer for instance to \cite{MR1379589,MR1211782,G2} for full details.

{\bf Step 3. Positivity}. 
We use the standard alternative approximating formula
$$
\left(\partial_t +\hat{p}\cdot \nabla_x \right) F^{n+1} + R(F^{n})F^{n+1}
=
\mathcal{Q}_+(F^{n},F^{n}).
$$
with the same initial conditions 
$
\left. F^{n+1}\right|_{t=0} =  F_0= J + \sqrt{J} f_0,
$
for $n\ge 1$ and for instance
$
 F^{1} \eqdef J + \sqrt{J} f_0.
$
Here we have used the standard decomposition of 
$
\mathcal{Q}
=
\mathcal{Q}_+ - \mathcal{Q}_-
$
into gain and loss terms with
$$
\mathcal{Q}_-(F^{n+1},F^{n}) = R(F^{n})F^{n+1},
$$
and
$
R(F^{n}) \eqdef \mathcal{Q}_-(1,F^{n}).
$
If we consider 
$
 F^{n+1}(t,x,p) = J + \sqrt{J} f^{n+1}(t,x,p),
$
then related to Step 1 we may show that 
$
f^{n+1}(t,x,p)
$
is convergent in $L^\infty_\ell$ on a local time interval $[0,T]$ where $T$ will generally depend upon the size of the initial data. 
In particular $f^{n+1}(t,x,p) = \frac{ F^{n+1} - J}{\sqrt{J}}$ satisfies the equation
$$
\left(\partial_t +\hat{p}\cdot \nabla_x  + \nu(p)\right)f^{n+1}
= K(f^{n})
+
\Gamma_+(f^{n},f^{n})
-
\Gamma_- (f^{n+1},f^{n}).
$$
We rewrite this equation using the solution formula to the system \eqref{rBlinWO} as
 $$
f^{n+1}
= 
G(t) f_0
+
\mathcal{L}(f^{n+1}, f^{n}).
$$
This solution formula $G(t)$ is defined just below \eqref{rBlinWO}.  Furthermore
\begin{multline*}
\mathcal{L}(f^{n+1}, f^{n})
\eqdef
 \int_0^t ~ ds ~ 
G(t-s) K(f^{n}) 
\\
+ \int_0^t ~ ds ~ 
G(t-s)  \Gamma_+(f^{n},f^{n})  
-
G(t-s)  \Gamma_- (f^{n+1},f^{n}).
\end{multline*}
For given $T>0$ and  $R > 0$ we consider the space $M_{k,\ell}^{R}([0,T])$ defined by
 $$
 \left\{f \in L^\infty ([0,T] \times \mathbb{T}^3 \times \mathbb{R}^3) :  \esssup_{0\le t \le T}(1+t)^k\| f(t) \|_{\infty,\ell} \le R \right\}.
$$
Now given 
$f^{n} \in  M_{k,\ell}^{R}([0,T])$
and
$
\left\|  f_0 \right\|_{\infty,\ell+k} \le \frac{R}{2 C_{k,\ell}}
$
with $R>0$ chosen sufficiently small, as in Step 1, we can prove the existence  of  
$f^{n+1} \in  M_{k,\ell}^{R}([0,T])$.

With the estimates established in this paper, it is now not hard to show that
$$
\sup_{0\le t \le T}\| f^{n+1} - f^{n} \|_{L^\infty_\ell }(t) \le C T  \sup_{0\le t \le T} \| f^{n} - f^{n-1} \|_{L^\infty_\ell }(t).
$$ 
Here
 $T>0$  is sufficiently small, and
 the constant $C>0$ can be chosen independent of any small $T$.  Therefore there exists a $T^*>0$ such that  $f^{n} \to f$ uniformly in $L^\infty_\ell$ on $[0,T^*]$. 
This will be sufficient to prove the positivity globally in time.  

Indeed if $F^{n} \ge 0$, then so is $\mathcal{Q}_+(F^{n},F^{n})\ge 0$.  With the  representation formula
\begin{multline*}
F^{n+1}(t,x,p) 
=
e^{-\int_0^t  ds  R(F^{n})(s, x- \hat{p}(t-s), p)} ~
F_0(x-\hat{p}t, p)
\\
+
\int_0^t  ds ~ e^{-\int_s^t  d\tau  R(F^{n})(\tau, x- \hat{p}(t-\tau), p)} ~ \mathcal{Q}_+(F^{n},F^{n})(s, x- \hat{p}(t-s), p).
\end{multline*}
Induction shows  $F^{n+1}(t,x,p)  \ge 0$ for all $n\ge 0$ if $F_0\ge 0$,
which implies in the limit $n\to\infty$ that $F(t,x,p)= J + \sqrt{J} f(t,x,p)\ge 0$.  Using our $L^\infty_\ell$ uniqueness, this is the same  $F$ as the one from Step 1 on the time interval $[0,T^*]$.  We  extend this positivity for all time intervals $[0,T^*] + T^* k$ for any $k\ge 1$ by repeating this  procedure and using the global uniform bound in $L^\infty_\ell$ from Step 1.
\qed \\

\section{Linear $L^\infty$ rapid Decay}
\label{s:linRdecay}

In this section we prove that the linear semi-group \eqref{Udef} exhibits rapid polynomial decay in $L^\infty_\ell$.  For any $k\ge 1$ we will discover a $k'=k'(k)\ge k$ such that
\begin{equation}
\| \{U(t)f_0\}\|_{\infty, \ell}
\le
C_{\ell,k} (1+t)^{-k} \| f_0\|_{\infty, \ell+k'}.
\label{Udecay}
\end{equation}
The main obstruction to proving such rapid decay in this low regularity $L^\infty_\ell$ framework was the term \eqref{r1def}
which came up during the course the proof of Lemma \ref{hESTIMATE3}.   In this section we perform a new high-order expansion of this remainder which allows one to prove rapid decay as follows.

\begin{proposition}
\label{r1rapidDecay}
Consider $R_1(f)(t)$ defined in \eqref{r1def}. 
Choose $\ell> 3/b$.
For any small $\eta>0$, and any $k\ge 1$, there exists a $k'=k'(k)\ge k$ such that  
\begin{gather*}
w_\ell(p)
\left|
R_1(f)(t)
\right|
\le
\eta (1+t)^{-k}||| f |||_{k,\ell}
+ C_{\ell,k',\eta} (1+t)^{-k}\| f_0\|_{\infty,\ell+k'}.
\end{gather*}
The power $k'$ can be explicitly computed from the proof.
\end{proposition}

The crucial difficulty with proving rapid decay for the soft potentials is caused by the high momentum values, for which the time decay is diluted by the momentum decay.  This causes the generation of weights on the initial data, typically one weight for each order of time decay.   In the proof below we are able to overcome this apparent obstruction by performing a new high order expansion for $R_1(f)$ which is explained in detail at the beginning of the proof.

We will first show that Proposition \ref{r1rapidDecay} implies \eqref{Udecay}.  We use the expansion \eqref{linearEXP}
and the semi-group notation
$
f(t) = \{U(t)f_0\}.
$
We  now see that 
Lemmas \ref{hESTIMATES}, \ref{hESTIMATE24}, and \ref{hESTIMATE3} together imply that for any $\eta >0$ and $k\ge 1/2$ we have
$$
\|f \|_{\infty,\ell}(t)
\le
C_{\ell,k} (1+t)^{-k}\| f_0\|_{\infty,\ell+ k}
+
\frac{\eta}{2} (1+t)^{-k} ||| f |||_{k,\ell}
+
w_\ell(p) \left|
R_1(f)(t)
\right|.
$$
Here we use
\eqref{dNORM}.
Then Proposition \ref{r1rapidDecay} further implies for some $k' \ge k$ that
$$
w_\ell(p)
\left|
R_1(f)(t)
\right|
\le
C_{\ell,k,\eta} (1+t)^{-k}\| f_0\|_{\infty,\ell+ k'}
+
\frac{\eta}{2} (1+t)^{-k} ||| f |||_{k,\ell}.
$$
Equivalently
$$
 ||| f |||_{k,\ell}
\le
C_{\ell,k'} \| f_0\|_{\infty,\ell+ k'}.
$$
This is the desired decay rate for the $L_\ell^\infty$ norm of mild solutions to the linear equation \eqref{rBlin},
which proves \eqref{Udecay} subject to  Proposition \ref{r1rapidDecay}.  In the rest of this section we prove this crucial new proposition.  \\

\noindent {\it Proof of Proposition \ref{r1rapidDecay}.} We will prove this lemma with a new high order expansion of \eqref{r1def}  by iterating the the semi-group \eqref{duhamelDOUBLE}.  For ease of exposition we write
$$
k_{\approx}^\chi(p,q_1)
\eqdef
k^\chi(p,q_1){\bf 1}_{\frac{1}{2}|p| \le |q_1|\le 2|p|}.
$$
Recall that ${\bf 1}_{\frac{1}{2}|p| \le |q_1|\le 2|p|}$ is the function which is one when ${\frac{1}{2}|p| \le |q_1|\le 2|p|}$ and zero elsewhere.
We will use similar expressions for $k_{\approx}^\chi$ with different arguments.  
Then we may split \eqref{r1def} as
\begin{gather*}
R_1(f)
=
S_1(f)  + L_1(f) + R_2(f).
\end{gather*}
For any small $\kappa>0$ we choose $\kappa_1 =\kappa$ and $\kappa_2= \kappa_1/2$ so that
\begin{gather*}
S_1(f) \eqdef
\int_0^{\kappa_1}  ~ds_1 ~ e^{-\nu(p) (t-s_1) }
 \int_{\mathbb{R}^3}  dq_1 ~ k_{\approx}^\chi(p,q_1)
~ {\bf 1}_{high} ~
\int_0^{s_1} ~ ds_2 ~ e^{-\nu(q_1) (s_1- s_2) }
\\
\times
\int_{\mathbb{R}^3} dq_2 ~ k_{\approx}^\chi(q_1,q_2) ~ 
\left\{  U(s_2) f_0 \right\}(y_2,q_2),
\\
L_1(f) \eqdef
\int_{\kappa_1}^t  ~ds_1 ~ e^{-\nu(p) (t-s_1) }
 \int_{\mathbb{R}^3}  dq_1 ~ k_{\approx}^\chi(p,q_1)
~ {\bf 1}_{high} ~
\int_0^{s_1-\kappa_2} ~ ds_2 ~ e^{-\nu(q_1) (s_1- s_2) }
\\
\times
\int_{\mathbb{R}^3} dq_2 ~ k_{\approx}^\chi(q_1,q_2) ~ 
\left\{  U(s_2) f_0 \right\}(y_2,q_2).
\end{gather*}
Then we may define the remainder term as
\begin{gather*}
R_2(f) \eqdef
\int_{\kappa_1}^t  ~ds_1 ~ e^{-\nu(p) (t-s_1) }
 \int_{\mathbb{R}^3}  dq_1 ~ k_{\approx}^\chi(p,q_1)
~ {\bf 1}_{high} ~
\int_{s_1-\kappa_2}^{s_1} ~ ds_2 ~ e^{-\nu(q_1) (s_1- s_2) }
\\
\times
\int_{\mathbb{R}^3} dq_2 ~ k_{\approx}^\chi(q_1,q_2) ~ 
\left\{  U(s_2) f_0 \right\}(y_2,q_2).
\end{gather*}
Our notation above is from the proof of Lemma \ref{hESTIMATE3}, in particular \eqref{1high}.
We will show that the first term $S_1(f)$ exhibits rapid decay in $L^\infty_\ell$.  The last term $L_1(f)$ is bounded in $L^2$ which further has rapid decay as in Theorem \ref{decay2}.

We will notice first of all that the term $R_2(f)$ naively exhibits second order polynomial decay.  However if we continue the expansion then we can obtain higher and higher order decay rates as follows.  We may expand
\begin{gather*}
R_2(f)
=
G_2(f)  + D_2(f) + N_2(f)  + L_2(f) + R_3(f).
\end{gather*}
Now each of the terms $G_2(f)$,   $D_2(f)$,  $N_2(f)$,   and $L_2(f)$-- to be defined below--will exhibit (for different reasons) high order polynomial decay right away again at a cost of momentum weights on the initial data.  The term $R_3(f)$ will clearly exhibit third order polynomial decay, however we may continue this expansion at each level so that at level $k$ we can again expand  
\begin{gather*}
R_k(f)
=
G_k(f)  + D_k(f) + N_k(f)  + L_k(f) + R_{k+1}(f).
\end{gather*}
As in the initial case each of the terms $G_k(f)$,   $D_k(f)$,  $N_k(f)$,   and $L_k(f)$--which are defined recursively--will exhibit high order polynomial decay.  The last term $R_{k+1}(f)$ will have $k+1$ order polynomial decay.  This expansion is well defined and can be continued to any order, which yields rapid decay.

We  define the $2nd$ order terms by plugging the iteration \eqref{duhamelDOUBLE}  into $R_2(f)$, using the expansion of
$
K
=
K^{1-\chi}
+
K^{\chi}
$ with \eqref{cut}, 
 and splitting the remaining time and momentum integrals in the following useful  way.
For $\kappa_3 = \kappa_2/2 = \kappa_1/2^2$,
we define
\begin{gather*}
G_2(f)
\eqdef
\int_{\kappa_1}^t  ~ds_1 ~ e^{-\nu(p) (t-s_1) }
 \int_{\mathbb{R}^3}  dq_1 ~ k_{\approx}^\chi(p,q_1)
~ {\bf 1}_{high} ~
\int_{s_1-\kappa_2}^{s_1} ~ ds_2 ~ e^{-\nu(q_1) (s_1- s_2) }
\\
\times
\int_{\mathbb{R}^3} dq_2 ~ k_{\approx}^\chi(q_1,q_2) ~ 
\left\{  G(s_2) f_0 \right\}(y_2,q_2),
\\
D_2(f) \eqdef 
\int_{\kappa_1}^t  ~ds_1 ~ e^{-\nu(p) (t-s_1) }
 \int_{\mathbb{R}^3}  dq_1 ~ k_{\approx}^\chi(p,q_1)
~ {\bf 1}_{high} ~
\int_{s_1-\kappa_2}^{s_1} ~ ds_2 ~ e^{-\nu(q_1) (s_1- s_2) }
\\
\times
\int_{\mathbb{R}^3} dq_2 ~ k_{\approx}^\chi(q_1,q_2) ~ 
\int_{0}^{s_2} ~ ds_3 ~
e^{-\nu(q_2) (s_2- s_3) }
K^{1-\chi}\left\{  U(s_3) f_0 \right\}(y_3,q_2).
\end{gather*}
Furthermore
\begin{gather*}
N_2(f) 
\eqdef
\int_{\kappa_1}^t  ds_1 ~ e^{-\nu(p) (t-s_1) }
 \int_{\mathbb{R}^3}  dq_1 ~ k_{\approx}^\chi(p,q_1)
~ {\bf 1}_{high} 
\int_{s_1-\kappa_2}^{s_1}  ds_2 ~ e^{-\nu(q_1) (s_1- s_2) }
\\
\times
\int_{\mathbb{R}^3} dq_2 ~ k_{\approx}^\chi(q_1,q_2) ~ 
\int_{0}^{s_2}  ds_3 ~
e^{-\nu(q_2) (s_2- s_3) }
\int_{\mathbb{R}^3} dq_3~ k_{\neq}^\chi(q_2,q_3)
\left\{  U(s_3) f_0 \right\}(y_3,q_3),
\\
L_2(f)
\eqdef
\int_{\kappa_1}^t  ds_1 ~ e^{-\nu(p) (t-s_1) }
 \int_{\mathbb{R}^3}  dq_1 ~ k_{\approx}^\chi(p,q_1)
~ {\bf 1}_{high} ~
\int_{s_1-\kappa_2}^{s_1} ~ds_2 ~ e^{-\nu(q_1) (s_1- s_2) }
\\
\times
\int_{\mathbb{R}^3} dq_2 ~ k_{\approx}^\chi(q_1,q_2) ~ 
\int_{0}^{s_2-\kappa_3} ~ ds_3 ~
e^{-\nu(q_2) (s_2- s_3) }
\\
\times
\int_{\mathbb{R}^3} dq_3~ k_{\approx}^\chi(q_2,q_3)
\left\{  U(s_3) f_0 \right\}(y_3,q_3).
\end{gather*}
And the remainder is given by
\begin{gather*}
 R_3(f)
\eqdef
\int_{\kappa_1}^t  ds_1 ~ e^{-\nu(p) (t-s_1) } 
 \int_{\mathbb{R}^3}  dq_1 ~ k_{\approx}^\chi(p,q_1)
~ {\bf 1}_{high} ~
\int_{s_1-\kappa_2}^{s_1}  ds_2 ~ e^{-\nu(q_1) (s_1- s_2) }
\\
\times
\int_{\mathbb{R}^3} dq_2 ~ k_{\approx}^\chi(q_1,q_2) ~ 
\int_{s_2-\kappa_3}^{s_2} ~ds_3 ~
e^{-\nu(q_2) (s_2- s_3) }
\\
\times
\int_{\mathbb{R}^3} dq_3~ k_{\approx}^\chi(q_2,q_3)
\left\{  U(s_3) f_0 \right\}(y_3,q_3),
\end{gather*}
where above $y_1$ and $y_2$ are defined in \eqref{zDEF}, and more generally 
\begin{multline}
\label{yiDEF}
y_{i+1}  \eqdef  
y_1-\hat{q}_1(s_1-s_2)-\cdots -\hat{q}_i(s_i-s_{i+1})
\\
=
x- \hat{p}(t-s_1)-\hat{q}_1(s_1-s_2)-\cdots -\hat{q}_i(s_i-s_{i+1}).
\end{multline}
So that in general for $i\ge 1$ we have
$$
y_{i+1}  = y_{i} -\hat{q}_i(s_i-s_{i+1}).
$$
Furthermore $k^\chi(p,q) = k_{\neq}^\chi(p,q) + k_{\approx}^\chi(p,q)$ with the notation
$$
k_{\neq}^\chi(p,q) \eqdef k^\chi(p,q)\left( {\bf 1}_{|p| \ge 2|q|}
+
{\bf 1}_{ |q|\ge 2|p|} \right).
$$
Now we will develop a collection of notations in order to put this expansion into a general framework and appropriately define the high order terms.  
We consider the sequence $\{\kappa\}$ where for $i\ge 1$ we define $\kappa_{i+1}=\kappa_{i}/2$
with a small $\kappa_1 = \kappa>0$ as above 
so that
$\kappa_{i}=\kappa/2^{i-1}$.
 For $i\ge 1$ we can now define 
\begin{gather*}
 A(f)(t,x,p,\{\kappa\})
\eqdef
\int_{\kappa_1}^t  ~ds_1 ~ e^{-\nu(p) (t-s_1) }
 \int_{\mathbb{R}^3}  dq_1 ~ k_{\approx}^\chi(p,q_1)
~ {\bf 1}_{high} ~
 f(s_1,y_1,q_1),
\\
B_{i+1}(f)(s_1,y_1, q_1,\{\kappa\})
\eqdef
\int_{\mathbb{R}^3} dq_2 ~ k_{\approx}^\chi(q_1,q_2)
\cdots
\int_{\mathbb{R}^3} dq_{i+1} ~ k_{\approx}^\chi(q_{i},q_{i+1}) ~ 
\\
\times \int_{s_1-\kappa_{2}}^{s_1} ~ ds_{2} ~ e^{-\nu(q_1) (s_1- s_2) }
\cdots
 \int_{s_i-\kappa_{i+1}}^{s_i} ~ ds_{i+1} ~ e^{-\nu(q_i) (s_i- s_{i+1}) }
 f(s_{i+1},y_{i+1},q_{i+1}).
\end{gather*}
To finish off our expansion we 
 further define
\begin{gather*}
\tilde{G}(f)(s,y,q)
\eqdef
  G(s) f_0 (y,q),
\\
D(f)(s_{i+1},y_{i+1},q_{i+1})  \eqdef 
\int_0^{s_{i+1}} ds_{i+2} ~ e^{-\nu(q_{i+1}) (s_{i+1}- s_{i+2}) } ~
\\
\times
K^{1-\chi}\left\{  U(s_{i+2}) f_0 \right\}(y_{i+2},q_{i+1}).
\end{gather*}
Above we recall that $G(s)$ is defined just below \eqref{rBlinWO}.  Additionally
\begin{gather*}
N(f)(s_{i+1},y_{i+1},q_{i+1}) 
\eqdef
\int_{0}^{s_{i+1}} ~ ds_{i+2} ~ 
e^{-\nu(q_{i+1}) (s_{i+1}- s_{i+2}) }
\\
\times
\int_{\mathbb{R}^3} dq_{i+2}~ k_{\neq}^\chi(q_{i+1},q_{i+2})
\left\{  U(s_{i+2}) f_0 \right\}(y_{i+2},q_{i+2}),
\\
\tilde{L}(f)(s_{i+1},y_{i+1},q_{i+1},\{\kappa\}) 
\eqdef
\int_{0}^{s_{i+1}-\kappa_{i+2}} ~ ds_{i+2} ~
e^{-\nu(q_{i+1}) (s_{i+1}- s_{i+2}) } ~
\\
\times
\int_{\mathbb{R}^3} dq_{i+2}~ 
k_{\approx}^\chi(q_{i+1},q_{i+2})~
\left\{  U(s_{i+2}) f_0 \right\}(y_{i+2},q_{i+2}).
\end{gather*}
Then we may use this notation to write
\begin{gather*}
G_2(f)
\eqdef
A(B_2(\tilde{G}(f)))(t,x,p,\{\kappa\}),
\quad
D_2(f)
\eqdef
A(B_2(D(f)))(t,x,p,\{\kappa\}),
\\
N_2(f)
\eqdef
A(B_2(N(f)))(t,x,p,\{\kappa\}),
\quad
L_2(f)
\eqdef
A(B_2(\tilde{L}(f)))(t,x,p,\{\kappa\}),
\\
R_3(f)
\eqdef
A(B_3(f))(t,x,p,\{\kappa\}).
\end{gather*}
We iteratively define the higher order terms of this expansion for $i\ge 2$ as
\begin{gather*}
G_i
(f)
\eqdef
A(B_i
(\tilde{G}(f)))(t,x,p,\{\kappa\}),
\quad
D_i
(f)
\eqdef
A(B_i
(D(f)))(t,x,p,\{\kappa\}),
\\
N_i
(f)
\eqdef
A(B_i(N(f)))(t,x,p,\{\kappa\}),
\quad
L_i(f)
\eqdef
A(B_i(\tilde{L}(f)))(t,x,p,\{\kappa\}),
\\
R_{i+1}(f)
\eqdef
A(B_{i+1}(f))(t,x,p,\{\kappa\}).
\end{gather*}
This expansion works to high order using \eqref{duhamelDOUBLE} which implies
$$
B_{i} = B_{i} \tilde{G}+B_{i} D+B_{i} N+B_{i} \tilde{L} + B_{i+1},
\quad
i\ge 2.
$$
This completes our general discussion of the expansion formula, and our strategy for obtaining the desired decay.  
In the following we prove the claimed time decay estimates for each term in a general framework.

We initially estimate the main term, $R_{k+1}$, with $k\ge 1$.  We {\it claim} that
\begin{multline}
\label{claimBk}
w_\ell(q_1) \left| B_{k+1}(f)(s_1,y_1, q_1,\{\kappa\})\right| 
\\
\lesssim
w_{-k(1+\zeta)}(q_1) ~
||| f |||_{k+1,\ell} ~ \int_{s_1-\kappa_{2}}^{s_1} ~ ds_{2} ~ e^{-c\nu(q_1)(s_1 - s_{2})} ~ 
\left( 1 + s_2 \right)^{-k-1}.
\end{multline}
This claim \eqref{claimBk} would imply with Lemma \ref{boundK2}  that 
\begin{multline*}
w_\ell(p) \left| R_{k+1}(f)(t,x,p,\{\kappa\})\right| 
=
w_\ell(p)\left| A(B_{k+1}(f))(t,x,p,\{\kappa\})\right| 
\\
\lesssim
\int_{\kappa_1}^t  ~ds_1 ~ e^{-\nu(p) (t-s_1) }
 \int_{\mathbb{R}^3}  dq_1 ~ \left| k_{\approx}^\chi(p,q_1) \right|
~ {\bf 1}_{high} ~
w_{-k(1+\zeta)}(q_1) ~
||| f |||_{k+1,\ell} ~ 
\\
\times
\int_{s_1-\kappa_{2}}^{s_1} ~ ds_{2} ~ e^{-c\nu(q_1)(s_1 - s_{2})} ~ 
\left( 1 + s_2 \right)^{-k-1}
\\
\lesssim
\frac{w_{-(k+1)(1+\zeta)}(p)}{M^\zeta} 
||| f |||_{k+1,\ell} ~ \int_{0}^{t} ~ ds_{1} ~ e^{-c\nu(p)(t - s_{1})} ~ \left( 1+ s_1 \right)^{-k-1}.
\end{multline*}
We have used that $s_1 \approx s_2$ for $\kappa$ sufficiently small, and also  $e^{-c\nu(q_1)(s_1 - s_{2})} \lesssim 1$ since $s_1 - s_{2} \ge 0$.  We have additionally used the fact that the momentum variables are comparable because of the support condition for $k^\chi_{\approx}$.  The large $M$ comes from the support of ${\bf 1}_{high}$ in \eqref{1high} and Lemma \ref{boundK2}.  Furthermore 
\begin{multline*}
\lesssim
\frac{1}{M^\zeta} 
~
||| f |||_{k+1,\ell} ~ \int_{0}^{t} ~ ds_{1} ~\left( 1+ t-s_1 \right)^{-k-1-\zeta} ~ \left( 1+ s_1 \right)^{-k-1}
\\
\lesssim
\frac{1}{M^\zeta} 
~ \left( 1 + t \right)^{-k-1} 
||| f |||_{k+1,\ell}.
\end{multline*}
We have additionally used  \eqref{ElemCalc}  and then Proposition \ref{BasicDecay}.
For $M\gg1$ chosen sufficiently large, this is the desired estimate for $R_{k+1}$.

We now prove the {\it claim} from \eqref{claimBk}.  Since all of the momentum variables are comparable in this operator we have the following iterated estimate
\begin{equation}
e^{-\nu(q_1) (s_1- s_2) }
\cdots
e^{-\nu(q_{k}) (s_{k}- s_{k+1}) }
\le
 e^{-C\nu(q_1) (t-s_{k+1}) }.
 \label{iteratedE}
\end{equation}
This uses in particular Lemma \ref{nuEST}.  We use \eqref{iteratedE} to obtain
\begin{multline*}
w_\ell(q_1) 
\left| B_{k+1}(f)(s_1,y_1, q_1,\{\kappa\}) \right|
\\
\lesssim
||| f |||_{k+1,\ell}
\int_{\mathbb{R}^3} dq_2 ~  \left| k_{\approx}^\chi(q_1,q_2) \right|
\cdots
\int_{\mathbb{R}^3} dq_{k+1} ~  \left| k_{\approx}^\chi(q_{k},q_{k+1})   \right|
\\
\times 
\int_{s_1-\kappa_{2}}^{s_1} ~ ds_{2} ~ e^{-\nu(q_1) (s_1- s_2) }
\cdots
 \int_{s_k-\kappa_{k+1}}^{s_k} ~ ds_{k+1} ~ e^{-\nu(q_k) (s_k- s_{k+1}) }
(1+s_{k+1})^{-k-1}
\\
\lesssim
||| f |||_{k+1,\ell}
\int_{\mathbb{R}^3} dq_2 ~ \left| k_{\approx}^\chi(q_1,q_2) \right|
\cdots
\int_{\mathbb{R}^3} dq_{k+1} ~  \left| k_{\approx}^\chi(q_{k},q_{k+1})  \right|
\\
\times 
\int_{s_1-\kappa_{2}}^{s_1} ~ ds_{2} ~ 
\cdots
 \int_{s_k-\kappa_{k+1}}^{s_k} ~ ds_{k+1} ~ e^{-C\nu(q_1) (s_1- s_{k+1}) }
(1+s_{k+1})^{-k-1}.
\end{multline*}
We {\it sub-claim} that  Lemma \ref{boundK2} can be used to control the momentum integrals as
\begin{equation}
\int_{\mathbb{R}^3} dq_2 ~  \left| k_{\approx}^\chi(q_1,q_2)  \right|
\cdots
\int_{\mathbb{R}^3} dq_{k+1} ~  \left| k_{\approx}^\chi(q_{k},q_{k+1})   \right|
\lesssim
 w_{-k(1+\zeta)}(q_1).
 \label{subclaim}
\end{equation}
The trick used in this estimate is Lemma \ref{boundK2} combined with
$
q_1 - q_2 \to q_2 
$
to obtain
$$
\int_{\mathbb{R}^3} dq_2 ~ \left|   k_{\approx}^\chi(q_1,q_2) \right|
\lesssim
 w_{-(1+\zeta)}(q_1)
 \int_{\mathbb{R}^3} dq_2 ~ 
e^{-c |q_2|}.
$$
We can do that $k$ times starting with the $dq_{k+1}$ integral and iterating backwards (using the essential point that all the momentum variables are comparable) to obtain the {\it sub-claim}.
Now by the definition of the sequence $\{\kappa\}$ we can say that 
$s_{k+1} \le s_2$ for $k\ge 1$ and more generally (on the integration region of $B_{k+1}$)
\begin{multline*}
s_{k+1} \ge s_{k} - \kappa_{k+1} \ge  s_{2} - \kappa_{3} - \cdots - \kappa_{k+1}
\\
=
s_{2} - \kappa\left(\frac{1}{4} + \cdots + \frac{1}{2^k}  \right)
\ge 
s_{2} - \frac{\kappa}{2}
\ge 
s_{2} - \frac{1}{4}.
\end{multline*}
Thus for $\kappa \le 1/2$ we use these estimates above to obtain
\begin{multline*}
\int_{s_1-\kappa_{2}}^{s_1} ~ ds_{2} ~ 
\cdots
 \int_{s_k-\kappa_{k+1}}^{s_k} ~ ds_{k+1} ~ e^{-c\nu(q_1) (s_1- s_{k+1}) }
(1+s_{k+1})^{-k-1}
\\
\lesssim
\kappa_{k+1} \cdots \kappa_{3} \int_{s_1-\kappa_{2}}^{s_1} ~ ds_{2} ~  e^{-c\nu(q_1) (s_1- s_{2}) }
\left(\frac{1}{2}+s_{2}\right)^{-k-1}.
\end{multline*}
Collecting the estimates above proves the claim \eqref{claimBk}.

To estimate the first term above, $S_1$, we 
obtain
\begin{multline*}
w_\ell(p)
\left|
S_1(f)(t)
\right|
\lesssim
w_\ell(p)
 \int_{\mathbb{R}^3}  dq_1 ~ \left| k_{\approx}^\chi(p,q_1) \right|
\int_{\mathbb{R}^3} dq_2~    \left|  k_{\approx}^\chi(q_1,q_2) \right|
~ {\bf 1}_{high}
\\
\times
\int_0^{\kappa_1}  ~ds_1 ~ e^{-\nu(p) (t-s_1) }
\int_0^{\kappa_1} ~ ds_2 ~
e^{-\nu(q_1) (s_1- s_2) }
\left|\left\{  U(s_2) f_0 \right\}(y_2,q_2)\right|.
\end{multline*}
Now we use \eqref{ElemCalc}, and Lemma \ref{nuEST}, to observe that
\begin{equation}
 e^{-\nu(p) t} \lesssim w_k(p) (1+t)^{-k}.
 \label{decayEstSoft22}
\end{equation}
Furthermore 
$
\int_0^{\kappa_1}  ~ds_1 ~ e^{\nu(p) s_1 }
\int_0^{\kappa_1} ~ ds_2 ~
e^{-\nu(q_1) (s_1- s_2) }
\lesssim \kappa_1^2.
$
We thus have
\begin{multline*}
w_\ell(p)
\left|
S_1(f)(t)
\right|
\lesssim
 \kappa_1^2(1+t)^{-k}
||| f |||_{0,\ell+k}
 \int_{\mathbb{R}^3}  dq_1 ~ \left| k_{\approx}^\chi(p,q_1) \right|
 \\
 \times
\int_{\mathbb{R}^3} dq_2~ \left| k_{\approx}^\chi(q_1,q_2) \right|
~
\frac{w_k(p){\bf 1}_{high}}{w_k(q_2)}
\\
\lesssim
\frac{\kappa_1^2 (1+t)^{-k}}{M^{2\zeta}} 
\| f_0\|_{\infty, \ell+k}.
\end{multline*}
We have  used the uniform bound from Theorem \ref{decay0l} with no decay and the bound for $\left| k^\chi  \right|$ from Lemma \ref{boundK2} and \eqref{1high}.  This is the desired estimate for $S_1$.  

We continue with an estimate for $L_i(f)$ with $i\ge 1$.  For all of the terms below we switch from the notation of $k$ to the notation of $i$ to indicate that the decay of each of these terms will not depend upon the index of the term, which is contrary to the decay of the $R_k$ terms above.
We estimate from above
\begin{multline*}
w_\ell(p)
\left|
L_i(f)
\right|
\lesssim
w_\ell(p)
 \int_{\mathbb{R}^3} dq_1 ~ k_{\approx}^\chi(p,q_1)
\cdots
\int_{\mathbb{R}^3} dq_{i+1}~ k_{\approx}^\chi(q_i,q_{i+1})
\\
\times
\int_{\kappa_1}^t  ds_1  e^{-\nu(p) (t-s_1) }  
\left(
\int_{s_1-\kappa_2}^{s_1}  ds_2 
e^{-\nu(q_1) (s_1- s_2) }
\cdots
\int_{s_{i-1}-\kappa_{i}}^{s_{i-1}}  ds_{i} 
e^{-\nu(q_{i-1}) (s_{i-1}- s_{i}) } \right)
\\
\times
{\bf 1}_{high}
\int_{0}^{s_{i}-\kappa_{i+1}} ~ ds_{i+1} ~
e^{-\nu(q_{i}) (s_{i}- s_{i+1}) } 
\left| \left\{  U(s_{i+1}) f_0 \right\}(y_{i+1},q_{i+1}) \right|.
\end{multline*}
The term in parenthesis above would be simply unity in the case of $L_1$.
Since all the momentum variables are comparable, we control the time integrals as
\begin{multline*}
\int_{\kappa_1}^t  ds_1  e^{-\nu(p) (t-s_1) }  
\int_{s_1-\kappa_2}^{s_1} ~ ds_2 ~
e^{-\nu(q_1) (s_1- s_2) }
\cdots
\int_{s_{i-1}-\kappa_{i}}^{s_{i-1}} ~ ds_{i} ~
e^{-\nu(q_{i-1}) (s_{i-1}- s_{i}) } 
\\
\times
\int_{0}^{s_{i}-\kappa_{i+1}} ~ ds_{i+1} ~
e^{-\nu(q_{i}) (s_{i}- s_{i+1}) }
\\
\lesssim 
\left(
{\kappa_2}
\cdots
{\kappa_i} \right)
\int_{\kappa_1}^t  ds_1  
\int_{0}^{t-\kappa_{i+1}} ~ ds_{i+1} ~
e^{-C\nu(p) (t-s_{i+1}) }.
\end{multline*}
We have used  $\nu(q_j) \ge C q_{j0}^{-b/2}$ from Lemma \ref{nuEST}, 
$\nu(p) \approx p_{0}^{-b/2}$, 
and $q_{j0}^{-b/2} \ge C p_{0}^{-b/2}$; these estimates hold for any $j\in\{1, \ldots, i\}$.
We have then used an
estimate analogous to \eqref{iteratedE}.
Furthermore 
\begin{multline*}
\left(
{\kappa_2}
\cdots
{\kappa_i} \right)
\int_{\kappa_1}^t  ds_1  
\int_{0}^{t-\kappa_{i+1}} ~ ds_{i+1} ~
e^{-C\nu(p) (t-s_{i+1}) }
\\
\lesssim
 \kappa^{i-1}
\left(
1+t \right)
\int_{0}^{t-\kappa_{i+1}} ~ ds_{i+1} ~
e^{-C\nu(p) (t-s_{i+1}) }
\\
\lesssim
 \kappa^{i-1}
\left(
1+t \right)
w_j(p)
\int_{0}^{t-\kappa_{i+1}} ~ ds_{i+1} ~
\left(
1+t-s_{i+1} \right)^{-j}.
\end{multline*}
These estimates follow from the definition of the sequence $\{\kappa\}$ as well as Lemma \ref{nuEST} together with \eqref{ElemCalc} 
in the form \eqref{decayEstSoft22}
for any $j \ge 0$.

Next we use Cauchy-Schwartz to estimate the following two integrals 
\begin{multline}
w_{\ell+j}(q_{i+1})
\int_{\mathbb{R}^3} dq_{i}  \left| k_{\approx}^\chi(q_{i-1},q_{i})\right|
\\
\times
\int_{\mathbb{R}^3} dq_{i+1} \left|  k_{\approx}^\chi(q_i,q_{i+1}) \right|
 \left| \left\{  U(s_{i+1}) f_0 \right\}(y_{i+1},q_{i+1}) \right|
\\
\lesssim
\left(
\int_{\mathbb{R}^3} dq_{i} 
\int_{\mathbb{R}^3} dq_{i+1}
\left|    w_{2}(q_{i+1})  k_{\approx}^\chi(q_{i-1},q_{i})k_{\approx}^\chi(q_i,q_{i+1}) \right|^2
\right)^{1/2}
\label{csESTell}
\\
\times
\left(
\int_{\mathbb{R}^3} dq_{i} 
\int_{Z_i} dq_{i+1}
\left| w_{\ell+j-2}(q_{i+1}) \left\{  U(s_{i+1}) f_0 \right\}(y_{i+1},q_{i+1})\right|^2
\right)^{1/2}.
\end{multline}
Above $Z_i \eqdef \{q_{i+1}: \frac{1}{2}|q_{i}| \le |q_{i+1}|\le 2|q_{i}|\}$.
Also in the case $i=1$ we consider $q_{i-1}=p$.
For now we focus on the second set of integrals involving the semi-group.
We apply the  change of variables $q_i\to y_{i+1}$ on the $dq_i$ integration with $y_{i+1}$ given by
 \eqref{yiDEF}.
Notice that similar to \eqref{jacobian11} the $3\times 3$ matrix Jacobian is
$$
\left(\frac{d y_{i+1} }{d q_{i}} \right)_{mn}
=
-(s_i-s_{i+1})\left(\frac{\delta_{mn}q_{i0}^2 - q_{im} q_{in}}{q_{i0}^3} \right).
$$
We recall $q_i = (q_{i1},q_{i2},q_{i3})$ with $q_{i0} = \sqrt{1+|q_i|^2}$.
This Jacobian matrix has two eigenvalues equal to $-\frac{(s_{i}-s_{i+1})}{q_{i0}}$,
and a third eigenvalue which is given by
$-(s_{i}-s_{i+1})\frac{q_{i0}^2-|q_i|^2}{q_{i0}^3} = -(s_{i}-s_{i+1})\frac{1}{q_{i0}^3}$.
Therefore the Jacobian determinant is
$$
\left| \frac{d y_{i+1} }{d q_{i}} \right| = \frac{\left| (s_{i}-s_{i+1})\right|^3}{q_{i0}^5}
 \ge 
 \frac{\kappa^3_{i+1}}{2^{5}q_{(i+1)0}^5}
 =
  \frac{\kappa^3}{2^{5+3i}q_{(i+1)0}^5}.
$$
This lower bound holds on the region $q_{i0}\le 2q_{(i+1)0}$, $s_{i+1} \in [0, s_{i}-\kappa_{i+1}]$. 
Furthermore  we have used that $s_i - s_{i+1} \ge \kappa_{i+1}$; this temporal estimate holds on the integration region of $L_i$.  These estimates explain the lower bound for the Jacobian.
Notice while the old variable $q_{i}$ occupies the whole space, the new variable $y_{i+1}$ satisfies the estimate
\begin{multline*}
\left| y_{i+1} - x\right| 
\le 
\left| \hat{p}(t-s_1)+\hat{q}_1(s_1-s_2)+\cdots +\hat{q}_i(s_i-s_{i+1})\right| 
\\
\le 
C \left( (t-s_1) + (s_1-s_2) +\cdots + (s_i-s_{i+1})\right)
\le
C \left( t-s_{i+1} \right). 
\end{multline*}
We remark that this procedure  would not hold in the non-relativistic situation, since in that case we do not have bounded velocities.
In particular, because of  relativity, the mapping 
$q_{i} \to y_{i+1}$ sends $\mathbb{R}^3$ into a bounded domain (for any finite $t$).
After application of this change of variables, denoting $y_{i+1} = y$, we have
\begin{multline*}
\int_{\mathbb{R}^3} dq_{i} 
\int_{Z_i} dq_{i+1}
\left| w_{\ell+j-2}(q_{i+1}) \left\{  U(s_{i+1}) f_0 \right\}(y_{i+1},q_{i+1})\right|^2
\\
\lesssim
\int_{|y-x|\le C(t-s_{i+1})}   dy
\left| \frac{d q_{i}}{d y } \right| 
\int_{Z_i} dq_{i+1}
\left| w_{\ell+j-2}(q_{i+1}) \left\{  U(s_{i+1}) f_0 \right\}(y,q_{i+1}) \right|^2
\\
\lesssim
\frac{\left( 1+t-s_{i+1} \right)^3}{\kappa^{3/2}}
\int_{{\mathbb{T}^3}}  dy ~
\int_{\mathbb{R}^3} dq_{i+1}~
w_{2\ell+2j-4+10/b}(q_{i+1})
\left| \left\{  U(s_{i+1}) f_0 \right\}(y,q_{i+1}) \right|^2
\\
= C(\kappa)\left( 1+t-s_{i+1}\right)^{3}  \| \left\{  U(s_{i+1}) f_0 \right\}\|_{2,\ell+j-2+5/b}^2 .
\end{multline*}
We have used that 
$
q_{(i+1)0}^5 = w_{10/b}(q_{i+1}).
$
This estimate above is the main one for the $L_i(f)$ terms which allows us to deduce high order decay.

Since we have used \eqref{csESTell}, for the momentum integrals in $L_i$, we are left to control
 the iteration of kernels.  We {\it claim} the following estimate
\begin{multline*}
 \left( \int_{\mathbb{R}^3} dq_1 ~ \left| k_{\approx}^\chi(p,q_1)  \right| 
\cdots
\int_{\mathbb{R}^3} dq_{i-1}~ \left| k_{\approx}^\chi(q_{i-2}, q_{i-1}) \right|
\right)
\\
\times
\left(
\int_{\mathbb{R}^3} dq_{i} 
\int_{\mathbb{R}^3} dq_{i+1}
\left|    w_{2}(q_{i+1})  k_{\approx}^\chi(q_{i-1},q_{i})k_{\approx}^\chi(q_i,q_{i+1}) \right|^2
{\bf 1}_{high}
\right)^{1/2}
\lesssim
\frac{1}{M^\zeta}.
\end{multline*}
Note that if $i=1$ then the first term in parenthesis above is simply unity.   
Firstly, from Lemma \ref{boundK2} we have
$$
w_{2}(q_{i+1})  k_{\approx}^\chi(q_{i-1},q_{i})k_{\approx}^\chi(q_i,q_{i+1})
{\bf 1}_{high}
\lesssim
\frac{1}{M^\zeta} 
e^{-c|q_{i-1}- q_{i}|}
e^{-c|q_i - q_{i+1}|}.
$$
This also uses \eqref{1high} and the fact that all the momentum variables are comparable.
They key point is then to employ the following series of changes of variables
which begins with
$q_{i}-q_{i+1} \to q_{i+1}$ on $dq_{i+1}$,
$q_{i-1}-q_{i} \to q_{i}$ on $dq_{i}$,
 and ends with 
  $p-q_1 \to q_1$ on $dq_1$.
The end result, with Lemma \ref{boundK2},  is that
\begin{multline*}
 \left( \int_{\mathbb{R}^3} dq_1 ~ \left| k_{\approx}^\chi(p,q_1)  \right| 
\cdots
\int_{\mathbb{R}^3} dq_{i-1}~ \left| k_{\approx}^\chi(q_{i-2}, q_{i-1}) \right|
\right)
\\
\times
\left(
\int_{\mathbb{R}^3} dq_{i} 
\int_{\mathbb{R}^3} dq_{i+1}
\left|  e^{-c|q_{i-1}- q_{i}|}   e^{-c|q_i - q_{i+1}|} \right|^2
\right)^{1/2}
\\
\lesssim
 \int_{\mathbb{R}^3} dq_1 ~ e^{-c|q_1|}
\cdots
\int_{\mathbb{R}^3} dq_{i-1}~ e^{-c|q_{i-1} |}
\left(
\int_{\mathbb{R}^3} dq_{i} 
\int_{\mathbb{R}^3} dq_{i+1}
e^{-2c|q_i|}
e^{-2c|q_{i+1}|}
\right)^{1/2}
\\
\lesssim 1.
\end{multline*}
Collecting the estimates in this paragraph establishes the {\it claim}.

We gather all of the estimates for $L_i(f)$ to obtain
\begin{multline*}
w_\ell(p)
\left|
L_i(f)
\right|
\lesssim
\left(1+t \right)
\int_{0}^{t} ~ ds_{i+1} ~
\left(
1+t-s_{i+1} \right)^{-j+3/2} 
\\
\times \| \left\{  U(s_{i+1}) f_0 \right\}\|_{2,\ell+j-2+5/b} 
\\
\lesssim
 \| f_0 \|_{2,\ell+j-2+5/b+k} 
\left(1+t \right)
\int_{0}^{t} ~ ds_{i+1} ~
\left(
1+t-s_{i+1} \right)^{-j+3/2} 
\left(1+s_{i+1} \right)^{-k-1}. 
\end{multline*}
Above we have used the decay of the linear solutions to \eqref{rBlin} from Theorem \ref{decay2}; these solutions are represented by \eqref{Udef}.  Then for $j\ge k+1 +3/2$ and $k' = j-2+5/b+k+\ell'$ (for any $\ell' > 3/b$) we  use Proposition \ref{BasicDecay} to show that
\begin{gather*}
w_\ell(p)
\left|
L_i(f)
\right|
\lesssim
 \| f_0 \|_{2,\ell+k'-\ell'} 
\left(1+t \right)^{-k}
\lesssim
 \| f_0 \|_{\infty,\ell+k'} 
\left(1+t \right)^{-k}. 
\end{gather*}
The last inequality above follows as in \eqref{ELLestINEQ}.
This is the desired estimate for $L_i(f)(t)$ which holds for any $i\ge 1$ and $k\ge 0$.

It remains to estimate 
$
G_{i+1}(f),
$
$
D_{i+1}(f),
$
and
$
N_{i+1}(f)
$
for $i\ge 1$.
First
\begin{multline*}
w_\ell (p) \left|G_{i+1}(f)\right|
\lesssim
 \int_{\mathbb{R}^3}  dq_1  \left| k_{\approx}^\chi(p,q_1)  \right|
 {\bf 1}_{high} \int_{\kappa_1}^t  ~ds_1 ~ e^{-\nu(p) (t-s_1) }
 \\
 \times 
\int_{\mathbb{R}^3} dq_2  \left| k_{\approx}^\chi(q_1,q_2) \right|
\cdots
\int_{\mathbb{R}^3} dq_{i+1} \left|  k_{\approx}^\chi(q_{i},q_{i+1})  \right|
\\
\times 
\int_{s_1-\kappa_{2}}^{s_1} ~ ds_{2} ~ e^{-\nu(q_1) (s_1- s_2) }
\cdots
 \int_{s_i-\kappa_{i+1}}^{s_i} ~ ds_{i+1} ~ e^{-\nu(q_i) (s_i- s_{i+1}) }
 \\
\times w_\ell (q_{i+1})  e^{-\nu(q_{i+1}) s_{i+1} } f_0(y_{i+1}-\hat{q}_{i+1}s_{i+1},q_{i+1}).
\end{multline*}
We have used that all the momentum variables are comparable and the trick from \eqref{iteratedE} to conclude that the upper bound above is further bounded as
\begin{multline*}
\lesssim
 \int_{\mathbb{R}^3}  dq_1  \left| k_{\approx}^\chi(p,q_1)\right|
 {\bf 1}_{high} 
\int_{\mathbb{R}^3} dq_2  \left| k_{\approx}^\chi(q_1,q_2)\right|
\cdots
\int_{\mathbb{R}^3} dq_{i+1}  \left| k_{\approx}^\chi(q_{i},q_{i+1})   \right|
\\
\times 
 e^{-C \nu(p) t } \frac{\| f_0 \|_{\infty, \ell + k}}{w_k(q_{i+1})} \int_{\kappa_1}^t  ~ds_1 ~ 
\int_{s_1-\kappa_{2}}^{s_1} ~ ds_{2} ~ 
\cdots
 \int_{s_i-\kappa_{i+1}}^{s_i} ~ ds_{i+1}.
\end{multline*}
Notice that using the definition of $\kappa_i = \kappa/ 2^{i-1}$ we have
\begin{gather*}
 \int_{\kappa_1}^t  ~ds_1 ~ 
\int_{s_1-\kappa_{2}}^{s_1} ~ ds_{2} ~ 
\cdots
 \int_{s_i-\kappa_{i+1}}^{s_i} ~ ds_{i+1}
 =
\kappa^{i} 2^{-i(i+1)/2}  \int_{\kappa_1}^t  ~ds_1 ~ 
\\
\lesssim
 (1+t).
\end{gather*}
Furthermore, as in \eqref{subclaim} and Lemma \ref{boundK2} with \eqref{1high} we have
\begin{multline*}
\int_{\mathbb{R}^3}  dq_1  \left| k_{\approx}^\chi(p,q_1)  \right|
 {\bf 1}_{high} 
 \\
 \times
\int_{\mathbb{R}^3} dq_2  \left| k_{\approx}^\chi(q_1,q_2)  \right|
\cdots
\int_{\mathbb{R}^3} dq_{i+1} \left|  k_{\approx}^\chi(q_{i},q_{i+1})  \right|  \frac{1}{w_k(q_{i+1})}
\\
 \lesssim 
\frac{1}{M^{2\zeta}}
\frac{1}{w_{k+i+1}(p)}
 \lesssim 
\frac{1}{w_{k+2}(p)}.
\end{multline*}
Moreover,  for any $k\ge 0$, by \eqref{ElemCalc} as in \eqref{decayEstSoft22} we have 
$$
 e^{-C \nu(p) t } 
 \lesssim 
 w_{k+2}(p) (1+t)^{-k-2}.
$$
Collecting the last few estimates we obtain
\begin{gather*}
w_\ell (p) \left|G_{i+1}(f)\right|
  \lesssim
 (1+t)^{-k-1}\| f_0 \|_{\infty, \ell + k}.
\end{gather*}
This is the desired estimate for $G_{i+1}(f)$.

We will now study  $D_{i+1}(f)$, which satisfies the following general estimate
\begin{gather*}
w_\ell (p) \left|D_{i+1}(f)\right|
\lesssim
 \int_{\mathbb{R}^3}  dq_1  \left| k_{\approx}^\chi(p,q_1) \right|
\cdots
\int_{\mathbb{R}^3} dq_{i+1} \left| k_{\approx}^\chi(q_{i},q_{i+1})   \right|
\\
\times 
\int_{\kappa_1}^t  ds_1 ~ e^{-\nu(p) (t-s_1) }
\int_{s_1-\kappa_{2}}^{s_1} ~ ds_{2} ~ e^{-\nu(q_1) (s_1- s_2) }
\cdots
 \int_{s_i-\kappa_{i+1}}^{s_i}  ds_{i+1} ~ e^{-\nu(q_i) (s_i- s_{i+1}) }
 \\
\times 
 {\bf 1}_{high} 
w_\ell (p) 
 \int^{s_{i+1}}_{0}  ds_{i+2} ~ e^{-\nu(q_{i+1}) (s_{i+1}- s_{i+2}) }
 \left| K^{1-\chi}\left(\{U(s_{i+2}) f_0\}\right)(y_{i+2},q_{i+1}) \right|.
\end{gather*}
Since all the momentum variables are comparable, with Lemma \ref{boundKinfX},  we have
\begin{multline*}
\lesssim  
 \int_{\mathbb{R}^3}  dq_1   \left| k_{\approx}^\chi(p,q_1) \right|
 {\bf 1}_{high} 
\int_{\mathbb{R}^3} dq_2   \left| k_{\approx}^\chi(q_1,q_2) \right|
\cdots
\int_{\mathbb{R}^3} dq_{i+1}  \left|  k_{\approx}^\chi(q_{i},q_{i+1})  \right|
\\
\times 
\int_{\kappa_1}^t  ~ds_1 ~ e^{-\nu(p) (t-s_1) }
\int_{s_1-\kappa_{2}}^{s_1} ~ ds_{2} ~ e^{-\nu(q_1) (s_1- s_2) }
\cdots
 \int_{s_i-\kappa_{i+1}}^{s_i} ~ ds_{i+1} ~ e^{-\nu(q_i) (s_i- s_{i+1}) }
 \\
 \times e^{-c p_0} 
||| f |||_{k,\ell}
 \int^{s_{i+1}}_{0} ~ ds_{i+2} ~ e^{-\nu(q_{i+1}) (s_{i+1}- s_{i+2}) }
(1+s_{i+2})^{-k}.
\end{multline*}
For this term Lemma \ref{boundKinfX} would allow  a better estimate for the momentum weight on $||| f |||_{k,\ell}$.
As in the estimate for $G_{i+1}$ above, with \eqref{1high}, we have
\begin{multline*}
w_\ell (p) \left|D_{i+1}(f)\right|
\lesssim
\frac{1}{M^{2\zeta}}
e^{-c p_0} ~ ||| f |||_{k,\ell}
\\
\times 
\int_{\kappa_1}^t  ds_1  e^{-C\nu(p) (t-s_1) }
\int_{s_1-\kappa_{2}}^{s_1}  ds_{2}  e^{-C\nu(p) (s_1- s_2) }
\cdots
 \int_{s_i-\kappa_{i+1}}^{s_i}  ds_{i+1}  e^{-C\nu(p) (s_i- s_{i+1}) }
 \\
\times 
 \int^{s_{i+1}}_{0} ~ ds_{i+2} ~ e^{-C\nu(p)(s_{i+1}- s_{i+2}) }(1+s_{i+2})^{-k}.
\end{multline*}
We have again used the crucial fact that all the momentum variables are comparable.  
Since we have exponential decay, we can iterate the estimates from \eqref{ElemCalc} and Proposition \ref{BasicDecay} as
in \eqref{decayEstSoft22} to obtain
\begin{gather}
 \notag
\int_{\kappa_1}^t  ~ds_1 ~ e^{-C\nu(p) (t-s_1) }
\\
 \notag
\int_{s_1-\kappa_{2}}^{s_1} ~ ds_{2} ~ e^{-C\nu(p)(s_1- s_2) }
\cdots
 \int_{s_i-\kappa_{i+1}}^{s_i} ~ ds_{i+1} ~ e^{-C\nu(p) (s_i- s_{i+1}) }
 \\
 \label{iterateINT}
\times 
 \int^{s_{i+1}}_{0} ~ ds_{i+2} ~ e^{-C\nu(p) (s_{i+1}- s_{i+2}) }~(1+s_{i+2})^{-k}
 \\
  \notag
 \le
 \int_{0}^t  ~ds_1 ~ e^{-C\nu(p) (t-s_1) }
\int_{0}^{s_1} ~ ds_{2} ~ e^{-C\nu(p) (s_1- s_2) }
\cdots
 \int_{0}^{s_i} ~ ds_{i+1} ~ e^{-C\nu(p) (s_i- s_{i+1}) }
\\
 \notag
\times
 \int^{s_{i+1}}_{0} ~ ds_{i+2} ~ e^{-C\nu(p) (s_{i+1}- s_{i+2}) }~(1+s_{i+2})^{-k}
  \\
   \notag
\lesssim
\int_{0}^t  ~ds_1 ~ (1+(t-s_1))^{-k-1}
\prod_{j=1}^{i}\int_{0}^{s_{j}}  ~ds_{j+1} ~ (1+(s_{j}-s_{j+1}))^{-k-1} 
\\
\times
w_{(k+1)(i+2)}(p) 
 \int^{s_{i+1}}_{0} ~ ds_{i+2} ~  (1+(s_{i+1}-s_{i+2}))^{-k-1} ~(1+s_{i+2})^{-k}.
 \notag
\end{gather}
After iteratively applying Proposition \ref{BasicDecay} we obtain an upper bound  of
\begin{gather*}
\lesssim
w_{(k+1)(i+2)}(p)(1+t)^{-k}.
\end{gather*}
Plugging this into the previous estimate we have
\begin{multline*}
w_\ell (p) \left|D_{i+1}(f)\right|
\le  
\frac{C}{M^{2\zeta}}
~w_{(k+1)(i+2)}(p) e^{-c p_0} (1+t)^{-k} ~ ||| f |||_{k,\ell}
\\
\le  
\frac{C}{M^{2\zeta}}
~ (1+t)^{-k} ~ ||| f |||_{k,\ell}.
\end{multline*}
This is the desired estimate for $D_{i+1}(f)$ when $M$ is chosen sufficiently large.

The final term to estimate is $N_{i+1}(f)$.  In this case we have the upper bound
\begin{multline*}
w_\ell (p) \left|N_{i+1}(f)\right|
\lesssim
w_\ell (p)  \int_{\mathbb{R}^3}  dq_1 \left|  k_{\approx}^\chi(p,q_1)  \right|
 {\bf 1}_{high} 
  \int_{\kappa_1}^t  ~ds_1 ~ e^{-\nu(p) (t-s_1) }
 \\
\times 
\int_{\mathbb{R}^3} dq_2 \left|  k_{\approx}^\chi(q_1,q_2) \right|
\cdots
\int_{\mathbb{R}^3} dq_{i+1} \left|  k_{\approx}^\chi(q_{i},q_{i+1})   \right|
\\
\times 
\int_{s_1-\kappa_{2}}^{s_1} ~ ds_{2} ~ e^{-\nu(q_1) (s_1- s_2) }
\cdots
 \int_{s_i-\kappa_{i+1}}^{s_i} ~ ds_{i+1} ~ e^{-\nu(q_i) (s_i- s_{i+1}) }
 \\
\times 
 \int_{0}^{s_{i+1}}  ds_{i+2}  e^{-\nu(q_{i+1}) (s_{i+1} - s_{i+2}) }
\int_{\mathbb{R}^3} dq_{i+2} ~ \left|  k_{\neq}^\chi(q_{i+1},q_{i+2})   \right|
 \\
\times 
 \left| \{U(s_{i+2}) f_0\}(y_{i+2},q_{i+2}) \right|.
\end{multline*}
We first estimate the time integrals above.  Since the relevant momentum variables are all comparable, as in \eqref{iterateINT} and \eqref{decayEstSoft22} we have
\begin{multline*}
\int_{\kappa_1}^t  ds_1 ~ e^{-\nu(p) (t-s_1) }
\int_{s_1-\kappa_{2}}^{s_1}  ds_{2} ~ e^{-\nu(q_1) (s_1- s_2) }
\cdots
 \int_{s_i-\kappa_{i+1}}^{s_i}  ds_{i+1} ~ e^{-\nu(q_i) (s_i- s_{i+1}) }
 \\
\times 
 \int_{0}^{s_{i+1}}  ds_{i+2}  ~ e^{-\nu(q_{i+1}) (s_{i+1} - s_{i+2}) }(1+ s_{i+2} )^{-k}
 \\
  \lesssim 
w_{(k+1)(i+2)}(p)(1+t)^{-k}.
 \end{multline*}
Now for the momentum integrals,
from Lemma \ref{boundK2} and \eqref{1high}
we have
\begin{multline*}
w_\ell (p)w_{(k+1)(i+2)}(p)  \int_{\mathbb{R}^3}  dq_1  k_{\approx}^\chi(p,q_1)
\int_{\mathbb{R}^3} dq_2  k_{\approx}^\chi(q_1,q_2)
\cdots
\int_{\mathbb{R}^3} dq_{i+1}  k_{\approx}^\chi(q_{i},q_{i+1})  
 \\
\times 
 {\bf 1}_{high} 
\int_{\mathbb{R}^3} dq_{i+2}  k_{\ne}^\chi(q_{i+1},q_{i+2})
\\
\le
\frac{C}{M^{(i+2)\zeta}}
\int_{\mathbb{R}^3}  dq_1  ~ e^{-c|p-q_1|} ~
 {\bf 1}_{high}  ~
\int_{\mathbb{R}^3} dq_2  ~ e^{-c|q_1-q_2|}
\cdots
\int_{\mathbb{R}^3} dq_{i+1}  ~
e^{-c|q_{i}-q_{i+1}|}
 \\
\times 
\int_{\mathbb{R}^3} dq_{i+2} ~
e^{-c|q_{i+1}-q_{i+2}|} \left( {\bf 1}_{|q_{i+2}| \ge 2|q_{i+1}|}
+
{\bf 1}_{ |q_{i+1}|\ge 2|q_{i+2}|} \right)~ w_{\ell+k(i+2)}(q_{i+1})
\\
\le
\frac{C}{M^{(i+2)\zeta}}.
\end{multline*}
In the last step we have used the following estimate
$$
\left( {\bf 1}_{|q_{i+2}| \ge 2|q_{i+1}|}
+
{\bf 1}_{ |q_{i+1}|\ge 2|q_{i+2}|} \right)
e^{-\frac{c}{2}|q_{i+1}-q_{i+2}|} ~ w_{\ell+k(i+2)}(q_{i+1})
\le C.
$$
Indeed in either of the regions 
$
 |q_{i+1}|\ge 2|q_{i+2}|
$
or
$
|q_{i+2}| \ge 2|q_{i+1}|
$
we can use estimates such as those in \eqref{split1314} or \eqref{split2314} to directly establish this bound.

Now by collecting the estimates in this paragraph, we have shown
\begin{gather*}
w_\ell (p) \left|N_{i+1}(f)\right|
\le  
\frac{C}{M^{(i+2)\zeta}} ~ (1+t)^{-k} ~
||| f |||_{k,0}.
\end{gather*}
Since $k\ge 0$ is arbitrary, we conclude our estimate and our proposition after choosing $M$ sufficiently large in this last upper bound.
 \qed \\

This concludes our proof of rapid linear decay.

\section{Nonlinear $L^\infty$ rapid Decay}
\label{s:rapidD}

In Section \ref{s:n0bd}, we have proven the existence of mild solutions \eqref{nonLinProbSol} to the non-linear relativistic Boltzmann equation \eqref{rBoltz} with the soft potentials.
For $\ell >3/b$ and $k\in(1/2,1]$
we have shown in Thereom \ref{existenceandminordecay} that these solutions, $f=f(t,x,p)$, satisfy 
\begin{gather*}
\| f\|_{\infty,\ell}(t) \le C_{\ell,k}(1+t)^{-k} \| f_0 \|_{\infty,\ell+k}.
\end{gather*}
Then in Section \ref{s:linRdecay} we prove high order ``almost exponential'' decay for the linear semi-group as in \eqref{Udecay}.
From these estimates and the solution formula, \eqref{nonLinProbSol}, we can prove the following non-linear almost exponential decay.

\begin{theorem}
\label{rapidDECAYb}
Given any $\ell> 3/b$ and $k\ge 0$, there is a $k'=k'(k)\ge 0$ such that the solutions which were proven to exist in Thereom \ref{existenceandminordecay}  further satisfy 
$$
\| f \|_{\infty, \ell}(t) \le 
C_{\ell, k} (1+t)^{-k} \| f_0 \|_{\infty, \ell+k'}.
$$
\end{theorem}

\noindent {\it Proof of Theorem \ref{rapidDECAYb}.} 
We use an induction which allows one to continually improve the decay.
The main point is to bound the non-linear term, since we already know this kind of rapid decay for the linear part of \eqref{nonLinProbSol} from \eqref{Udecay} which follows from the crucial Proposition \ref{r1rapidDecay}.
In the first step we note that by 
Thereom \ref{existenceandminordecay}, Theorem \ref{rapidDECAYb} is true for $k\in(1/2,1]$.  Then 
 given any $j\ge 0$, from \eqref{Udecay} we have
 \begin{multline*}
w_{\ell}(p)
\left|
N[ f, f](t,x,p)
\right|
\lesssim
\int_0^t ~ ds ~
w_{\ell}(p)
\left|
 \{U(t-s) \Gamma [ f(s), f(s)]\}(x,p)
\right|
\\
\lesssim
\int_0^t ~ ds ~ (1+t-s)^{-j}
\left\|
  \Gamma [ f(s), f(s)]
\right\|_{\infty,\ell+j'}.
\end{multline*}
Above $j'\ge j$ is the number corresponding to $j$ in \eqref{Udecay}.  From Lemma \ref{nonlin0}  we have
\begin{gather*}
\lesssim
\int_0^t ~  ~ ds ~ (1+t-s)^{-j}
\left\|  f(s) \right\|_{\infty,\ell+j'-1}
\left\|  f(s) \right\|_{\infty,\ell+j'-1}.
\end{gather*}
Next we use the non-linear decay from Theorem \ref{existenceandminordecay} to see
\begin{multline*}
\lesssim
\left\|  f_0 \right\|_{\infty,\ell+j'+i-1}^2
\int_0^t ~  ~ ds ~ (1+t-s)^{-j}(1+s)^{-2i}
\\
\lesssim
(1+t)^{-\rho}\left\|  f_0 \right\|_{\infty,\ell+j'+i-1}^2.
\end{multline*}
The last estimate follows from Proposition \ref{BasicDecay}
with
 $\rho=\min\{j+2i -1,\min\{j,2i\} \}$.  
 In the above estimates we can choose $j\in (1,2]$ and then $i\in(1/2,1]$ such that $j=2i>1$.
Then we have shown Theorem \ref{rapidDECAYb} for $k\in (1,2]$ by choosing $\rho = j = k$ and 
$k' = \max\{j', j'+i-1\} = j'$.

Next suppose the theorem is correct for some $k> 2$, we will show that we may go beyond this $k$.  Indeed similar to the initial case we have
\begin{gather*}
w_{\ell}(p)
\left|
N[ f, f](t,x,p)
\right|
\le
C
(1+t)^{-\rho}\left\|  f_0 \right\|_{\infty,\ell+j'+i'-1}^2,
\end{gather*}
with $\rho=\min\{j+2i -1,\min\{j,2i\} \}=\min\{j,2i\}$.  
Above $j'$ corresponds to the power of the weight coming from decay level $j$ in \eqref{Udecay} and $i'$ corresponds to the power of the weight generated by decay level $i\in (0,k]$ in this Theorem \ref{rapidDECAYb}.

Choose  
$i\in(k/2, k]$  and 
$j=2i\in (k,2k]$. This is always possible.  Then we  have
$\rho = j$ so that we have proven Theorem \ref{rapidDECAYb}
for any $\tilde{k} \in (k,2k]$ and the corresponding $\tilde{k}' =\max\{ j', j'+i'-1\}$.
We conclude by induction.
\qed \\

%

\section*{Appendix:  Derivation of the Compact Operator}
\label{secAPP:HSr}

In this section, we give a complete exposition of the derivation of the Hilbert-Schmidt form \eqref{k2def} for the compact operator
from  \eqref{compactK}.
The linearized collision operator takes the form
\eqref{L}.
In that formulation we have the multiplication operator as in \eqref{nuDEF}.
The remaining ``compact'' part of the linearized operator is given by \eqref{compactK}
with $K= K_2-K_1$ and  in particular
\begin{eqnarray*}
K_2( f ) &\eqdef &
\frac{1}{p_0}\int_{\mathbb{R}^3} \frac{dq}{q_0}\int_{\mathbb{R}^3}\frac{dq^\prime}{q^\prime_0}\int_{\mathbb{R}^3}\frac{dp^\prime}{p^\prime_0}W(p, q | p^\prime, q^\prime)  
\sqrt{J(q)}
\left\{\sqrt{J(q^\prime)} f(p^\prime)\right\}
\\
&
&+ 
\frac{1}{p_0}\int_{\mathbb{R}^3} \frac{dq}{q_0}\int_{\mathbb{R}^3}\frac{dq^\prime}{q^\prime_0}\int_{\mathbb{R}^3}\frac{dp^\prime}{p^\prime_0}W(p, q | p^\prime, q^\prime)  
\sqrt{J(q)}
\left\{\sqrt{J(p^\prime)} f(q^\prime)\right\}.
\end{eqnarray*}
We are using the original notation from the top of this paper, which includes the delta functions.  We will outline in detail a procedure which is sketched in \cite[p.277]{MR635279} (see also \cite{MR0471665}), that allows a nice reduction of the term $K_2$ as in 
\eqref{k2def}.  In particular we give the exact form of the Lorentz transformation.
This reduction for the $K_1$ term can be reduced to the form \eqref{k1def} using much simpler methods than the ones we use below,
see e.g. \cite{MR635279,strainPHD,strainNEWT}.

We recall the definition of the transition rate, $W$, from the top of this paper.  We plug the definition of $W$ into $K_2$ above to obtain
\begin{multline*}
K_2( f ) \eqdef 
\frac{1}{p_0}\int_{\mathbb{R}^3} \frac{dq}{q_0}\int_{\mathbb{R}^3}\frac{dq^\prime}{q^\prime_0}\int_{\mathbb{R}^3}\frac{dp^\prime}{p^\prime_0}s\sigma(g, \theta) \delta^{(4)}(p^\mu+q^\mu-p^{\mu\prime}-q^{\mu\prime})
\\
\quad
\times
\sqrt{J(q)}
\left\{\sqrt{J(q^\prime)} f(p^\prime)+\sqrt{J(p^\prime)} f(q^\prime)\right\}.
\end{multline*}
We will first reduce this to a Hilbert-Schmidt form and second carry out the delta function integrations in the  kernel. 

In preparation, we write down some invariant quantities.  By (\ref{collisionalCONSERVATION}) and (\ref{sDEFINITION})
\begin{equation*}
\begin{split}
(p^{\mu \prime}-q^{\mu \prime})(p_\mu-q_{\mu})
&= 
2p^{\mu \prime} p_\mu +2 q^{\mu \prime}q_\mu
-
 (p^{\mu \prime}+q^{\mu \prime})(p_\mu+q_{\mu})
\\
&= 2p^{\mu \prime} p_\mu +2 q^{\mu \prime}q_\mu+s.
\end{split}
\end{equation*}
Further notice that (\ref{collisionalCONSERVATION}) implies
$$
(p^{\mu}-p^{\mu \prime}) (p_{\mu}-p^{ \prime}_{\mu}) 
= 
(q^{\mu \prime}-q^{\mu}) (q_{\mu}^{ \prime}-q_{\mu}).
$$
Expanding this we have
$$
-2-2p^{\mu} p_{\mu}^{\prime} = -2-2q^{\mu \prime}q_{\mu}.  
$$
We thus have another invariant $p^{\mu} p_{\mu}^{\prime} = q^{\mu}q_{\mu}^{\prime}$. Define $\bar{g}=g(p^{\mu},p^{\mu \prime})$ as in (\ref{gDEFINITION}).  We will always use $g$ without the bar to exclusively denote $g=g(p^{\mu},q^{\mu})$.

From \eqref{sDEFINITION} and \eqref{gDEFINITION} we know
$
s=g^2+4.
$
We may re-express $\theta$ from \eqref{angle} as
\begin{equation}
\cos\theta
= 
(p^\mu - q^\mu) (p_\mu^\prime -q_\mu^\prime)/g^2
=1-2\left(\frac{\bar{g}}{g}\right)^2.
\nonumber
\end{equation}
This follows from the invariant calculations in the previous paragraph and
$$
(p^\mu - q^\mu) (p_\mu^\prime -q_\mu^\prime)
=
g^2 + 4 +4p^{\mu} p_{\mu}^{\prime}
=
g^2 - 2\bar{g}^2.
$$
We further {\it claim} that
\begin{equation}
g^2=\bar{g}^2-\frac 12 (p^{\mu}+p^{\mu \prime}) (q_{\mu}+q_{\mu}^{ \prime}-p_{\mu}-p_{\mu}^{ \prime}).
\label{gbarIDENTc}
\end{equation}
Let $\bar{s}=s(p^{\mu},p^{\mu \prime})=\bar{g}^2+4$.  Then (\ref{gbarIDENTc}) is equivalent to 
\begin{eqnarray*}
g^2&=&\bar{g}^2-\frac{1}{2}\bar{s}-\frac 12 (p^{\mu}+p^{\mu \prime}) (q_{\mu}+q_{\mu}^{ \prime})
\\
&=&\frac 12\bar{g}^2 
-2-\frac 12 (p^{\mu}+p^{\mu \prime}) (q_{\mu}+q_{\mu}^{ \prime})
\\
&=&\frac 12\bar{g}^2 
+g^2+2p^{\mu}q_{\mu}-\frac 12 (p^{\mu}+p^{\mu \prime}) (q_{\mu}+q_{\mu}^{ \prime}).
\end{eqnarray*}
We thus prove (\ref{gbarIDENTc}) by  showing that 
$$
\frac{1}{2}\bar{g}^2 +2p^{\mu} q_{\mu}-\frac 12 (p^{\mu}+p^{\mu \prime})(q_{\mu}+q_{\mu}^{ \prime})=0.
$$
Expanding this expression we obtain
$$
-p^{\mu}p_{\mu}^{ \prime}-1+2p^{\mu} q_{\mu}-\frac 12 p^{\mu} q_{\mu}
-\frac 12 p^{\mu} q_{\mu}^{ \prime}-\frac 12 p^{\mu \prime} q_{\mu}-\frac{1}{2} p^{\mu \prime} q_{\mu}^{ \prime}.
$$
Notice further that
$
p^{\mu} q_{\mu}=p^{\mu \prime} q_{\mu}^{ \prime}
$ 
and 
$
p^{\mu} p_{\mu}^{ \prime}=q^{\mu} q_{\mu}^{ \prime}
$ 
as a result of (\ref{collisionalCONSERVATION}) and \eqref{sDEFINITION}.   We thus
obtain
$$
p^{\mu} q_{\mu}-1
-
\frac 12 p^{\mu} p_{\mu}^{ \prime}
-
\frac 12 q^{\mu} q_{\mu}^{ \prime}
-
\frac 12 p^{\mu} q_{\mu}^{ \prime}
-
\frac 12 p^{\mu \prime} q_{\mu},
$$
which by (\ref{collisionalCONSERVATION}) is
$$
p^{\mu} q_{\mu}-1
-\frac 12 (p^{\mu}+q^{\mu}) (p_{\mu}^{\prime}+q_{\mu}^{\prime})
=p^{\mu} q_{\mu}-1
+\frac 12 s=0.
$$
This establishes the claim (\ref{gbarIDENTc}).

Now we establish the Hilbert-Schmidt form.  First consider
$$
\frac{1}{p_0}\int_{\mathbb{R}^3} \frac{dq}{q_0}\int_{\mathbb{R}^3}\frac{dq^\prime}{q^\prime_0}\int_{\mathbb{R}^3}\frac{dp^\prime}{p^\prime_0}s \sigma(g, \theta) \delta^{(4)}(p^{\mu}+q^{\mu}-p^{\mu \prime}-q^{\mu \prime}) 
\sqrt{J(q)}
\sqrt{J(q^\prime)} f(p^\prime).
$$
Exchanging $q$ with $p^\prime$ the integral above is equal to
$$
\frac{1}{p_0}\int \frac{dq}{q_0}f(q)\left\{ \int  \frac{dq^\prime}{q^\prime_0}
\int \frac{dp^\prime}{p^\prime_0}
\bar{s} \sigma(\bar{g}, \theta) \delta^{(4)}(p^{\mu}+p^{\mu \prime}-q^{\mu}-q^{\mu \prime})
\sqrt{J(p^\prime)}
\sqrt{J(q^\prime)} \right\},
$$
where  $\theta$ is now defined by
\begin{equation}
\cos\theta= 
1-2\left(\frac{g}{\bar{g}}\right)^2,
\label{theta1}
\end{equation}
and from (\ref{gbarIDENTc}), with the new argument in the delta function above, we have
\begin{equation}
\bar{g}^2=g^2+\frac 12 (p^{\mu}+q^{\mu}) (p_{\mu}+q_{\mu} - p_{\mu}^{ \prime} - q_{\mu}^{ \prime}),
\label{gbarIDENT}
\end{equation}
and further $\bar{s}$ is defined by
$
\bar{s}=\bar{g}^2+4.
$
We do a similar calculation for the second term in $K_2f$, e.g. exchange $q$ with $q^\prime$ and then swap the $q^\prime$ and $p^\prime$ notation.  
The result is that we can define
\begin{equation}
k_2(p,q)
\eqdef
\frac{2}{p_0q_0}
\int  \frac{dq^\prime}{q^\prime_0}\int \frac{dp^\prime}{p^\prime_0}
\bar{s} \sigma(\bar{g}, \theta) \delta^{(4)}(p^{\mu}+p^{\mu \prime}-q^{\mu}-q^{\mu \prime})
\sqrt{J(p^\prime)}
\sqrt{J(q^\prime)}.
\label{k2KERNEL}
\end{equation}
We now write the Hilbert-Schmidt form
$
K_2 (f)
=
\int k_2(p,q) f(q) dq.
$
We will carry out the delta function integrations in $k_2(p,q)$ using a special Lorentz matrix.

We first translate \eqref{k2KERNEL} into an expression involving the total and relative momentum variables, 
$p^{\mu \prime}+q^{\mu \prime}$ and $p^{\mu \prime}-q^{\mu \prime}$ respectively.  Define $u$ by $u(r) = 0$ if $r<0$ and $u(r) =1$ if $r\ge 0$.  Let $\underline{g}=g(p^{\mu \prime},q^{\mu \prime})$ and $\underline{s}=s(p^{\mu \prime},q^{\mu \prime})$ .  We 
{\it claim} that 
\begin{equation}
\int_{\mathbb{R}^3}  \frac{dq^\prime}{q^\prime_0}\int_{\mathbb{R}^3} \frac{dp^\prime}{p^\prime_0}
G(p^{\mu},q^{\mu},p^{\mu \prime},q^{\mu \prime})
=
\frac{1}{16}\int_{\mathbb{R}^4\times\mathbb{R}^4}  d\Theta(p^{\mu \prime}, q^{\mu \prime})
G(p^{\mu},q^{\mu},p^{\mu \prime},q^{\mu \prime}),
\label{claimINT}
\end{equation}
where we are now integrating over the eight vector $(p^{\mu \prime}, q^{\mu \prime})$ and
$$
d\Theta(p^{\mu \prime}, q^{\mu \prime})= dp^{\mu \prime} dq^{\mu \prime} u(p_0^\prime+q_0^\prime)
u(\underline{s}-4)
\delta(\underline{s}-\underline{g}^2-4)
\delta((p^{\mu \prime}+q^{\mu \prime}) (p_{\mu}^{ \prime}-q_{\mu}^{ \prime})).
$$
To establish the claim, first notice that 
\begin{equation*}
\begin{split}
-(p^{\mu \prime}+q^{\mu \prime}) (p_{\mu}^{ \prime}-q_{\mu}^{ \prime})
&=
-p^{\mu \prime}p_{\mu}^{ \prime}+q^{\mu \prime}q_{\mu}^{ \prime}
\\
&=
(p_0^\prime)^2-|p^\prime|^2-(q_0^\prime)^2+|q^\prime|^2
=
A_p-A_q,
\end{split}
\end{equation*}
where now $p_0^\prime$ and $q_0^\prime$ are integration variables and we have defined
$$
A_p=(p_0^\prime)^2-(|p^\prime|^2+1),
\quad
A_q=(q_0^\prime)^2-(|q^\prime|^2+1).
$$
Integrating first over $dp^{\mu \prime}$, we see that alternatively
\begin{equation*}
\begin{split}
-(p^{\mu \prime}+q^{\mu \prime}) (p_{\mu}^{ \prime}-q_{\mu}^{ \prime})
&=
(p_0^\prime)^2-(|p^\prime|^2+1+A_q)
\\
&=
\left\{p_0^\prime-\sqrt{|p^\prime|^2+1+A_q}\right\}\left\{p_0^\prime+\sqrt{|p^\prime|^2+1+A_q}\right\}.
\end{split}
\end{equation*}
Furthermore, by \eqref{sDEFINITION} and \eqref{gDEFINITION} we have \begin{equation*}
\begin{split}
\underline{s}-\underline{g}^2-4
&=
-(p^{\mu \prime}+q^{\mu \prime}) (p_{\mu}^{ \prime}+q_{\mu}^{ \prime})
-
(p^{\mu \prime}-q^{\mu \prime}) (p_{\mu}^{ \prime}-q_{\mu}^{ \prime})-4
\\
&=
-2p^{\mu \prime} p_{\mu}^{ \prime}
-
2q^{\mu \prime} q_{\mu}^{ \prime}-4
\\
&
=
2A_p+2A_q.
\end{split}
\end{equation*}
Then similarly
\begin{equation*}
\begin{split}
\underline{s}-\underline{g}^2-4
&=
2(q_0^\prime)^2-2[|q^\prime|^2+1-A_p]
\\
&=
2\left\{q_0^\prime-\sqrt{|q^\prime|^2+1-A_p}\right\}\left\{q_0^\prime+\sqrt{|q^\prime|^2+1-A_p}\right\}.
\end{split}
\end{equation*}
Further note that 
$
p_0^\prime+q_0^\prime\ge 0
$
and
$
\underline{s}-4\ge 0
$
together imply 
$
p_0^\prime\ge 0
$
and 
$
q_0^\prime\ge 0.
$
With these expressions and standard  calculations we establish \eqref{claimINT}.

We thus conclude that
\begin{equation*}
k_2(p,q)
=
\frac{1}{p_0q_0}\frac{1}{8}
\int_{\mathbb{R}^4\times\mathbb{R}^4}  d\Theta(p^{\mu \prime}, q^{\mu \prime})
\bar{s} \sigma(\bar{g}, \theta) \delta^{(4)}(p^{\mu}+p^{\mu \prime}-q^{\mu}-q^{\mu \prime})
\sqrt{J(q^\prime)J(p^\prime)}.
\end{equation*}
Now apply the change of variables
$$
\bar{p}^{\mu}=p^{\mu \prime}+q^{\mu \prime}, 
\quad
\bar{q}^{\mu}=p^{\mu \prime}-q^{\mu \prime}.
$$
This transformation has Jacobian $=16$ and inverse tranformation as
$$
p^{\mu \prime} =\frac 12 \bar{p}^{\mu} +\frac 12 \bar{q}^{\mu}, 
\quad
q^{\mu \prime} =\frac 12 \bar{p}^{\mu} -\frac 12 \bar{q}^{\mu}.
$$
With this change of variable, for some $c'>0$, the integral becomes
\begin{equation*}
k_2(p,q)
=
\frac{c'}{ p_0q_0}
\int_{\mathbb{R}^4\times\mathbb{R}^4}  d\Theta(\bar{p}^{\mu}, \bar{q}^{\mu})
\bar{s} \sigma(\bar{g}, \theta) \delta^{(4)}(p^{\mu}-q^{\mu}+\bar{q}^{\mu})
\sqrt{J(\bar{p})},
\end{equation*}
with 
$
\sqrt{J(\bar{p})}
=
e^{- \bar{p}_0/2}
$
(ignoring constants).
Above the measure is now
$$
d\Theta(\bar{p}^{\mu}, \bar{q}^{\mu})= d\bar{p}^{\mu} ~ d\bar{q}^{\mu} ~ u(\bar{p}_0)
u(-\bar{p}^{\mu} \bar{p_{\mu}}-4)
\delta(-\bar{p}^{\mu} \bar{p_{\mu}}-\bar{q}^{\mu} \bar{q_{\mu}}-4)
\delta(\bar{p}^{\mu} \bar{q_{\mu}}).
$$
Also $\bar{g}\ge 0$ from (\ref{gbarIDENT}) is now given by
$$
\bar{g}^2=g^2+\frac 12 (p^{\mu}+q^{\mu}) (p_{\mu}+q_{\mu}-\bar{p}_{\mu}).
$$
And $\theta$ and $\bar{s}$ can be defined through the new $\bar{g}$ with (\ref{theta1}).

We next carry out the delta function argument for $\delta^{(4)}(p^{\mu}-q^{\mu}+\bar{q}^{\mu})$ to obtain
\begin{equation*}
k_2(p,q)
=
\frac{c'}{p_0q_0}
\int_{\mathbb{R}^4}  d\Theta(\bar{p}^{\mu})
\bar{s} \sigma(\bar{g}, \theta)
e^{- \bar{p}_0/2},
\quad 
\exists c' > 0.
\end{equation*}
where the measure is 
$$
d\Theta(\bar{p}^{\mu})= d\bar{p}^{\mu} ~ u(\bar{p}_0) ~ 
u(-\bar{p}^{\mu} \bar{p}_{\mu}-4) ~
\delta(-\bar{p}^{\mu}\bar{p}_{\mu}-g^2-4)
\delta(\bar{p}^{\mu} (q_{\mu}-p_{\mu})).
$$
Since
$
s = g^2 + 4
$
 we have
\begin{eqnarray*}
u(\bar{p}_0)
\delta(-\bar{p}^{\mu}\bar{p}_{\mu} -g^2-4)
&=&
u(\bar{p}_0)
\delta(-\bar{p}^{\mu}\bar{p}_{\mu} -s)
\\
&=&
u(\bar{p}_0)
\delta((\bar{p}_0)^2-|\bar{p}|^2-s)
\\
&=&
\frac{\delta(\bar{p}_0-\sqrt{|\bar{p}|^2+s})}{2\sqrt{|\bar{p}|^2+s}}.
\end{eqnarray*}
We then carry out one integration using the delta function to get
\begin{equation*}
k_2(p,q)
=
\frac{c^\prime }{p_0q_0}
\int_{\mathbb{R}^3}  \frac{d\bar{p}}{ \bar{p}_0}u(-\bar{p}^{\mu}\bar{p}_{\mu} -4)
\delta(\bar{p}^{\mu} (q_{\mu}-p_{\mu}))
\bar{s} \sigma(\bar{g}, \theta)
e^{- \bar{p}_0 /2},
\end{equation*}
with
$
\bar{p}_0\eqdef \sqrt{|\bar{p}|^2+s}.
$
Using
$
s = g^2 + 4
$
we have
$$
-\bar{p}^{\mu}\bar{p}_{\mu} -4=s-4=g^2\ge 0.
$$
So always $u(-\bar{p}^{\mu}\bar{p}_{\mu} -4)=1$
and the integral reduces to
\begin{equation*}
k_2(p,q)
=
\frac{c^\prime }{p_0q_0}
\int_{\mathbb{R}^3}  \frac{d\bar{p}}{ \bar{p}_0}
\delta(\bar{p}^{\mu} (q_{\mu}-p_{\mu}))
\bar{s} \sigma(\bar{g}, \theta)
e^{- \bar{p}^{\mu} \bar{U}_{\mu}/2},
\end{equation*}
where 
$\bar{U}=(-1,0,0,0)^t$
and
$
e^{- \bar{p}_0/2}
=
e^{- \bar{p}^{\mu} \bar{U}_{\mu}/2}.
$

We finish off our reduction by moving to a new Lorentz frame.  
We consider the Lorentz Transformation $\Lambda$ which maps 
$$
-\Lambda^{\mu\nu} (p_\mu + q_\mu) = (\sqrt{s},0,0,0), 
\quad 
-\Lambda^{\mu\nu} (p_\mu - q_\mu) = (0,0,0,g).
$$
Our notation here for raising and lowering indicies is standard in the sense that 
$
p_\mu = \sum_{\nu = 0}^3 p^\nu g^{\nu\mu}
$
where $\left( g^{\nu\mu} \right) = \text{diag}(-1 ~ 1 ~ 1  ~ 1)$.  
We recall that $p^\mu = (-p_0, p)$.
Then we use the Einstein summation convention  as 
$
\Lambda^{\mu\nu} p_\mu
= \sum_{\mu = 0}^3 \Lambda^{\mu\nu} p_\mu.
$
From this information, we have derived in \cite{strainPHD} and exposited in \cite{strainNEWT} that
$$
\Lambda
=
\left(\Lambda^{\mu\nu} \right)
=
\left(
\begin{array}{cccc}
\frac{p_0+q_0}{\sqrt{s}} & -\frac{p_1+q_1}{\sqrt{s}} & -\frac{p_2+q_2}{\sqrt{s}} &-\frac{p_3+q_3}{\sqrt{s}}
\\
\Lambda^{01} & \Lambda^{11} & \Lambda^{21} & \Lambda^{31}
\\
0 & \frac{(p\times q)_1}{|p\times q|}& \frac{(p\times q)_2}{|p\times q|}& \frac{(p\times q)_3}{|p\times q|}
\\
\frac{p_0-q_0}{g} & -\frac{p_1-q_1}{g} & -\frac{p_2-q_2}{g} &-\frac{p_3-q_3}{g}
\end{array}
\right).
$$
With second row given by
$$
\Lambda^{01}= \frac{2|p\times q|}{g\sqrt{s}},
$$
and
$$
\Lambda^{i1}=
\frac{2\left( p_i\left\{p_0+q_0 p^\mu q_\mu\right\}+q_i\left\{q_0+p_0 p^\mu q_\mu \right\}\right)}{g\sqrt{s}|p\times q|} 
\quad (i=1,2,3).
$$
We have a complete description of this Lorentz transformation in terms of $p, q$.

Define $U^\mu=\Lambda^{\nu\mu}\bar{U}_\nu$, notice that
$$
U^\mu
=\left(\frac{p_0+q_0}{\sqrt{s}},\frac{2|p\times q|}{g\sqrt{s}},~0,~\frac{p_0-q_0}{g} \right).
$$
Then
$$
 \int \frac{d\bar{p}}{\bar{p}_0} \bar{s} \sigma(\bar{g}, \theta) e^{- \bar{p}_0/2}
\delta(\bar{p}^{\mu} (q_{\mu}-p_{\mu})) 
=
 \int \frac{d\bar{p}}{\bar{p}_0} \bar{s}_\Lambda \sigma(\bar{g}_\Lambda, \theta_\Lambda) 
e^{-\frac {1}2\bar{p}^{\mu} U_{\mu}} \delta(\bar{p}^{\mu} \Lambda_{~\mu}^{\nu~}(q_\nu-p_\nu)).
$$
where $\bar{g}_\Lambda, \bar{s}_\Lambda\ge 0$ are now given by
\begin{equation}
\begin{split}
\bar{g}_\Lambda^2
&=g^2+\frac 12 \Lambda^{\nu\mu}(p_\nu+q_\nu) \{\Lambda_{~\mu}^{\gamma~}(p_\gamma+q_\gamma)-\bar{p}_{\mu}\}
=g^2+\frac 12 \sqrt{s} \{\bar{p}_0-\sqrt{s}\},
\\
\bar{s}_\Lambda
&=4+\bar{g}_\Lambda^2,
\\
\cos\theta_\Lambda
&=
1-2\left(\frac{g}{\bar{g}_\Lambda}\right)^2.
\label{angle2}
\end{split}
\end{equation}
The equality of the two integrals holds because 
$d\bar{p}/\bar{p}_0$ is a Lorentz invariant.

We work with the integral on the right hand side above. Now
$$
\bar{p}^{\mu} \Lambda_{~\mu}^{\nu~}(q_\nu-p_\nu) =\bar{p}_3g.
$$
We switch to polar coordinates in the form
$$
d\bar{p}=|\bar{p}|^2d|\bar{p}|  \sin\psi  d\psi d\varphi, ~~~ 
\bar{p}\equiv |\bar{p}|(\sin\psi\cos\varphi, \sin\psi\sin\varphi, \cos\psi).
$$
Then we can write $k_2(p,q)$ as
\begin{equation*}
\frac{c^\prime}{p_0q_0}
 \int_0^{2\pi}d\varphi\int_0^{\pi}\sin\psi d\psi
\int_0^{\infty} \frac{|\bar{p}|^2  d|\bar{p}| }{\bar{p}_0}
\bar{s}_\Lambda ~
 \sigma(\bar{g}_\Lambda, \theta_\Lambda)
 ~ e^{- \bar{p}^{\mu} U_{\mu} /2}  ~ \delta(|\bar{p}|g\cos\psi).
\end{equation*}
We evaluate the last delta function at $\psi=\pi/2$ to write $k_2(p,q)$ as
\begin{equation}
\frac{c^\prime}{g p_0q_0}
\int_0^{2\pi}d\varphi\int_0^{\infty} \frac{|\bar{p}|  d|\bar{p}| }{\bar{p}_0}
\bar{s}_\Lambda\sigma(\bar{g}_\Lambda, \theta_\Lambda) 
e^{- \bar{p}_0\frac{p_0+q_0}{2\sqrt{s}} }e^{\frac{|p\times q|}{g\sqrt{s}}|\bar{p}|\cos\varphi}.
\label{lastINT}
\end{equation}
This is already a useful reduced form for $k_2(p,q)$.

We recall,
for
$
I_0\left(\frac{|p\times q|}{g\sqrt{s}}|\bar{p}|\right),
$
 the modified Bessel function of index zero given in \eqref{bessel0}. 
We further re-label the integration as $|\bar{p}| = y$.
Notice that (\ref{angle2}) implies
\begin{equation*}
g=\bar{g}_\Lambda \sqrt{\frac{1-\cos\theta_\Lambda}{2}}=\bar{g}_\Lambda\sin\frac{\theta_\Lambda}{2},
\end{equation*}
with
\begin{equation*}
\sin\frac{\theta_\Lambda}{2}
=
\frac{g}{\bar{g}_\Lambda}
=
\frac{g}{\sqrt{g^2- s/2 + \frac{\sqrt{s}}{2}\sqrt{y^2 + s}}}
=
\frac{\sqrt{2}g}{\sqrt{g^2- 4 + s\sqrt{y^2/s + 1}}}.
\end{equation*}
We may rewrite \eqref{lastINT} as
\begin{equation*}
k_2(p,q)
=
\frac{c^\prime}{g p_0q_0}
\int_0^{\infty} \frac{y  dy }{\sqrt{y^2+s}}
\bar{s}_\Lambda
\sigma\left(\frac{g}{\sin\frac{\theta_\Lambda}{2}}, \theta_\Lambda \right) 
e^{- \frac{p_0+q_0}{2\sqrt{s}} \sqrt{y^2+s}}
I_0\left(\frac{|p\times q|}{g\sqrt{s}}y\right).
\end{equation*}
From \eqref{angle2} we have
$$
\bar{s}_\Lambda
=
4+g^2+\frac 12 \sqrt{s} \{\sqrt{y^2+s}-\sqrt{s}\}
=
\frac 12 s+\frac 12 s \sqrt{y^2/s+1}.
$$
We apply the change of variables $y \to y/\sqrt{s}$ to obtain that
$k_2(p,q)$
is given by
\begin{equation*}
\frac{c^\prime s^{3/2}}{g p_0q_0}
\int_0^{\infty} \frac{y  \left(1+ \sqrt{y^2+1}\right)dy }{\sqrt{y^2+1}}
\sigma\left(\frac{g}{\sin\frac{\psi}{2}}, \psi\right) 
e^{- \frac{p_0+q_0}{2} \sqrt{y^2+1}}
I_0\left(\frac{|p\times q|}{g}y\right).
\end{equation*}
where $\sin\frac{\psi}{2}$ is given by
\eqref{sinPSI}.
This is the expression from \eqref{k2def} once we incorporate the cut-off function \eqref{cut} which is insignificant for the purposes of this calculation.

Significant simplifications can be performed in the case of the ``hard ball'' cross section where $\sigma = constant$.
The relevant integral  is then a Laplace transform and a known integral, which can be calculated exactly via a taylor expansion \cite[p.134]{MR0350075}.   For instance, it is well known that for any $R>r\ge 0$ 
\begin{equation*}
\begin{split}
\int_0^\infty \frac{e^{-R\sqrt{1+y^2}}yI_0(ry)}{\sqrt{1+y^2}} dy
&=\frac{e^{-\sqrt{R^2-r^2}}}{\sqrt{R^2-r^2}},
\\
\int_0^\infty e^{-R\sqrt{1+y^2}}yI_0(ry) dy
&=\frac{R}{R^2-r^2}\left\{1+\frac{1}{\sqrt{R^2-r^2}}\right\}
e^{-\sqrt{R^2-r^2}}.
\end{split}
\end{equation*}
See for instance \cite{MR950173}, \cite[p.383]{MR1407750}, or  \cite[p.322]{MR1211782}. 

Using these formulas we can express the integral as
$$
k_2(p,q)=
\frac {c^\prime s^{3/2}}{gp_0q_0}
\tilde{U}_1(p,q)
 \exp\left( -\tilde{U}_2(p,q)\right),
$$
where $\tilde{U}_2(p,q)\eqdef \sqrt{\{(p_0+q_0)/2\}^2-(|p\times q|/g)^2}$ and
$$
\tilde{U}_1(p,q)\eqdef \left(1+\frac{p_0+q_0}{2}\left(\tilde{U}_2(p,q)\right)^{-1}+\frac{p_0+q_0}{2}\left(\tilde{U}_2(p,q)\right)^{-2}\right)
\left(\tilde{U}_2(p,q)\right)^{-1}.
$$
  Further, 
$$
\tilde{U}_2(p,q)  = \frac{ \sqrt{s}}{2g}|p-q|=|p-q|\sqrt{\frac{g^2+4}{4g^2}}.
$$
Therefore,
$
\tilde{U}_2(p,q)\ge \frac{1}{2}|p-q|+1.
$
This completes our discussion of the Hilbert-Schmidt form for the linearized collision operator.

\subsection*{Acknowledgments} 
 The author gratefully thanks the anonymous referees for their lengthy detailed comments which helped to substantially improve the presentation of this research paper.


\bibliographystyle{siam}
\bibliography{soft}


%

%

\end{document}